\documentclass[11pt]{article}
\usepackage{amssymb}
\usepackage{mathrsfs}
\usepackage{latexsym,amsmath,amssymb,amsfonts,epsfig,graphicx,cite,psfrag}
\usepackage{eepic,color,colordvi,amscd}
\usepackage{ebezier}
\usepackage{verbatim}

\newtheorem{theo}{Theorem}[section]

\newtheorem{ques}[theo]{Question}
\newtheorem{lem}[theo]{Lemma}

\newtheorem{conj}[theo]{Conjecture}

\def\qed{\hfill \rule{4pt}{7pt}}

\def\pf{\noindent {\it Proof. }}

\begin{document}

\title{The  Kelmans-Seymour conjecture IV: a proof\footnote{Partially supported by NSF grants DMS-1265564 and DMS-1600738}}

\author{
Dawei He\footnote{dhe9@math.gatech.edu; this work was started when DH was a student at East China Normal University}, 
Yan Wang\footnote{yanwang@gatech.edu}, Xingxing Yu\footnote{yu@math.gatech.edu; 
Part of the work was done while XY was visiting East China Normal University and received
partial support from ECNU under the Science and Technology Commission of Shanghai Municipality (STCSM) grant No. 13dz2260400}\\School of Mathematics \\
Georgia Institute of Technology\\
Atlanta, GA 30332-0160, USA }

\date{December 20, 2016}

\maketitle

\begin{abstract} 
A well known theorem of Kuratowski in 1932 states that a graph is planar if, and only if, it does not
contain a subdivision of $K_5$ or $K_{3,3}$. Wagner proved in 1937 that if a
graph other than $K_5$ does not contain any subdivision 
of $K_{3,3}$ then it is planar or it admits a cut of size at most 2.   
Kelmans and, independently, Seymour conjectured in the 1970s that if a graph does
not contain any subdivision of $K_5$ then it is planar or it admits a cut of size at most 4. 
In this paper, we give a proof of the Kelmans-Seymour conjecture. We
also discuss several related results and problems.

\bigskip
AMS Subject Classification: 05C10, 05C40, 05C83

Keywords: $K_5$-subdivision, independent paths, separation, connectivity, discharging, contraction.  
\end{abstract}


\newpage
\section{Introduction}
For a graph $G$, we use $TG$ to denote a subdivision of $G$, and the vertices in $TG$ that correspond to the vertices of $G$ are said to be its {\it branch} vertices. 
Thus, $TK_5$ denotes a subdivision of $K_5$, and  the vertices in a
$TK_5$ of degree four are its branch vertices. 

The well known result of Kuratowski \cite{Ku30} states that a graph is
planar if, and only if, it does not contain $TK_5$ or $TK_{3,3}$. A simple application of Euler's formula
for planar graphs shows that, for $n\ge 3$,  if an $n$-vertex graph has at least $3n-5$ edges then it must be nonplanar and, hence, contains $TK_5$ or $TK_{3,3}$. 
Dirac \cite{Di64} conjectured that  for $n\ge 3$,  if an $n$-vertex graph has at least $3n-5$ edges then it must contain $TK_5$. 
This conjecture was also reported by Erd\H{o}s and Hajnal \cite{EH64}. 
 K\'{e}zdy and McGuiness \cite{KM91} showed that a minimal
 counterexample to Dirac's conjecture must be 5-connected and contains
 $K_4^-$ (obtained from the complete graph $K_4$ by deleting an edge).  After some partial results in \cite{Sk68, Th74,Th96, Th97}, 
Dirac's conjecture was proved by Mader \cite{Ma98}, where 
he also showed that every 5-connected $n$-vertex graph with at least $3n-6$ edges contains $TK_5$ or $K_4^-$.

Seymour \cite{Se77} (also see \cite{Ma98, Th97}) and, independently, Kelmans \cite{Ke79} conjectured that every $5$-connected nonplanar graph contains $TK_5$. 
Thus, the Kelmans-Seymour conjecture implies Mader's theorem. This conjecture is also related to several interesting problems, which we will discuss in Section 7. 

The authors  \cite{HWY15I, HWY15II, HWY15III} produced lemmas needed for proving  this Kelmans-Seymour conjecture, and 
we  are now ready to prove it in this paper.

\begin{theo}
\label{main}
Every $5$-connected non-planar graph contains  $TK_5$. 
\end{theo}

The starting point of our work is  the following result of Ma and Yu \cite{MY10, MY13}:
Every $5$-connected nonplanar graph containing $K_4^-$ has a
$TK_5$. This result, combined with the result of  K\'{e}zdy and McGuiness \cite{KM91} on minimal
 counterexamples to Dirac's conjecture, gives an alternative proof
 of Mader's theorem. 
Also using this result,  Aigner-Horev \cite{AH12} proved that every $5$-connected nonplanar apex graph contains $TK_5$. A simpler proof of Aigner-Horev's result 
using discharging argument was obtained by Ma, Thomas and Yu, and, independently, by Kawarabayashi, see \cite{KMY15}.

We now briefly describe the process for proving Theorem~\ref{main}. For a more detailed version, we recommend the
reader to read Section 6 first, which  should also give motivation
to some of the  technical lemmas listed in Sections 2, 3, 4 and 5.

 Suppose $G$ is a 5-connected non-planar graph not containing $K_4^-$.
We fix a vertex $v\in V(G)$, and let $M$ be a maximal connected
subgraph of $G$ such that $v\in V(M)$, $G/M$ (the graph obtained from $G$ by contracting $M$) is 
nonplanar, $G/M$ contains no $K_4^-$, and 
$G/M$ is 5-connected (i.e., $M$ is contractible).
Note that $V(M)=\{v\}$ is possible. 
Let $x$ denote the vertex of $H:=G/M$ resulting from the contraction of $M$.  Then, for each subgraph $T$ of $H$ with $v\in V(T)$ and $T\cong K_2$ or $T\cong K_3$, 
$H/T$ is planar, or $H/T$ contains $K_4^-$, or $H/T$ is not
5-connected. If, for some $T$, $H/T$ is planar or contains
$K_4^-$ then we can find a $TK_5$ in $G$ using results  from \cite{HWY15I,HWY15II,HWY15III}.
Thus, in this paper, our main work is to deal with the final case: for any $T\subseteq H$ with $x\in V(T)$ and $T\cong K_2$ or $T\cong K_3$, 
$H/T$ is nonplanar, $H/T$ contains no $K_4^-$, and $H/T$ is not 5-connected. In this case,  there exists
$S_T\subseteq V(H)$ such that $V(T)\subseteq S_T$, $|S_T|=5$ or
$|S_T|=6$,  and $H-S_T$ is not connected. We will be using such cuts
to divide the graph into smaller parts and use them to find a special 
$TK_5$ in $H$.  The reason to also include the case $T\cong K_3$ is to avoid the situation 
when $T\cong  K_2$, $|S_T|=5$, and $H-S_T$ has exactly two components, 
one of which is trivial. This does not cause problem when $T\cong
K_3$,  as the graph $H$ would then contain $K_4^-$, and we could use
results from \cite{HWY15I, HWY15II, HWY15III}.

We will need a number of results from \cite{HWY15I,HWY15II,HWY15III}, which are given in Section 2. 
In Section 3, we derive a simplified version of a result on disjoint paths from \cite{Yu03I, Yu03II, Yu03III}, which will be 
used several times in  Section 4. 
For each subgraph $T$ of $H$ with $v\in V(T)$ and $T\cong K_2$ or
$T\cong K_3$, we will associate to it a quadruple $(T,S_T,A,B)$,
where, roughly, $A\cap B=\emptyset$, $H-S_T=A\cup B$, and $H$ has no
edge between $A$ and $B$. (A precise definition of a quadruple is given in Section 4.)
In Section 4, we prove some basic properties of quadruples, and take
care of two special cases involving  quadruples (using disjoint paths results from Section 3). 
In Section 5, we take care of other cases involving quadruples. 
We complete the proof of Theorem~\ref{main} in Section 6, and discuss several related problems in Section 7. 

\medskip

We end this section with some notation and terminology. Let $G$ be a graph. 
By $S\subseteq G$ we mean that $S$ is a subgraph of $G$. 
We may   view $S\subseteq V(G)$ as a subgraph of $G$ with vertex set
$S$ and no edges. 
For $S\subseteq G$, we use $G[S]$ to denote the subgraph of $G$ induced by $V(S)$. 
For any $x\in V(G)$ we use $N_G(x)$ to denote the 
neighborhood of $x$ in $G$, and for $S\subseteq G$ let $N_G(S)=\{x\in V(G)\setminus V(S): N_G(x)\cap V(S)\ne\emptyset\}$.
When understood, the reference to $G$ may be dropped. For $S\subseteq E(G)$, $G-S$ denotes the graph obtained from $G$ 
by deleting all edges in $S$; 
and for $K,L\subseteq G$, $K-L$ denotes the graph obtained from $K$ by deleting $V(K\cap L)$ and all edges of $K$ 
incident with $V(K\cap L)$.

A {\it separation} in a graph $G$ consists of a pair of subgraphs $G_1, G_2$ of $G$, denoted as $(G_1,G_2)$, such that
$E(G_1) \cup E(G_2)=E(G)$, $E(G_1\cap G_2)=\emptyset$, and $E(G_1)\cup (V(G_1)\setminus V(G_2)\ne \emptyset\ne E(G_2)\cup (V(G_2)\setminus V(G_1))$. 
The {\it order} of this separation is $|V(G_1)\cap V(G_2)|$, and $(G_1,G_2)$ is said to be a {\it $k$-separation} if
its order is $k$. Thus, a set $S\subseteq V(G)$ is a {\it $k$-cut} (or
a {\it cut} of size $k$)
in $G$, where $k$ is a positive integer,  if $|S|=k$ and $G$ has a
separation $(G_1,G_2)$ such that $V(G_1)\cap V(G_2)=S$ and $V(G_1-S)\ne
\emptyset\ne V(G_2-S)$. If $v\in V(G)$ and $\{v\}$ is a cut
of $G$, then $v$ is said to be a {\it cut vertex} of $G$.
For $A\subseteq V(G)$ and for a positive integer $k$, we say that $G$ is {\it $(k,A)$-connected} if, for any cut $S$ with $|S|<k$, every 
component of $G-S$ contains a vertex from $A$.

Given a path $P$ in a graph and $x,y \in V(P)$, $xPy$ denotes the
subpath of $P$ between $x$ and $y$ (inclusive). 
The {\it ends} of
the path $P$ are the vertices of the minimum degree in $P$, and all
other vertices of $P$ (if any) are its {\it
  internal} vertices. A path $P$ with ends $u$ and $v$ (or an $u$-$v$ path) is also said to be {\it
from $u$ to $v$} or {\it between $u$ and $v$}.
A collection of paths are said to be {\it independent} if no vertex
of any path in this collection is an internal vertex of any other path in the collection.

Let $G$ be a graph. 
Let $K\subseteq G$, $S\subseteq V(G)$, and $T$ a collection of 2-element subsets of $V(K)\cup S$. Then $K+(S\cup T)$ denotes the graph
with vertex set $V(K)\cup S$ and edge set $E(K)\cup T$, and if
$T=\{\{x,y\}\}$  we write $K+xy$ instead of $K+\{\{x,y\}\}$.

 For any positive integer $k$, let $[k]:=\{1,\ldots, k\}$.  A {\it 3-planar graph} $(G,{\cal A})$ consists
of a graph $G$ and  a set ${\cal A}=\{A_1,\ldots,
A_k\}$ of pairwise disjoint subsets of $V(G)$ (possibly ${\cal
  A}=\emptyset$ when $k=0$) such that
\vspace*{-1ex}
\begin{itemize} 
\item [(a)] for distinct $i, j\in [k]$, $N(A_i)\cap
A_j=\emptyset$, 
\vspace*{-1ex}
\item [(b)] for $i\in [k]$,  $|N(A_i)|\le 3$, and 
\vspace*{-1ex}
\item [(c)] if $p(G,{\cal A})$ denotes the 
graph obtained from $G$ by (for each $i$) deleting $A_i$ and
adding  edges joining every pair of distinct vertices in  $N(A_i)$,
then $p(G,{\cal A})$ may be drawn in a closed disc $D$ with no pair of edges
crossing such that, for each $A_i$ with $|N(A_i)|=3$, $N(A_i)$ induces a
facial triangle in $p(G,{\cal A})$. 
\end{itemize}
\vspace*{-1ex}
\noindent If, in addition,  $b_1, \ldots, b_n$ are  vertices of
$G$ such that $b_i\notin A_j$ 
for any $i\in [n]$ and  $j\in [k]$ and $b_1,$ $\ldots,b_n$ occur on the boundary of
the disc $D$ in that cyclic order,
then we say that $(G,{\cal A}, b_1,\ldots,b_n)$ is
{\it 3-planar}. If there is no need to specify ${\cal A}$, we will simply say 
that $(G,b_1,\ldots,b_n)$ is 3-planar.  If there is no need to specify the order of $b_1,\ldots, b_n$ then we simply 
say that  $(G,\{b_1,\ldots,b_n\})$ is 3-planar. When ${\cal A}=\emptyset$, we say that $(G, b_1,\ldots,b_n)$ and 
$(G,\{b_1,\ldots,b_n\})$ are planar.

\section{Previous results}

In this section, we list a number of previous results which we will use  as lemmas in our proof of Theorem~\ref{main}. 
We begin with the main result of \cite{MY10, MY13}.
\begin{lem} \label{K4-} 
Every $5$-connected nonplanar graph containing $K_4^-$ has \textsl{a} $TK_5$. 
\end{lem}

We also need the main result of \cite{HWY15II} to take care of the case when the vertex $x$ in $H=G/M$ (see Section 1) is a degree 2 vertex in a $K_4^-$ (which is 
$y_2$ in the lemma below).  

\begin{lem}\label{y_2}
Let $G$ be $\textsl{a}$ $5$-connected nonplanar graph and $\{x_1 , x_2 , y_1 , y_2\}\subseteq V(G)$ 
such that $G[\{x_1 , x_2 , y_1 , y_2\}]\cong K_4^-$ with $y_1y_2\notin E(G)$. 
Then one of the following holds:
\begin{itemize}
\item[$(i)$] $G$ contains $\textsl{a}$ $TK_5$ in which $y_2$ is not $\textsl{a}$ branch vertex.
\item[$(ii)$] $G-y_2$ contains $K_4^-$.
\item[$(iii)$] $G$ has $\textsl{a}$ $5$-separation $(G_1,G_2)$ such that $V(G_1\cap G_2)=\{y_2, a_1,a_2,a_3,a_4\}$, and  $G_2$ is 
the graph obtained from the edge-disjoint union of the $8$-cycle $a_1b_1a_2b_2a_3b_3a_4b_4\allowbreak a_1$ 
and the $4$-cycle $b_1b_2b_3b_4b_1$ by adding $y_2$ and the edges $y_2b_i$ for $i\in [4]$.
\item[$(iv)$] For $w_1,w_2,w_3\in N(y_2)-\{x_1,x_2\}$, $G-\{y_2v: v\notin \{w_1,w_2,w_3,x_1,x_2\}\}$ contains  $TK_5$.
\end{itemize}
\end{lem}

To deal with conclusion $(iii)$ of Lemma~\ref{y_2}, we need Proposition 1.3 from \cite{HWY15I} in which $a$ plays the role of $y_2$ in Lemma~\ref{y_2}. 

\begin{lem}\label{8cycle}
Let $G$ be \textsl{a} 5-connected nonplanar graph, $(G_1,G_2)$ \textsl{a} 5-separation in $G$, $V(G_1 \cap G_2)=\{a,a_1,a_2,a_3,a_4\}$ such that  
$G_2$ is the graph obtained from the edge-disjoint union of the $8$-cycle $a_1b_1a_2b_2a_3b_3a_4b_4a_1$ 
and the $4$-cycle $b_1b_2b_3b_4b_1$ by adding $a$ and the edges $ab_i$, $i\in [4]$. Suppose $|V(G_1)|\ge 7$. 
Then, for any $u_1, u_2 \in N(a) - \{b_1,b_2,b_3\}$, $G - \{av: v \not\in \{b_1,b_2,b_3,u_1,u_2\} \}$ contains  $TK_5$.
\end{lem}

Next we list a few results from \cite{HWY15I, HWY15II, HWY15III}. For convenience, we state their versions from \cite{HWY15III}. 
First, we need Theorem 1.1 in  \cite{HWY15III} to take care of the case when the vertex $x$ in 
$H=G/M$ (see Section 1) is a degree 3 vertex in a $K_4^-$ (which is 
$x_1$ in the lemma below).

\begin{lem}\label{x1a}
Let $G$ be \textsl{a} $5$-connected nonplanar graph and $x_1,x_2,y_1,y_2 \in V(G)$ be distinct such that 
$G[\{x_1 , x_2 , y_1 , y_2 \}]\cong K_4^-$ and $y_1y_2\notin E(G)$. Then one of the following holds:
\begin{itemize}
\item [$(i)$] $G$ contains \textsl{a} $TK_5$ in which $x_1$ is not \textsl{a}  branch vertex.
\item [$(ii)$] $G-x_1$ contains $K_4^-$, or $G$ contains \textsl{a}  $K_4^-$ in which $x_1$ is of degree 2. 
\item [$(iii)$] $x_2,y_1,y_2$ may be chosen so that for any distinct $z_0, z_1\in N(x_1)-\{x_2,y_1,y_2\}$, 
$G-\{x_1v:v\notin \{ x_2,y_1,y_2, z_0, z_1\}\}$ contains $TK_5$. 
\end{itemize}
\end{lem}

When applying the next three lemmas, the vertex $a$ will correspond to the vertex $x$ in $H=G/M$ in Section 1. 
The following result is Lemma 2.7 in \cite{HWY15III}, which deals with 5-separations with an apex side.

\begin{lem}\label{apex}
Let $G$ be \textsl{a} $5$-connected nonplanar graph and let $(G_1, G_2)$ be \textsl{a} $5$-separation in $G$. 
Suppose $|V(G_i)|\geq 7$ for $i\in [2]$, $a\in V(G_1\cap G_2)$, and $(G_2-a,V(G_1\cap G_2)-\{a\})$ is planar. Then one of the following holds:
\begin{itemize}
\item [$(i)$]  $G$ contains \textsl{a} $TK_5$ in which $a$ is not \textsl{a} branch vertex.
\item [$(ii)$] $G-a$ contains $K_4^-$, or $G$ contains \textsl{a} $K_4^-$ in which $a$ is of degree 2.
\end{itemize}
\end{lem}

The next result is Lemma 2.8 in  \cite{HWY15III},  which will be used to take care of 
5-cuts containing the vertices of a triangle. 

\begin{lem}
\label{5cut_triangle}
Let $G$ be \textsl{a} $5$-connected graph and $(G_1,G_2)$ be \textsl{a} $5$-separation in $G$. Suppose that $|V(G_i)|\ge 7$ for $i\in [2]$ and 
$G[V(G_1\cap G_2)]$ contains \textsl{a} triangle $aa_1a_2a$. Then one of the following holds:
\begin{itemize}
\item [$(i)$]  $G$ contains \textsl{a} $TK_5$ in which $a$ is not \textsl{a} branch vertex.
\item [$(ii)$] $G-a$ contains $K_4^-$, or $G$ contains \textsl{a} $K_4^-$ in which $a$ is of degree 2.
\item [$(iii)$] For any distinct $u_1,u_2,u_3\in N(a)-\{a_1,a_2\}$, $G-\{av: v \not\in \{a_1,a_2,u_1,u_2,u_3\}\}$ contains $TK_5$. 
\end{itemize}
\end{lem}

The following is Lemma 2.9 in \cite{HWY15III}. 

\begin{lem}\label{6-cut2}
Let $G$ be \textsl{a} graph, $A\subseteq V(G)$, and $a\in A$  such that $|A|=6$, $|V(G)|\geq 8$, 
$(G-a, A-\{a\})$ is planar, and $G$ is $(5, A)$-connected.
Then one of the following holds:  
\begin{itemize}
 \item[$(i)$] $G-a$ contains $K_4^-$, or $G$ contains \textsl{a} $K_4^-$ in which the degree of $a$ is $2$.
 \item[$(ii)$] $G$ has \textsl{a} $5$-separation $(G_1,G_2)$ such that $a\in V(G_1\cap G_2)$, $|V(G_2)|\ge 7$, $A\subseteq V(G_1)$, and 
$(G_2-a, V(G_1\cap G_2)-\{a\})$ is planar. 
\end{itemize}
\end{lem}

We need Theorem 1.4 in \cite{HWY15I}. This will be used to show that, for a quadruple $(T,S_T,A,B)$ in $H=G/M$ with  $x\in V(T)$ (see  Section 1), 
$x$ has a neighbor in $A$ (which corresponds to $G_1-G_2$ in the statement). 

\begin{lem} \label{N(x1)capA}
Let $G$ be \textsl{a} $5$-connected graph and $x\in V(G)$, and let 
$(G_1,G_2)$ be \textsl{a} $6$-separation in $G$ such that 
$x\in V(G_1\cap G_2)$,  $G[V(G_1\cap G_2)]$ contains a triangle
$xx_1x_2x$, $|V(G_i)|\ge 7$ for $i\in [2]$. Moreover, assume that
$(G_1,G_2)$ is chosen so that, subject to $\{x,x_1,x_2\} \subseteq
V(G_1\cap G_2)$ and $|V(G_i)|\ge 7$ for $i\in [2]$, $G_1$ is minimal.  
Let $V(G_1 \cap G_2) = \{x,x_1,x_2,v_1,v_2,v_3\}$. 
Then $N(x)\cap V(G_1-G_2)\ne \emptyset$, or  one of the following holds: 
\begin{itemize}
\item[$(i)$] $G$ contains \textsl{a} $TK_5$ in which $x$ is not \textsl{a} branch vertex.
\item[$(ii)$] $G$ contains $K_4^-$. 
\item[$(iii)$] There exists $x_3\in N(x)$ such that for any distinct $y_1,y_2\in N(x)-\{x_1,x_2,x_3\}$, 
$G-\{xv:v\notin \{x_1,x_2,x_3,y_1,y_2\}\}$ contains $TK_5$.
\item [$(iv)$] For some $i\in [2]$ and some $j \in [3]$, 
$N(x_i)\subseteq V(G_1-G_2)\cup \{x,x_{3-i}\}$, and  any three
independent paths in $G_1-x$ from $\{x_1,x_2\}$ to $v_1,v_2,v_3$, respectively, 
with two from $x_i$ and one from $x_{3-i}$, 
must contain \textsl{a} path from $x_{3-i}$ to $v_j$.
\end{itemize}
\end{lem}
We remark that  conclusion $(iv)$ in Lemma~\ref{N(x1)capA} will be
dealt with in Section 4, using a result on disjoint paths from 
 \cite{Yu03I, Yu03II, Yu03III}.
We also need Proposition 4.1 from \cite{HWY15I} to deal with the case
when $H/T$ is planar (see Section 1) for some $T\subseteq H$ with
$x\in V(T)$ and $T\cong K_2$ or $T\cong K_3$. 

\begin{lem}\label{contraction}
Let $G$ be \textsl{a} $5$-connected nonplanar graph, $x\in V(G)$, $T\subseteq G$ such that $x\in V(T)$, $T\cong K_2$ or $T\cong K_3$,   
$G/T$ is $5$-connected and  planar. Then $G-T$ contains $K_4^-$.
\end{lem}

We conclude this section with three additional results, first of which is a result of Seymour \cite{Se80};  equivalent versions are proved in \cite{Th80, CR79, Sh80}.

\begin{lem} \label{seymour-2dp}  
Let $G$ be \textsl{a} graph and let $s_1, s_2, t_1, t_2\in V(G)$ be distinct.
Then either $G$ contains disjoint paths from $s_1$ to $t_1$ and from $s_2$ to $t_2$, or $(G, s_1,s_2,t_1,t_2)$ 
is 3-planar.
\end{lem}

The second result is due to Perfect \cite{Pe68}. 

\begin{lem} \label{Perfect} 
Let $G$ be \textsl{a} graph, $u\in V(G)$, and $A\subseteq V(G-u)$. Suppose
there exist $k$ independent paths from $u$ 
to distinct $a_1,\ldots, a_k\in A$, respectively, and internally disjoint from $A$. Then for any $n\ge k$, if there exist $n$
independent paths  $P_1,\ldots, P_n$ in $G$ from $u$ to $n$ distinct
vertices in $A$ and internally  disjoint from $A$ then  $P_1,\ldots,
P_n$ may be chosen so that $a_i\in V(P_i)$ for $i\in [k]$. 
\end{lem}

The third result is due to Watkins and Mesner \cite{WM67}.
\begin{lem}
\label{Watkins}
Let $G$ be \textsl{a} $2$-connected graph and let $y_1, y_2, y_3$ be
three distinct vertices of $G$. Then $G$ has no cycle 
containing $\{y_1, y_2, y_3\}$ 
if, and only if, one of the following holds:
\begin{itemize}
\item [$(i)$] There exists \textsl{a} 2-cut $S$ in $G$ and there exist pairwise disjoint subgraphs $D_{y_i}$ of $G - S$, $i\in [3]$, such that 
$y_i\in V(D_{y_i})$ and each $D_{y_i}$ is \textsl{a} union of components of $G - S$. 
\item [$(ii)$] There exist 2-cuts $S_{y_i}$ in $G$, $i\in [3]$, and pairwise disjoint subgraphs $D_{y_i}$ of $G$, such that 
$y_i \in V(D_{y_i})$, each $D_{y_i}$ is \textsl{a} union of components of $G-S_{y_i}$, there exists $z\in S_{y_1} \cap S_{y_2} \cap S_{y_3}$, 
and $S_{y_1} - \{z\}, S_{y_2} - \{z\}, S_{y_3} - \{z\}$ are pairwise disjoint. 
\item [$(iii)$] There exist pairwise disjoint $2$-cuts $S_{y_i}$ in $G$ and pairwise disjoint subgraphs $D_{y_i}$ of 
$G - S_{y_i}$,  $i\in [3]$, such that $y_i \in V(D_{y_i})$, 
$D_{y_i}$ is \textsl{a} union of components of $G - S_{y_i}$, and $G - V(D_{y_1} \cup D_{y_2} \cup D_{y_3})$ has precisely two components, 
each containing exactly one vertex from $S_{y_i}$ for $i\in [3]$.
\end{itemize}
\end{lem}


\section{Obstruction to three paths}

In order to deal with $(iv)$ of Lemma~\ref{N(x1)capA}, we need a
result of the third author  \cite{Yu03I, Yu03II, Yu03III}, which characterizes 
graphs $G$ in which any three 
disjoint paths from $\{a,b,c\}\subseteq V(G)$ to $\{a',b',c'\}\subseteq V(G)$ must contain a path from $b$ to $b'$.
The objective of this section is to derive a much simpler version of that characterization by imposing extra
conditions on $G$. This result will be used several times in the proofs
of Lemmas~\ref{T_K_3_T'_K_2} and \ref{lemma_butterfly}.  To state the
result from \cite{Yu03I, Yu03II, Yu03III}, we need to describe {\it rungs}
and {\it ladders}.

Let $G$ be a graph, $\{a,b,c\}\subseteq V(G)$,
 and $\{a',b',c'\}\subseteq V(G)$. Suppose $\{a,b,c\}\ne \{a',b',c'\}$,  and
 assume that  $G$ has no separation $(G_1,G_2)$ such that $|V(G_1\cap G_2)|\le 3$, $\{a,b,c\}\subseteq V(G_1)$, and $\{a',b',c'\}\subseteq V(G_2)$. 
We say that $(G,(a,b,c), (a',b',c'))$ is a
 {\it rung} if one of the following holds:
\begin{enumerate}
\addtolength{\baselineskip}{-1ex}

\item [(1)] $b=b'$ or $\{a,c\}=\{a',c'\}$. 
\item [(2)] $a=a'$ and $(G-a,c,c',b',b)$ is 3-planar, or $c=c'$ and $(G-c,a,a',b',b)$ is 3-planar.

\item [(3)] $\{a,b,c\}\cap \{a',b',c'\}=\emptyset$ and $(G,a',b',c',c,b,a)$ is 3-planar.

\item [(4)] $\{a,b,c\}\cap \{a',b',c'\}=\emptyset$,
$G$ has a 1-separation $(G_1,G_2)$ such that $(i)$ 
$\{a,a',b,b'\}\subseteq V(G_1)$, $\{c,c'\}\subseteq V(G_2)$, and
$(G_1,a,a',b',b)$ is 3-planar, or $(ii)$ $\{c,c',b,b'\}\subseteq V(G_1)$, $\{a,a'\}\subseteq V(G_2)$, and
$(G_1,c,c',b',b)$ is 3-planar.

\item [(5)]  $\{a,b,c\}\cap \{a',b',c'\}=\emptyset$,  and $G$ has a separation $(G_1,G_2)$ such that  $V(G_1\cap 
G_2)=\{z,b\}$ (or $V(G_1\cap G_2)=\{z,b'\}$), and $(i)$ $(G,a,a',b',b)$ is
3-planar, $\{a,a',b,b'\}\subseteq V(G_1)$,
$\{c,c'\}\subseteq V(G_2)$, and $(G_2,c,c',z,b)$ (or
$(G_2,c,c',b',z)$) is 3-planar, or $(ii)$  $(G,c,c',b',b)$ is
3-planar, $\{c,c',b,b'\}\subseteq V(G_1)$,
$\{a,a'\}\subseteq V(G_2)$, and $(G_2,a,a',z,b)$ (or 
$(G_2,a,a',b',z)$) is 3-planar. 

\item [(6)] $\{a,b,c\}\cap \{a',b',c'\}=\emptyset$, and 
there are pairwise edge disjoint subgraphs $G_a,G_c,M$ of $G$ such that 
$G=G_a\cup G_c\cup M$, $V(G_a\cap M)=\{u,w\}$, $V(G_c\cap M)=\{p,q\}$, $V(G_a\cap G_c)=
\emptyset$, and $(i)$ $\{a,a',b'\}\subseteq V(G_a)$, $\{c,c',b\}\subseteq
V(G_c)$, and $(G_a,a,a',b',w,u)$ and 
$(G_c,c',c,b,p,q)$ are 3-planar, or $(ii)$ $\{a,a',b\}\subseteq V(G_a)$, $\{c,c',b'\}\subseteq
V(G_c)$, $(G_a,b,a,a',w,u)$ and 
$(G_c,b',c',c,p,q)$ are 3-planar.

\item [(7)] $\{a,b,c\}\cap \{a',b',c'\}=\emptyset$, and 
there are pairwise edge disjoint subgraphs $G_a,G_c,M$ of $G$ such that 
$G=G_a\cup G_c\cup M$, $V(G_a\cap M)=\{b,b',w\}$, $V(G_c\cap M)=\{b,b',p\}$, $V(G_a\cap G_c)=
\{b,b'\}$, $\{a,a'\}\subseteq V(G_a)$, $\{c,c'\}\subseteq
V(G_c)$, and $(G_a,a,a',b',w,b)$ and $(G_c,c',c,b,p,b')$ are 3-planar. 
\end{enumerate}

Let $L$ be a graph and let $R_1,\ldots, R_m$ be edge disjoint 
subgraphs of $L$ such that 

\begin{itemize}
\item [$(i)$] 
$(R_i,(x_{i-1},v_{i-1},y_{i-1}),(x_i,v_i,y_i))$ is a rung for
each $i\in [m]$, 
\item [$(ii)$] $V(R_i\cap R_j)=\{x_i,v_i,y_i\}\cap \{x_{j-1},v_{j-1},y_{j-1}\}$
for $i,j\in [m]$ with $i<j$, 
\item [$(iii)$] for any distinct $i,j\in [m]$, if $x_i=x_j$ 
then $x_k=x_i$  for all $i\le k\le
j$, if $v_i=v_j$ then  $v_k=v_i$ for all $i\le k\le
j$, and if $y_i=y_j$ then $y_k=y_i$ for all $i\le k\le
j$, and 
\item [$(iv)$]   $L=(\bigcup_{i=1}^mR_i)+S$, where
$S$ consists of those edges of $L$ each of which has both ends in 
$\{x_i,v_i,y_i\}$  for some $i\in [m]$.
\end{itemize}
Then  $(L,(x_0,v_0,y_0), (x_m,v_m,y_m))$ is
a {\it ladder with rungs}
$(R_i,(x_{i-1},v_{i-1},y_{i-1})$, $(x_i,v_i,y_i))$, $i\in [m]$,  
or simply, a {\it ladder along} $v_0\ldots v_m$. 

By the definition of a rung, we see that a ladder $(L,(x_0,v_0,y_0),
(x_m,v_m,y_m))$ has three disjoint paths from $\{x_0,v_0,y_0\}$ to $\{x_m,v_m,y_m\}$.

For a sequence $W$, the {\it reduced sequence} of $W$ is the sequence obtained from 
$W$ by removing  all but one consecutive identical elements. For example, the reduced 
sequence of $aaabcca$ is $abca$. We can now state the main result in \cite{Yu03III}.

\begin{lem}\label{characterization}
Let $G$ be \textsl{a} graph,  $\{a,b,c\}\subseteq V(G)$,
 and $\{a',b',c'\}\subseteq V(G)$ such that $\{a,b,c\}\ne
 \{a',b,'c'\}$. Assume that, for any $T\subseteq V(G)$ with $|T|\le 3$, every
component of $G-T$ contains some element of $\{a,b,c\}\cup \{a',b',c'\}$. 
Then any three disjoint paths in $G$ from $\{a,b,c\}$ to
$\{a',b',c'\}$ must include one from $b$ to $b'$  if, and only if, one of the
following statements holds:
\vspace*{-1ex}
\begin{itemize}
\item [$(i)$] $G$ has \textsl{a} separation $(G_1,G_2)$ of order at most 2 such that
$\{a,b,c\}\subseteq V(G_1)$ and $\{a',b',c'\}\subseteq V(G_2)$.
\vspace*{-1ex}
\item [$(ii)$] $(G,(a,b,c),(a',b',c'))$ is a ladder.
\vspace*{-1ex}
\item [$(iii)$] $G$ has \textsl{a} separation $(J,L)$ such that $V(J\cap
L)=\{w_0,\ldots, w_n\}$, 
$(J,w_0,\ldots,w_n)$ is 3-planar, $\{a,b,c\}\cup \{a',b',c'\}\subseteq V(L)$,
$(L,(a,b,c),(a',b',c'))$ is \textsl{a} ladder along \textsl{a} sequence $v_0\ldots
v_m$, where $v_0=b$, $v_m=b'$, and 
$w_0\ldots w_n$ is the reduced sequence of $v_0 \ldots v_m$.
\end{itemize}
\end{lem}

We may view $(ii)$ as a special case of $(iii)$ by letting $J$ be a
subgraph of $L$.   In the applications of Lemma~\ref{characterization}
in this paper, we will consider rungs and ladders in a 5-connected graph without  $TK_5$. 
With such extra conditions, the rungs have much simpler structure, as
given in the next two  lemmas.

\begin{lem}\label{simplerungs-1}
Let $G$ be \textsl{a} 5-connected graph and $(R,R')$ \textsl{a}
separation in $G$ such that  $|V(R')|\ge 8$, $V(R\cap R')
=\{a,b\}\cup \{a',b',c'\}$, $a\ne b$, and $a',b',c'$ are pairwise distinct. Let $R^*$ be obtained from $R$ by adding the new
vertex $c$ and joining $c$ to each neighbor of $a$ in $R$ with an
edge, and assume  $(R^*, (a,b,c), (a',b',c'))$ is a rung.
Then $b=b'$, $V(R)=\{a,b,a',c'\}$ and $E(R)=\{aa',ac'\}$.
\end{lem}

\pf Since $a$ and $c$ have the same set of neighbors in $R^*$ and
$(R^*, (a,b,c), (a',b',c'))$ is a rung, it follows from the definition
of a rung that  $(R^*, (a,b,c), (a',b',c'))$ is of type (1) or
(2). Then, since $G$ is 5-connected, $V(R)=\{a,b\}\cup
\{a',b',c'\}$. 

Suppose $(R^*, (a,b,c), (a',b',c'))$ is of type (2). By symmetry, we
may assume that $c=c'$ and $(G-c,a,a',b',b)$ is 3-planar. Then
$ab'\notin E(G)$ or $a'b\notin E(G)$. Hence, $\{a',b,c\}$ or
$\{a,b',c\}$ would be a cut in $R^*$ separating $\{a,b,c\}$ from
$\{a',b',c'\}$, a contradiction.

So  $(R^*, (a,b,c), (a',b',c'))$ is of type (1). Then, since $R^*$ has
no separation of order at most 3 separating $\{a,b,c\}$ from
$\{a',b',c'\}$, we deduce that $a\ne a'$, $c\ne
c'$, and $E(R)=\{aa',ac'\}$. \qed

\medskip
Note that the conclusion of Lemma~\ref{simplerungs-1} is a special
case of $(i)$ of the next lemma.

\begin{lem}\label{simplerungs}
Let $G$ be \textsl{a} 5-connected graph and $(R,R')$ \textsl{a} separation in $G$ such that  $|V(R')|\ge 8$, $V(R\cap R') =\{a,b,c\} \cup \{a',b',c'\}$,
$\{a,b,c\}\ne \{a',b',c'\}$, and $(R, (a,b,c), (a',b',c'))$ is a rung. 
Then $G$ contains $TK_5$ or $K_4^-$, or 
one of the following holds: 
\begin{itemize} 
\item [$(i)$] $b=b'$.
\item [$(ii)$] $\{a,c\}=\{a',c'\}$, $V(R) = \{a,c,b,b'\}$, and $E(R)=\{bb'\}$. 
\item [$(iii)$]  $V(R)-(\{a,b,c\}\cup \{a',b',c'\})=\{v\}$ and $N(v)=\{a,b,c\}\cup \{a',b',c'\}$, and either $a=a'$ and $E(R-v)=\{bb',cc'\}$ or 
$c=c'$ and $E(R-v)=\{bb',aa'\}$.
\item [$(iv)$] $\{a,b,c\}\cap \{a',b',c'\}=\emptyset$, $V(R)-\{a,a',b,b',c,c'\}=\{v\}$, $N(v)=\{a,a',b,b',c,c'\}$, and  $E(R-v)=\{aa', bb', cc'\}$. 
\end{itemize}
\end{lem}

\pf Without loss of generality, let  $A,B,C$ be disjoint paths in $R$ from $a,b,c$ to $a',b',c'$, respectively. First, we consider the case when  $\{a,b,c\}\cap \{a',b',c'\}\ne \emptyset$. If $b=b'$ then ($i$) holds; so we may  assume 
$b\ne b'$. If $a=a'$ and $c=c'$ then,  
 since $G$ is 5-connected, $V(R)=\{a,b,b',c\}$; so $bb'\in E(R)$
 (because of the paths $A,B,C$), and we have ($ii$). Thus by symmetry between $\{a,a'\}$ and $\{c,c'\}$, we may assume $c\ne c'$. 
Suppose $a=a'$. Then by the definition of  a rung, $R-a$ has no
disjoint paths from $b,c$ to $c',b'$, respectively. So by Lemma~\ref{seymour-2dp},
$(R-a,c,c',b',b)$ is 3-planar.  Since $G$ is 5-connected,
$(R-a,c,c',b',b)$ is in fact planar.
   If $|V(R)|\ge 7$ then $G$ contains $TK_5$ or $K_4^-$ by
   Lemma~\ref{apex}, using
the separation $(R,R')$. If 
$V(R)=\{a,b,b',c,c'\}$ then, since $(R-a,c,c',b',b)$ is planar, either $\{a,b,c'\}$ or $\{a,b',c\}$ is a 3-cut in $R$ separating 
$\{a,b,c\}$ from $\{a',b',c'\}$, contradicting the definition of
a rung. Thus, we may assume $|V(R)|=6$ and let $v\in 
V(R)- \{a,b,b', c,c'\}$. Since $G$ is 5-connected, $N(v) =
\{a,b,b',c,c'\}$. Since $(R-a, c,c',b',b)$ is planar, $bc',cb'\notin E(R)$. So 
$bb',cc'\in E(R)$, as otherwise $\{a,v,c\}$ or $\{a,v,b\}$ would be a
3-cut in $R$ separating $\{a,b,c\}$ from $\{a',b',c'\}$, 
contradicting the definition of a rung. Hence, ($iii$) holds. 

Thus, we may assume that $\{a,b,c\}\cap \{a',b',c'\}=\emptyset$. We
need to deal with (3) -- (7) in the definition of a rung. We deal with
(4)--(7) in order, and treat  (3) last (which  is the most complicated case where we use the discharging technique).

Suppose (4) holds for $(R, (a,b,c), (a',b',c'))$. By symmetry, assume that  $R$ has a 1-separation $(G_1,G_2)$ such that 
$\{a,a',b,b'\}\subseteq V(G_1)$, $\{c,c'\}\subseteq V(G_2)$, and
$(G_1,a,a',b',b)$ is 3-planar. Let $V(G_1 \cap G_2) = \{v\}$. Since $G$ is 5-connected, $(G_1, a,a',b',b)$ is planar and  
$V(G_2)= \{v,c,c'\}$. Moreover, $vc,vc',cc'\in E(G)$; for otherwise $R$ would have a separation $(R_1,R_2)$ such that 
$\{a,b,c\}\subseteq V(R_1)$, $\{a',b',c'\}\subseteq V(R_2)$, and $V(R_1\cap R_2)\in \{\{a,b,c'\}, \{a',b',c\},\{a,b,v\}\}$. 
If $|V(G_1)|\ge 7$ then the assertion follows from Lemma~\ref{apex}, 
using the separation $(G_1,G_2\cup R')$. So we may assume $|V(G_1)|\le
6$. If $|V(G_1)|=6$ then let $t\in
V(G_1)-\{a,a',b,b',v\}$;  now $N(t)=\{a,a',b,b',v\}$ and $|(N(v)-\{c,c'\})\cap N(t)|\ge 2$ (since $G$ is 5-connected), and hence
$R$ (and therefore $G$) contains $K_4^-$.  So we may assume $V(G_1) = \{a,a',b,b',v\}$. Then $va'\in E(G)$; otherwise $N(v)=\{a,b,b',c,c'\}$ and, hence, $a'b\notin E(G)$ (as
$(G_1,a,a',b',b)$ is planar), which implies that $\{a,b',c'\}$ is a cut in $R$ separating $\{a,b,c\}$ from $\{a',b',c'\}$, a contradiction. 
Similarly, $va,vb,vb'\in E(G)$. Then by planarity of
$(G_1,a,a',b',b)$, we have $ab',ba'\notin E(G)$. So $aa',bb'\in E(G)$ as
$\{c,v,b'\}$ and $\{a,v,c\}$ are not  3-cuts in $R$ separating $\{a,b,c\}$ from $\{a',b',c'\}$.  Thus we have ($iv$).

Suppose (5) holds for $(R, (a,b,c), (a',b',c'))$, and assume by symmetry that $(R,a,a',b',b)$ is
3-planar, and $R$ has a separation $(G_1,G_2)$ such that $V(G_1\cap 
G_2)=\{z,b\}$, $\{a,a',b,b'\}\subseteq V(G_1)$, $\{c,c'\}\subseteq V(G_2)$, and $(G_2,c,c',z,b)$ is 3-planar. 
Since $G$ is 5-connected, $V(G_2) = \{b,z,c,c'\}$. Then $cz,cc'\in E(G)$ as otherwise, $\{a,b,c'\}$ or $\{a,b,z\}$ would be a 3-cut in $R$ separating 
$\{a,b,c\}$ from $\{a',b',c'\}$. Hence, since $(G_2,b,z,c',c)$ is
planar, $bc'\notin E(G)$. Since $(R, a, a',b',b)$ is 3-planar,
$(G_1,a,a',b',b)$ is 3-planar. Thus, the separation $(G_1,G_2-b)$ shows  that $(R, (a,b,c), (a',b',c'))$ is of type 
(4); so  we may assume that
($iv$) holds by the argument in the previous paragraph.

Now suppose (6) holds for $(R, (a,b,c), (a',b',c'))$, and, by symmetry, assume that 
there are pairwise edge disjoint subgraphs $G_a,G_c,M$ of $R$ such that 
$R=G_a\cup G_c\cup M$, $V(G_a\cap M)=\{u,w\}$, $V(G_c\cap M)=\{p,q\}$, $V(G_a\cap G_c)=
\emptyset$, $\{a,a',b'\}\subseteq V(G_a)$, $\{c,c',b\}\subseteq
V(G_c)$, and $(G_a,a,a',b',w,u)$ and $(G_c,c',c,b,p,q)$ are 3-planar. Since $G$ is 5-connected,  we have $V(M)=\{p,q,u,w\}$, and $(G_a,a,a',b',w,u)$ and $(G_c,c',c,b,$ $p,q)$ are planar. We may assume that 
$|V(G_c) - \{b,c,c',p,q\}| \leq 1$ and $|V(G_a) - \{a,a',b',u,w\}|
\leq 1$, as otherwise the assertion follows from  Lemma~\ref{apex}
with the separation $(G_c,G_a\cup M\cup R')$ or $(G_a,G_c\cup M\cup R')$. 
If there exists $v\in V(G_c) - \{b,c,c',p,q\}$ then, since $G$ is 5-connected, $N(v)=  \{b,c,c',p,q\}$ and $|(N(p)-\{u,w\})\cap\{b,c,c',q\}|\ge 2$; 
so $R$ (and hence $G$) contains $K_4^-$. Thus we may assume  $V(G_c) = \{b,c,c',p,q\}$. Since $G$ is 5-connected, $p$ and $q$ each have at least five neighbors in $G_c\cup M$. 
Hence, since $(G_c,b,c,c',q,p)$ is planar, $N(p)=\{u,w,b,c,q\}$ and
$N(q)=\{u,w,c,c',p\}$; so $G[\{p,q,u,w\}]$ (and hence $G$) contains $K_4^-$.

Suppose (7) holds for $(R, (a,b,c), (a',b',c'))$. Then 
there are pairwise edge disjoint subgraphs $G_a,G_c,M$ of $R$ such that 
$R=G_a\cup G_c\cup M$, $V(G_a\cap M)=\{b,b',w\}$, $V(G_c\cap M)=\{b,b',p\}$, $V(G_a\cap G_c)=
\{b,b'\}$, $\{a,a'\}\subseteq V(G_a)$, $\{c,c'\}\subseteq V(G_c)$, and $(G_a,a,a',b',w,b)$  and 
$(G_c,c',c,b,p,b')$ are 3-planar.  Since $G$ is 5-connected, we have  $V(M)=\{b,b',p,w\}$, and $(G_a,a,a',b',w,b)$  and 
$(G_c,c',c,b,p,b')$ are actually planar. 
If  $|V(G_c)|\ge 7$ then the assertion follows from Lemma~\ref{apex}
with the separation $(G_c,G_a\cup M\cup R')$. So we may assume $|V(G_c)|\le 6$. 
If there exists $q\in V(G_c)-\{b,b',c,c',p\}$ then $N(q)=\{b,b',c,c',p\}$ (as $G$ is 5-connected); therefore, since $(G_c,c',c,b,p,b')$ 
is planar, $N(p)\subseteq \{b,b',w,q\}$, a contradiction. Thus  $V(G_c)=\{b,b',c,c',p\}$ and, hence, $N(p)=\{b,b',c,c',w\}$. 
Similarly, by considering $G_a$, we may assume
$N(w)=\{a,a',b,b',p\}$. Thus $G[\{b,b',p,w\}]$ (and hence $G$) contains $K_4^-$.

Finally, assume that (3) holds for $(R, (a,b,c), (a',b',c'))$. So $(R,a',b',c',c,b,a)$ is planar (as $G$ is 5-connected), and we may assume that 
$R$ is embedded in a closed disc with no edge crossings such that $a,b,c,c',b',a'$ occur on the boundary of the disc in clockwise order. 
We apply the discharging method. For convenience, let $A=\{a,b,c,a',b',c'\}$, $F(R)$ denote the set of faces of $R$, and $f_{\infty}$ denote 
the outer face of $R$ (which is incident with all vertices in $A$). For each $f\in F(R)$, let $d_R(f)$ denote the number of incidences of the edges of $R$ with $f$, and 
$\partial f$ denote the set of vertices of $R$ incident with $f$. 
For $x\in V(R)\cup F(R)$, let $\sigma(x) = d_R(x) - 4$ be the charge
of $x$. Note that $R$
is connected as in $R$ there is no separation $(R_1,R_2)$ of order at
most 3 such that $\{a,b,c\}\subseteq V(R_1)$ and
$\{a',b',b'\}\subseteq V(R_2)$. Hence, by Euler's formula, 
$\sum_{x \in V(R)\cup F(R)} \sigma(x)=-8$.

We redistribute charges according to the following rule: For each  $v \in V(R)-A$, 
$v$ sends $1/2$ to each $f\in F(R)$ that is incident with $v$ and has
$d_R(f)=3$. Let $\tau(x)$ denote the new charge for all $x\in V(R)\cup
F(R)$. Then $$\sum_{x \in V(R)\cup F(R)} \tau(x)=\sum_{x\in V(R)\cup
  F(R)}\sigma(x)=-8.$$ 

Note that we may assume $K_4^-\not\subseteq G$. Thus, each $v\in V(R)-A$ is incident with at most $\lfloor d_R(v)/2 \rfloor$ faces 
$f\in F(R)$ with $d_R(f)=3$; so $\tau(v) \geq 0$ (as $d_R(v)\ge 5$). 
Moreover, for $f \in F(R)$, $\tau(f) \geq 0$ unless $d_R(f)=3$  and $f$ is incident with at least two vertices in $A$.

Since $R$ has no separation $(R_1,R_2)$ of order at most 3 such that
$\{a,b,c\}\subseteq V(R_1)$ and $\{a',b',c'\}\subseteq V(R_2)$, we see that 
$\{a,b,c\}$ and $\{a',b',c'\}$ are independent in $R$. Moreover, since
$(R,a,a',b',c',c,b)$ is planar, 
$ab',ac',ba',bc',ca',cb'\notin E(R)$, and $d_R(v)\ge 2$ for $v\in A$. 
Hence, $bb'\notin E(R)$; otherwise, since $G$ is 5-connected, $V(R)=A$ (to avoid 4-cuts $\{a,a',b,b'\}$ and $\{b,b',c,c'\}$), which in turn would 
force $d_R(v)\le 1$ for some $v\in A$. 

Therefore, $d_R(f_{\infty})\ge 10$, and if $f\in F(R)$ with $d_R(f)=3$ and $|\partial f\cap A|\ge 2$ then $\partial f\cap A=\{a,a'\}$ 
or  $\partial f\cap A=\{c,c'\}$.
Hence, 
\begin{eqnarray*}
\sum_{x \in V(R)\cup F(R)} \tau(x) & \geq & \sum_{v\in  V(R)} \tau(v) + \sum_{f \in F(R), |\partial f \cap A| \geq 2}  \tau(f) \\
& \geq & \sum_{v\in  A} (d_R(v) - 4) + ( d_R(f_{\infty}) - 4 )+ \sum_{d_R(f)=3, |\partial f \cap A| \geq 2}  (d_R(f)-4)   \\
& \geq & (-12)+(10-4)+(-1)\times 2 \\
& =& -8. 
\end{eqnarray*} 
Thus, all the inequalities above hold with equality. In particular,  $d_R(f_{\infty})=10$, $d(x)=2$ for $x\in A$, 
and there exist $u,v\in V(R)-A$ such that $uaa'u$ and $vcc'v$ are triangles and 
$aa'ub'vc'cvbua$ is the outer walk of $R$. Since $G$ is 5-connected and $(R,a,b,c,c',b',a')$ is planar, 
$V(R)=A\cup \{u,v\}$ and $uv\in E(R)$. Hence,
$G[\{b,b',u,v\}]\cong K_4^-$, a contradiction. \qed

\section{Quadruples and special structure} 

As mentioned in Section 1, we  need to deal with 5-connected
graphs in which  every edge or triangle at a given vertex is contained
in a cut of size 5 or 6. Thus,   for convenience, we introduce the
following concept of quadruple. 
  
Let $G$ be a graph. For $x\in V(G)$, let $\mathcal{Q}_x$ denote the set of all quadruples $(T,S_T,A,B)$,  such that 
\begin{itemize}
\item[(1)] $T\subseteq G$, $x\in V(T)$, and $T\cong K_2$ or $T\cong K_3$, 
\item[(2)] $S_T$ is a cut of $G$ with $V(T) \subseteq S_T$, $A$ is a
  nonempty union of components of $G-S_T$, and $B = G-A-S_T\ne \emptyset$,
\item[(3)] if $T\cong K_3$ then  $5\le |S_T| \le 6$, and
\item [(4)] if $T\cong K_2$ then $|S_T|=5$, $|V(A)|\ge 2$,  and
  $|V(B)|\ge 2$. 

 \end{itemize}

 The purpose of this section is  to  derive useful  properties of
 quadruples, in particular, those $(T,S_T,A,B)$ that minimize
 $|V(A)|$. We begin with a few simple properties, first of which gives a
 bound on  $|V(A)|$. 
\begin{lem}
\label{cut_4-}
Let $G$ be \textsl{a} $5$-connected graph, $x\in V(G)$, and
$(T,S_T,A,B) \in \mathcal{Q}_x$. Then $G$ contains $K_4^-$, or
$|V(A)|\ge 5\le |V(B)|$.
\end{lem}

\pf Suppose 
there exists $(T,S_T,A,B) \in \mathcal{Q}_x$ such that $|V(A)|\le 4$
or $|V(B)|\le 4$. We choose such $(T,S_T,A,B) \in \mathcal{Q}_x$ with
$|V(A)|$ minimum. Then $|V(A)|\le 4$. 
Let $\delta$ denote the minimum degree of $A$, and let $u\in V(A)$ such that $u$ has degree $\delta$ in  $A$. 

We may assume $\delta \ge 1$. For, suppose $\delta=0$. 
 If $T\cong K_3$ then, since $G$ is 5-connected,  $|N(u)\cap S_T|\ge
 5$; so $G[T+u]$ contains $K_4^-$. Hence we may assume
 $T\cong K_2$. Then $|V(A)|\ge 2$. In fact, by the minimality of
 $|V(A)|$, $|V(A)|=2$ and $A$ consists of two isolated vertices. Now
 $G[A\cup T]$ contains $K_4^-$. 

\medskip
{\it Case} 1. $\delta =1$. 

Then  $|N(u)\cap S_T|\ge 4$.  Let $v$ be the
unique neighbor of $u$ in $A$. Since $|V(A)|
\le 4$ and $G$ is 5-connected, $|N(v)\cap S_T|\ge 2$.
We may assume  $|N(u)\cap N(v)\cap S_T|\le 1$; for, otherwise,  $G[S_T\cup \{u,v\}]$ 
contains $K_4^-$. 

Suppose $|N(v)\cap S_T|\ge 3$ or  $N(u)\cap N(v)\cap S_T=\emptyset$. Then $|S_T|=6$ and,
hence, $T\cong K_3$. Therefore, $|N(u)\cap V(T)|\ge 2$ or $|N(v)\cap
V(T)|\ge 2$; so $G[T+u]$ or $G[T+v]$ contains $K_4^-$. 

Hence, we may assume that  $|N(v)\cap S_T|\le 2$ and  $|N(u)\cap
N(v)\cap S_T| =1$. Then, since
$|V(A)|\le 4$ and $G$ is 5-connected,  $|N(v)\cap S_T|=2$, $|N(v)\cap
V(A)|=3$, and $|V(A)|=4$. Let $v_1,v_2\in
V(A)-\{u,v\}$, and let $w\in N(u)\cap N(v)\cap S_T$. Since $G$ is 5-connected, $|N(v_i)\cap S_T|\ge 3$ for $i\in [2]$.

We may assume $w\notin V(T)$; for, if $w\in V(T)$ then $|V(T)\cap
N(u)|\ge 2$ or $|V(T)\cap N(v)|\ge 2$, and  $G[T+\{u,v\}]$ contains
$K_4^-$.  We
may also assume $w\notin N(v_i)$ for $i\in [2]$, as otherwise
$G[\{u,v,w,v_i\}]$ contains $K_4^-$. 

If $v_1v_2\notin E(G)$ then $|N(v_i)\cap S_T|\ge 4$ for $i\in [2]$; so
$|N(v_i)\cap V(T)|\ge 2$ for $i\in [2]$ (since $w\notin N(v_i)$ and $w\notin V(T)$), and hence, $G[T+\{v_1,v_2\}]$ contains $K_4^-$. So assume
$v_1v_2\in E(G)$. Since $G$ is 5-connected and $w\notin N(v_i)$ for
$i\in [2]$,  there exists $w'\in N(v_1)\cap N(v_2)\cap S_T$. Now $G[\{v,v_1,v_2,w'\}]$ contains $K_4^-$.

\medskip

{\it Case} 2. $\delta\ge 2$.

 If $|V(A)|=3$ then $A\cong K_3$ and, since $G$ is 5-connected, $|N(a)\cap S_T|\ge 3$
 for all $a\in V(A)$; hence, since $|S_T|\le 6$,  $G[V(A)\cup S_T]$ 
contains $K_4^-$. So assume  $|V(A)|=4$. We may further assume
that $A$ is a cycle as otherwise $A$ contains $K_4^-$.  
 Moreover, we may assume that for any $st\in E(A)$, $|N(s)\cap N(t)\cap
S_T|\le 1$; for otherwise $G[\{s,t\}\cup S_T]$ contains $K_4^-$. Let $A=uvwru$. 

Suppose $T\cong K_2$. Then for any $st\in E(A)$, $(N(s)\cup N(t))-V(A)=S_T$
and  $|N(s)\cap N(t)\cap S_T|= 1$. Let $S_T=\{x_1,x_2,x_3,x_4,x_5\}$
and, without loss of generality, let $N(u)\cap A=\{x_1,x_2,x_3\}$ and
$N(v)\cap A=\{x_3,x_4,x_5\}$. Since  $(N(w)\cup N(r))-V(A)= S_T$,
$wx_3\in E(G)$ or $rx_3\in E(G)$. Then $G[\{u,v,w,x_3\}]\cong K_4^-$ or $G[\{r,u,v,x_3\}]\cong
K_4^-$.

Now assume $T\cong K_3$. Let $S_T=\{x_1,x_2,x_3,x_4,x_5,x_6\}$ such
that $V(T)=\{x_1,x_2,x_3\}$. 
We may assume $|N(a)\cap V(T)|\le 1$ for each $a \in V(A)$, for, otherwise, $G[T+a]$ contains $K_4^-$. 
Hence, let $x_4,x_5\in N(u)$, $x_5,x_6\in N(v)$, and $x_6,x_4\in N(w)$. Note that $N(r)\cap \{x_4,x_6\}\ne \emptyset$. 
 If $x_4\in N(r)$ then $G[\{u,w,r,x_4\}]\cong K_4^-$, and if $x_6\in N(r)$ then 
$G[\{v,w,r,x_6\}]\cong K_4^-$. \qed

\medskip

Next, we show that if a graph $G$ has no contractible edge or triangle
at some vertex $x$
then every edge of $G$ at $x$ is associated with a quadruple in ${\cal Q}_x$.

\begin{lem}
\label{existence_of_T'}
Let $G$ be \textsl{a} 5-connected graph and $x \in V(G)$.
Suppose for any $T\subseteq G$ with $x\in V(T)$ and $T\cong K_2$ or
$T\cong K_3$, $G/T$ is not 5-connected. 
Then for any $ax \in E(G)$, there exists $(T',S_{T'},C,D) \in \mathcal{Q}_x$ such that $\{a,x\} \subseteq V(T')$.
\end{lem}

\pf 
Let $T_1 = ax$. By assumption, $G/T_1$ is not 5-connected. So there exists a $5$-cut $S_{T_1}$ in $G$ with $V(T_1)\subseteq S_{T_1}$.
We may assume that $G - S_{T_1}$ has a trivial component; for otherwise, let $C$ be a component of $G-S_{T_1}$ and $D=(G-S_{T_1})-C$. 
Then $(T_1,S_{T_1},C,D)\in {\cal Q}_x$ is the desired quadruple. 

So let $y\in V(G)$ such that $y$ is a component of $G-S_{T_1}$. 
Let $T_2: = G[T_1 + y] \cong K_3$. By assumption, $G/T_2$ is not
$5$-connected. So there exists a cut $S_{T_2}$ in $G$ such that 
$V(T_2) \subseteq S_{T_2}$ and $|S_{T_2}|\in \{5,6\}$.  Let $C$ be a component of $G-S_{T_2}$ and $D=(G-S_{T_2})-C$. 
Then $(T_2,S_{T_2},C,D)\in {\cal Q}_x$ is the desired quadruple.
\qed

\medskip
We now show that if $(T,S_T,A,B)$ is chosen to minimize $|V(A)|$ then
we may assume $T\cong K_3$.

\begin{lem}
\label{T_K_2}
Let $G$ be \textsl{a} 5-connected graph and $x \in V(G)$.
Suppose for any $T\subseteq G$ with $x\in V(T)$ and $T\cong K_2$ or
$T\cong K_3$, $G/T$ is not 5-connected. 
Then $G$ contains $K_4^-$, or for any $(T,S_T,A,B)\in {\cal Q}_x$ with $|V(A)|$ minimum, $T\cong K_3$.
\end{lem}

\pf Let $(T,S_T,A,B)\in {\cal Q}_x$ with $|V(A)|$ minimum, and assume $T \cong K_2$. Then $|S_T| = 5$. Let $a\in N(x)\cap V(A)$.
By Lemma \ref{existence_of_T'}, there exists $(T',S_{T'},C,D) \in \mathcal{Q}_x$ such that $\{a,x\} \subseteq V(T')$. 
Note that $T' \cong K_2$ and $|S_{T'}| = 5$, or $T' \cong K_3$ and
$|S_{T'}| \in \{5,6\}$. 
We may assume $|V(A)|\ge 5$; for, if not, then $G$ contains $K_4^-$ by Lemma~\ref{cut_4-}.

We may assume that if $A\cap C\ne \emptyset$ then  $|(S_{T'} \cup S_T) -
V(B\cup D) | \geq |S_{T'}|+1$. For, suppose $A\cap C\ne \emptyset$ and $|(S_{T'} \cup S_T) -
V(B\cup D) | \le |S_{T'}|$. 
 If $|V(A\cap C)|\ge 2$ or $T'\cong K_3$ then $(T', (S_{T'} \cup S_T) - V(B\cup D), A \cap C, B \cup
D) \in \mathcal{Q}_x$ and $|V(A \cap C)| \leq |V(A) - \{a\}| <
|V(A)|$, contradicting the choice of $(T,S_T,A,B)$ that $|V(A)|$ is
minimum. So assume 
$|V(A\cap C)|=1$ and $T'\cong K_2$.  Then $|(S_{T'} \cup S_T) -V(B\cup D) |
=|S_{T'}|=5$ and $|V(C)|\ge 2\le |V(D)|$.
Assume for the moment $A \cap D =
\emptyset$. By Lemma~\ref{cut_4-},  we may assume $|S_{T'} \cap V(A)|
= 4$ (as $|S_{T'}|=5$ and $|V(A)|\ge 5$); so
$|S_{T'} \cap V(B)| = 0$, $|S_T\cap V(C)|=0$, and 
$|S_{T'} \cap S_T| = 1$.  Since $|V(C)|\ge 2$,   $B \cap C \neq \emptyset$. 
So $S_T \cap S_T'$ is a 1-cut in $G$, contradicting the assumption
that $G$ is 5-connected.  Hence,  $A\cap D\ne \emptyset$. 
We may assume  $|V(A\cap D)|\ge 2$; as otherwise,  since $G$ is
5-connected, $G[(A\cap C)\cup (A\cap D)\cup \{a,x\}]\cong K_4^-$. 
Then $|(S_{T'} \cup S_T) - V(B\cup C)| \geq |S_{T'}|+1$; otherwise, 
$(T', (S_{T'} \cup S_T) - V(B \cup C), A \cap D, B \cup C) \in
\mathcal{Q}_x$ and $2\le |V(A \cap D)| < |V(A)|$, contradicting 
the choice of $(T,S_T,A,B)$ that $|V(A)|$ is minimum.
Hence,  $|(S_{T'} \cup S_T) - V(A \cup D) | = |S_T| + |S_{T'}| - |(S_{T'}
\cup S_T) - V(B \cup C) | \leq 4$.
Since $G$ is $5$-connected, $B \cap C = \emptyset$. Since $|(S_{T'}
\cup S_T) -V(B\cup D) |=5$, $|S_T\cap V(C)|\le 3$.  Therefore,
$|V(C)|\le 4<|V(A)|$, a contradiction. 

Similarly, we may assume that if $A\cap D\ne \emptyset$ then $|(S_{T'} \cup S_T) - V(B\cup
C)| \geq |S_{T'}|+1$. 

Suppose $A\cap C=A\cap D=\emptyset$. Then, since $|V(A)|\ge 5$ and
$|S_{T'}|\le 6$,  $|S_{T'}\cap
V(A)|=|V(A)|=5$, $|S_T\cap S_{T'}|=1$, and $|S_{T'}\cap
V(B)|=0$. Since $|S_T|=5$ and $G$ is 5-connected, we see that $B\cap
C=\emptyset$ or $B\cap D=\emptyset$. However, this implies $|V(C)|\le
4$ or $|V(D)|\le 4$, contradicting the choice of $(T,S_T,A,B)$ that
$|V(A)|$ is minimum.

We may thus assume $A\cap C\ne \emptyset$. 
Then  $|(S_{T'} \cup S_T) - V(B\cup D) | \geq |S_{T'}|+1$. So 
 $|(S_{T'} \cup S_T) - V(A \cup C)| = |S_T| + |S_{T'}| - |(S_{T'} \cup
 S_T) - V(B \cup D) | \leq 4$. 
Since $G$ is $5$-connected, $B \cap D = \emptyset$.
In addition, $A \cap D \neq \emptyset$; as otherwise, $|V(D)| \leq 4<|V(A)|$, contradicting the choice of $(T,S_T,A,B)$
that $|V(A)|$ is minimum.
Therefore, $|(S_{T'} \cup S_T) - V(B \cup C) | \geq |S_{T'}|+1$. 
Hence, $|(S_{T'} \cup S_T) - V(A \cup D) | = |S_T| + |S_{T'}| - |(S_{T'}
\cup S_T) - V(B \cup C) | \leq 4$.
Since $G$ is $5$-connected, $B \cap C = \emptyset$.
Thus, $|V(B)| \leq |S_{T'} - V(T')| = 4$, contradicting the fact
$|V(A)|\ge 5$ and $|V(A)|$ is minimum. 
\qed

\medskip

The next lemma will allow us to assume that if $(T,S_T,A,B)\in {\cal
  Q}_x$ with $|V(A)|$ minimum and  $(T',S_{T'},C,D)\in {\cal
  Q}_x$  with $T'\cap A\ne \emptyset$ then  $T\cong K_3$ and $T'\cong K_3$. 

\begin{lem}
\label{T_K_3_T'_K_2}
Let $G$ be \textsl{a}  5-connected graph and $x \in V(G)$.
Suppose for any $T\subseteq G$ with $x\in V(T)$ and $T\cong K_2$ or
$T\cong K_3$, $G/T$ is not 5-connected. 
Let $(T,S_T,A,B) \in {\cal Q}_x$ with $|V(A)|$ minimum and $(T',S_{T'},C,D) \in {\cal Q}_x$ with $T' \cap A \neq \emptyset$.
Suppose $T'\cong K_2$. Then one of the following holds:
\begin{itemize}
\item[$(i)$] $G$ contains \textsl {a} $TK_5$ in which  $x$ is not \textsl{a} branch vertex.
\item[$(ii)$] $G$ contains $K_4^-$.
\item[$(iii)$] There exist distinct  $x_1,x_2,x_3\in N(x)$ such that for any distinct $y_1,y_2\in N(x)-\{x_1,x_2,x_3\}$, $G':=G-\{xv:v\notin \{x_1,x_2,x_3,y_1,y_2\}\}$ contains $TK_5$.
\end{itemize}
\end{lem}

\pf By Lemma~\ref{T_K_2}, we may assume $T\cong K_3$. By Lemma~\ref{5cut_triangle}, we may  further assume 
$|S_T| = 6$. Note the symmetry between $C$ and $D$, and assume that $V(T) \subseteq S_T - V(D)$. 
Since $|V(T')|=2$,  $|S_{T'}|=5$.

Suppose $A \cap C \neq \emptyset$. Then $|(S_{T'} \cup S_T) - V(B \cup
D) | \geq 7$; otherwise, $(T, (S_{T'} \cup S_T) - V(B \cup D), A \cap C,
B \cup D) \in \mathcal{Q}_x$ and $0<|V(A \cap C)|  < |V(A)|$, contradicting the choice of $(T,S_T,A,B)$ that $|V(A)|$ is minimum. 
Hence, $|(S_{T'} \cup S_T) - V(A \cup C) | = |S_T| + |S_{T'}| - |(S_{T'}
\cup S_T) - V(B \cup D) | \leq 4$. 
Since $G$ is $5$-connected, $B \cap D = \emptyset$. We may assume 
$A \cap D \neq \emptyset$; otherwise, $|V(D)| \leq 4$ and, by Lemma \ref{cut_4-}, $(ii)$ holds. 
We may also assume $|V(D)|>|V(A)|$; otherwise, $(T',S_{T'},D,C) \in
\mathcal{Q}_x$ and, by Lemma~\ref{T_K_2}, $G$ contains $K_4^-$. Hence, 
 $|V(D) \cap S_T| > |V(A \cap C)|+|V(A) \cap S_{T'}| \geq |V(A) \cap
 S_{T'}| + 1$. Then, since $|S_T|=6$ and $V(T)\subseteq S_T-V(D)$,
 $|V(D)\cap S_T|=3$ and $|V(A)\cap S_{T'}|=1$. Hence, $|(S_{T'} \cup
 S_T) - V(B \cup D )| \le 4 $, a contradiction as $G$ is 5-connected. 

Now assume  $A \cap C = \emptyset$. Then, since $|S_{T'}\cap V(A)|\le 4$,  we may
assume  $A \cap D \neq \emptyset$ by Lemma \ref{cut_4-}. 

Suppose $|(S_{T'} \cup S_T) - V(B \cup C) |=5$. Then, since  $|V(A \cap
D)|  < |V(A)|$, $|V(A\cap D)|=1$; otherwise,
$(T', (S_{T'} \cup S_T) - V(B \cup  C), A \cap D, B \cup C) $ contradicts the choice of $(T,S_T,A,B)$ that $|V(A)|$ is minimum. 
Hence by Lemma \ref{cut_4-}, we may assume  $|V(A) \cap S_{T'}| = 4$;  so $V(B) \cap S_{T'} = V(D) \cap S_T = \emptyset$.
Since $G$ is 5-connected, $B \cap D = \emptyset$. 
So $|V(D)|=1$, a contradiction.

Hence, we may assume $|(S_{T'} \cup S_T) - V(B \cup C) | \geq 6$. Then $S_T \cap V(D) \neq \emptyset $ because $|S_{T'}| = 5$. 
By  Lemma \ref{cut_4-}, we may assume $B \cap C \neq \emptyset$
(otherwise $|V(C)|\le 4$).
Hence, since $G$ is $5$-connected, $|(S_{T'} \cup S_T) - V(A \cup D) | \geq 5$. 
Since $|(S_{T'} \cup S_T )- V(A \cup D) | + |(S_{T'} \cup S_T)- V(B \cup C) | = |S_T| + |S_{T'}| = 11$, 
 $|(S_{T'} \cup S_T) - V(A \cup D) | = 5$. If $|V(B\cap C)|=1$ then,
 since $G$ is 5-connected, $G[T\cup (B\cap C)]\cong K_4^-$. If
 $|V(B\cap C)|\ge 2$ then, since  $V(T)\subseteq (S_{T'}\cup S_T)-V(A\cup D)$, the assertion follows from Lemma~\ref{5cut_triangle}.
\qed

\medskip

 The proofs of the remaining two results in this section  use
 Lemmas~\ref{characterization}, \ref{simplerungs-1} and \ref{simplerungs}.
The following result  will allow us to assume  that if $(T,S_T,A,B)\in {\cal Q}_x$ is chosen to
minimize $|V(A)|$ then $N(x)\cap V(A)\ne \emptyset$, which in turn
will allow us to choose another quadruple at $x$.

\begin{lem}
\label{existence_of_a}
Let $G$ be \textsl{a} 5-connected nonplanar graph and $x\in V(G)$. Suppose for any $H\subseteq G$ with $x\in V(H)$ and $H\cong K_2$ 
or $H\cong K_3$, $G/H$ is not 5-connected. 
Let $(T,S_T,A,B) \in \mathcal{Q}_x$  minimizing $|V(A)|$. Then $N(x)\cap V(A)\ne \emptyset$, or  
one of the following holds: 
\begin{itemize}
\item[$(i)$] $G$ contains \textsl{a} $TK_5$ in which $x$ is not \textsl{a} branch vertex.
\item[$(ii)$] $G$ contains $K_4^-$.
\item[$(iii)$] There exist distinct  $x_1,x_2,x_3\in N(x)$ such that for any distinct $u_1,u_2\in N(x)-\{x_1,x_2,x_3\}$, 
$G':=G-\{xv: v\notin \{x_1,x_2,x_3,u_1,u_2\}\}$ contains $TK_5$.
\end{itemize}
\end{lem}

\pf Suppose $N(x)\cap V(A)=\emptyset$. Then, since $G$ is 5-connected,
$|S_T|=6$ and $T\cong K_3$. Let $V(T)=\{x,x_1,x_2\}$ and 
$S_T=\{x,x_1,x_2,v_1,v_2,v_3\}$. By Lemma~\ref{N(x1)capA}, we may
assume $N(x_1)\subseteq V(A)\cup \{x,x_2\}$, and any three independent
paths in $G_A:=G[A+ (S_T-\{x\})] - E(S_T)$ from $\{x_1,x_2\}$ to $v_1,v_2,v_3$, respectively, with two from
$x_1$ and one from $x_2$, must include a path from $x_2$ to $v_1$. 

We wish to apply Lemma~\ref{characterization}.  
Let $G_A'$ be obtained from $G_A$ by adding  a new vertex  $x_1'$ and joining $x_1'$ to each vertex in $N(x_1)\cap V(G_A)$ with an edge.
Thus, in $G_A'$, $x_1$ and $x_1'$ have the same set of neighbors. Note
that $\{x_1,x_1',x_2\}$ and $\{v_1,v_2,v_3\}$ are independent sets in $G_A'$. 
\medskip

{\it Claim} 1. There is no separation $(A_1,A_2)$ in $G_A'$ such that  $|V(A_1\cap A_2)|\le 3$,
 $\{x_1,x_1',x_2\}\subseteq V(A_1)$ and $\{v_1,v_2,v_3\}\subseteq V(A_2)$.

For, suppose such $(A_1,A_2)$ does exist.  Then $\{x_1,x_1'\}\not\subseteq
V(A_1\cap A_2)$; for, otherwise, $A_1-\{x_1,x_1',x_2\}\ne 
\emptyset$ (as $\{x_1,x_1',x_2\}$ is independent in $G_A'$ and $x_2$
has a neighbor in $V(A)$) and,
hence,   $(V(A_1\cap A_2)-\{x_1'\})\cup \{x,x_2\}$ is a cut in $G$ of size at most 4, a contradiction.

Thus, we may assume by
symmetry that  $x_1\notin V(A_1\cap A_2)$. Then $(A_1,A_2)$ may be chosen so that $x_1'\notin V(A_1\cap A_2)$ (as $x_1'$ has the same set of neighbors as 
$x_1$ in $G_A'$). Moreover, $V(A_1)-V(A_2) \subseteq
\{x_1,x_1',x_2\}$; otherwise $S_T':=V(A_1\cap A_2)\cup V(T)$ is a cut
in $G$ with $|S_T'|\le 6$, and $G-S_T'$ has a component strictly contained in $A$, 
contradicting the choice of $(T,S_T,A,B)$ that $|V(A)|$ is minimum. 

Since $G$ is 5-connected and $N(x_1)\subseteq V(A)\cup \{x,x_2\}$,
$V(A_1\cap A_2)\cup \{x,x_2\}$ is not a 4-cut in $G$. So $x_2\in
V(A_1)-V(A_2)$ and $|V(A_1\cap A_2)|=3$. 
Since $G$ is 5-connected and  $V(A_1)-V(A_2) \subseteq
\{x_1,x_1',x_2\}$, $N(x_1)=\{x,x_2\}\cup V(A_1\cap A_2)$. Since 
$N(x_2)\cap V(A_1)\ne \emptyset$,  there exists $v\in V(A_1\cap A_2)$
such that $vx_2\in E(G)$. Now $G[\{v,x,x_1,x_2\}] \cong K_4^-$ and ($ii$) holds.  
This completes the proof of Claim 1.

\medskip
 
Since any three disjoint paths in $G_A'$ from $\{x_1,x_2,x_1'\}$ to
$\{v_1,v_2,v_3\}$ contains a path from $x_2$ to $v_1$, it follows from 
Claim 1 and Lemma~\ref{characterization} that $G_A'$ has a separation
$(J,L)$ such that $V(J\cap  L)=\{w_0,\ldots, w_n\}$, 
$(J,w_0,\ldots,w_n)$ is 3-planar, $(L,(x_1,x_2,x_1'),(v_2,v_1,v_3))$
is a ladder along some sequence $b_0 \ldots b_m$, 
where $b_0=x_2$, $b_m=v_1$, and $w_0\ldots w_n$ is the reduced
sequence of $b_0 \ldots b_m$. (Note that if $(ii)$ of
Lemma~\ref{characterization} holds then, by Claim 1, 
$(G_A',(x_1,x_2,x_1'),(v_2,v_1,v_3))$ is a rung, and we let $L=G_A'$
and $J$ consist of $v_1$ and $x_2$.)

Since $L$ is a ladder, $L$ contains three disjoint paths $P_1,P_2,P_3$ from $x_1,x_2, x_1'$, respectively, 
to $\{v_1,v_2,v_3\}$, with $v_1\in V(P_2)$. Without loss of generality, 
we may further assume that $v_2\in V(P_1)$ and $v_3\in V(P_3)$. Let
$(R_i,(a_{i-1},b_{i-1},c_{i-1}), (a_i,b_i,c_i))$,  $i\in [m]$, be the rungs in $L$, 
with $a_i\in V(P_1)$, $b_i\in V(P_2)$ and $c_i\in V(P_3)$ for
$i=0,\ldots, m$. Since $G$ is 5-connected, $(J,w_0,\ldots,w_n)$ is
planar and,  by Lemmas~\ref{simplerungs-1} and \ref{simplerungs}, 
we may assume that the rungs in $L$ have the simple structures as in Lemma~\ref{simplerungs}.

\medskip

{\it Claim} 2. There exist $t\in V(A)$ and  independent paths
$Q_1,Q_2,Q_3,Q_4,Q_5$ in $G_A$ such that $Q_1,Q_2,Q_3,Q_4$ are  
from $t$ to $x_1,x_2,v_1,v_2$, respectively, and 
$Q_5$ is from $x_1$ to $v_3$; and there exist $t\in V(A)$ and
independent paths $Q_1',Q_2',Q_3',Q_4',Q_5'$ in $G_A$ such that $Q_1',Q_2',Q_3',Q_4'$ are from 
$t$ to $x_1,x_2,v_1,v_3$, respectively, and
$Q_5'$ is from $x_1$ to $v_2$. 

We may assume that for $i\in [m]$, $(R_i,(a_{i-1},b_{i-1},c_{i-1}),
(a_i,b_i,c_i))$ is not of type ($iv$) as in Lemma~\ref{simplerungs}. 
For, suppose $(R_i,(a_{i-1},b_{i-1},c_{i-1}),
(a_i,b_i,c_i))$ is of type ($iv$) for some $i\in [m]$, and let $v\in V(R_i)- (\{a_{i-1},b_{i-1},c_{i-1}\}\cup
\{a_i,b_i,c_i\})$. Then Claim 2 holds with
$v, va_{i-1}\cup a_{i-1}P_1x_1,vb_{i-1}\cup b_{i-1}P_2x_2, vb_{i}\cup
b_{i}P_2v_1, va_{i}\cup a_{i}P_1v_2, P_3$ as $t,Q_1,Q_2,Q_3,Q_4,Q_5$,
respectively, and with $v, vc_{i-1}\cup c_{i-1}P_3x_1, v b_{i-1}\cup b_{i-1}P_2x_2$, 
$vb_{i}\cup
b_{i}P_2v_1$, $vc_{i}\cup c_{i}P_3v_3, P_1$ as $t,Q_1',Q_2',Q_3',Q_4',Q_5'$,
respectively.

We claim that there exists $q\in [m]$, such that $x_1b_q\in E(G)$. 
Let $q\ge 1$ be the smallest integer such that 
$(R_q,(a_{q-1},b_{q-1},c_{q-1}), (a_q,b_q,c_q))$ is not of type
($ii$) as in Lemma~\ref{simplerungs}, which must exist as $x_1\notin \{v_1,v_2,v_3\}$.
Then  $a_{q-1} = x_1$ and $c_{q-1} = x_1'$. Since $G$ is 5-connected,
$(R_q,(a_{q-1},b_{q-1},c_{q-1}),
(a_q,b_q,c_q))$ cannot be of type ($iii$) (thus, must be of type
$(i)$) as in Lemma~\ref{simplerungs}. Since $x_1$ and $x_1'$ have the same set of neighbors in $G_A'$, 
$a_q\ne x_1$ and $c_q\ne x_1'$. Since $G$ is 5-connected,  $V(R_q)=\{x_1,x_1',a_q,b_q,c_q\}$.
Since $N(x_1)\subseteq V(A)\cup \{x,x_2\}$ and $G$ is 5-connected,
$x_1b_q\in E(G)$. 

 We choose such $q$ to be maximum. Note that $q\ne 0$ as $x_1b_0\notin
 E(G_A')$. We now show the existence of $t$ and $Q_i, i\in [5]$; the proof of the
existence of $t$ and $Q_i',i\in [5]$, is symmetric (by switching the roles of
$v_2,P_1$ and $v_3,P_3$). 

We may assume that for any choice of $P_1,P_3$ there does not exist $r$, with $q<r\le m$, such that $L$
has disjoint paths  $S,S'$ from $b_r,x_1$ to $v_2,v_3$, respectively, 
and internally disjoint from $J\cup P_2$. For, suppose for some choice
of $P_1,P_3$ such $r,S,S'$ exist. 
By Claim 1, $J\cup P_2$ is 2-connected. So let $P_2'$ denote the path
between $x_2$ and $v_1$ in $J\cup P_2$ such that the cycle $P_2'\cup
P_2$ bounds the infinite face of $J\cup P_2$. 
Let $t\in V(P_2')$ such that $x_2t\in E(P_2')$. 
If there exist independent paths $L_1,L_2$ in $J\cup P_2$ from $t$ to $b_q,b_r$, respectively, and internally disjoint from $P_2'$, then 
$L_1\cup b_qx_1,L_2\cup S,tx_2,tP_2'v_1,S'$ give  the desired 
$Q_1,Q_2,Q_3,Q_4,Q_5$, respectively. 
Thus we may assume that such $L_1,L_2$ do not exist. So $J\cup P_2$ has a separation $(J_1,J_2)$ 
such that $|V(J_1 \cap J_2)|\le 3$, $t\in V(J_1)-V(J_2)$, and 
$\{b_q,b_r,v_1,x_2\}\subseteq V(J_2)$. By planarity of $J\cup P_2$,
$V(J_1\cap J_2)$ contains $x_2$ and a vertex $t'\in V(tP_2'v_1)$. 
Since  $V(J_1\cap J_2)$ cannot be a cut in $G$, we must have $|V(J_1\cap J_2)|=3$, $t'=v_1$, and $V(J_1\cap J_2)-\{t',x_2\}\subseteq V(b_rP_2v_1)$. 
Let $b_s\in V(J_1\cap J_2)-\{t',x_2\}$. Then $V(T)\cup \{a_s,b_s,c_s\}$ is a
cut in $G$ separating $\bigcup_{i=1}^sR_s$ from $B+t$, contradicting the choice of  $(T,S_T,A,B)$ that $|V(A)|$
is minimum.

Hence, for any $j>q$, $(R_j,(a_{j-1},b_{j-1},c_{j-1}), (a_j,b_j,c_j))$
must be of type ($i$) or ($ii$) as in Lemma~\ref{simplerungs} and
there is no edge in $G_A'$ from $P_2$ to $P_1-x_1$.  Also
notice that, for $j\le q$ with $b_{j-1}\ne b_q$, because of edges $x_1b_q,x_1'b_q$ in $G_A'$,
$(R_j,(a_{j-1},b_{j-1},c_{j-1}), (a_j,b_j,c_j))$ must be of type $(ii)$
as in Lemma~\ref{simplerungs}.  For $j\le q$ with $b_{j-1}= b_q$, we
see that $V(R_j)=\{x_1,x_1',a_j,b_q,c_j\}$ as $G$ is 5-connected, and
we may assume that $b_qa_j\notin E(G)$ (otherwise,
$b_q,b_qx_1,b_qP_2x_2,b_qP_qv_1,b_qa_q\cup a_qP_1v_2,P_3$ give the desired
$t,Q_1,Q_2,Q_3,Q_4,Q_5$).

Thus, we may assume that for some $j>q$, $\{a_{j-1},c_{j-1}\}\cap
\{a_j,c_j\}=\emptyset$. For, otherwise, $(G_A, x_1,x_2,v_1,v_2,v_3)$
is planar, and the assertion follows from Lemma~\ref{apex}.

If $R_j-a_{j-1}$ contains disjoint paths $S_1,S_2$ from $b_j, c_{j-1}$
to $a_j,c_j$, respectively, then $b_j$ and the paths
$S_1\cup a_jP_1v_2, x_1P_3c_{j-1}\cup S_2\cup c_jP_3v_3$ contradict
the nonexistence of $b_r, S, S'$. 
So assume $S_1,S_2$ do not exist. Then by Lemma~\ref{seymour-2dp}, $(R_j-a_{j-1}, a_j,c_j,b_j,c_{j-1})$ is planar.  
By Lemma~\ref{apex}, 
we may assume $|V(R_j-a_{j-1})|\le 5$. 

If $|V(R_j-a_{j-1})|=5$ then there exists $v\in V(R_j)-\{a_{j-1},a_j,b_j, c_{j-1},c_j\}$ such that $v$ is 
adjacent to all of $\{a_{j-1},a_j,b_j, c_{j-1},c_j\}$; so $b_j$ and the paths $b_jva_j\cup
a_jP_1v_2,P_3$ contradict
the nonexistence of $b_r, S, S'$. 

Hence, we may assume $|V(R_j-a_{j-1})|=4$. 
Then, since $R_j$ has no cut of size at most 3 separating
$\{a_{j-1},b_{j-1},c_{j-1}\}$ from $\{a_j,b_j,c_j\}$, we must have
$a_{j-1}c_j,a_jc_{j-1}\in E(G)$. Note that there exists $t>q$ such
that $L$ has a path $Z$ from $b_t$ to $z\in V(x_1P_1a_{j-1}-x_1)\cup
V(x_1'P_3c_{j-1}-x_1')$ and internally disjoint from $J\cup P_1\cup
P_2\cup P_3$; for
otherwise, $\{a_j,b_j,c_j,x_1\}$ would be a cut in $G$.  If $z\in
V(x_1P_1a_{j-1}-x_1)$ then $b_t$ and the paths $Z\cup zP_1v_2,P_3$ contradict
the nonexistence of $b_r, S, S'$. 
So assume $z\in V(x_1P_3c_{j-1}-x_1)$. Then $b_t$ and the paths  $Z\cup zP_3c_{j-1}\cup
c_{j-1}a_j\cup a_jP_1v_2, x_1P_1a_{j-1}\cup a_{j-1}c_j\cup c_jP_3v_3$ contradict
the nonexistence of $b_r, S, S'$,  with  $x_1'P_3c_{j-1}\cup
c_{j-1}a_j\cup a_jP_1v_2, x_1P_1a_{j-1}\cup a_{j-1}c_j\cup c_jP_3v_3$
as $P_1,P_3$, respectively. This  completes the proof of Claim 2.

\medskip

Now that we have the paths in Claim 2, we turn to $G_B:=G[B+S_T-x_1]$. Choose $x_3\in N(x)\cap V(B)$, let $u_1:=x_3$
and let $u_2\in N(x)-\{x_1,x_2,x_3\}$ be arbitrary. Note that $u_2\in
S_T\cup V(B)$.  We wish to
prove ($iii$) by attempting to find a $TK_5$ in $G':=G-\{xv: v\notin \{u_1,u_2,x_1,x_2\}\}$. 
Since $G$ is $5$-connected and $N(x_1)\cap V(B)=\emptyset$, $G_B$ has four independent paths
$B_1,B_2,B_3,B_4$ from $u_1$ to $v_1,v_2,v_3,x_2$, respectively, and
we may assume that these paths are induced. 

\medskip
{\it Claim} 3. We may assume $u_2\notin S_T$.

For, suppose $u_2\in S_T$. If $u_2=v_1$ then $T\cup Q_1\cup Q_2\cup (Q_3\cup v_1x)\cup u_1x\cup B_4\cup (B_2\cup Q_4)\cup (B_3\cup Q_5)$ 
is a $TK_5$ in $G'$ with branch vertices $t,u_1,x,x_1,x_2$. If $u_2=v_2$ then  
$T\cup Q_1\cup Q_2\cup (Q_4\cup v_2x)\cup u_1x\cup B_4\cup (B_1\cup Q_3)\cup 
(B_3\cup Q_5)$ is a $TK_5$ in $G'$ with branch vertices $t,u_1,x,x_1,x_2$. Now assume $u_2=v_3$. Then 
$T\cup Q_1'\cup Q_2'\cup (Q_4'\cup v_3x)\cup u_1x\cup B_4\cup (B_1\cup Q_3')\cup 
(B_2\cup Q_5')$ is a $TK_5$ in $G'$ with branch vertices $t,u_1,x,x_1,x_2$. This completes the proof of Claim 3. 

\medskip

Let $P$ be a path in $G_B$ from $u_2$ to some $w_2\in V(B_1\cup B_2\cup
B_3\cup B_4)-\{u_1\}$ and internally disjoint from $B_1\cup B_2\cup
B_3\cup B_4$.

\medskip
{\it Claim} 4.  We may assume that for any choice of $P$, $w_2 \in V(B_4)$. 

For, if $w_2 \in V(B_1)$ then $T \cup Q_1 \cup Q_2\cup (Q_3\cup
v_1B_1w_2\cup P\cup u_2x)
\cup u_1x\cup B_4  \cup (B_2\cup Q_4) \cup (B_3\cup Q_5)$ 
is a $TK_5$ in $G'$ with branch vertices $t,u_1,x,x_1,x_2$.  
If $w_2 \in V(B_2)$ then $T\cup  Q_1 \cup Q_2\cup (Q_4\cup
v_2B_2w_2\cup P\cup u_2x)
\cup u_1x\cup B_4  \cup (B_1\cup Q_3) \cup (B_3\cup Q_5)$ 
is a $TK_5$ in $G'$ with branch vertices $t,u_1,x,x_1,x_2$. 
If $w_2\in V(B_3)$ then  $T\cup Q_1' \cup Q_2'\cup (Q_4'\cup
v_3B_3w_2\cup P\cup u_2x)\cup u_1x\cup  B_4 \cup (B_1\cup Q_3') \cup (B_2\cup Q_5')$ 
is a $TK_5$ in $G'$  with branch vertices $t,u_1,x,x_1,x_2$.
This completes the proof of Claim 4.

\medskip

Let $U_2$ denote the $(B_1\cup B_2\cup B_3)$-bridge of $G_B$
containing $B_4+u_2$. That is, $U_2$ is the subgraph of $G_B$ induced by the edges in the
component of $G_B-(B_1\cup B_2\cup B_3)$ containing $B_4+u_2$
and the edges from that component to $B_1\cup B_2\cup B_3$.

\medskip
{\it Claim} 5. We may assume that  $V(U_2)\cap V(B_2\cup B_3)=\{u_1\}$. 

For, suppose there exists $w\in V(U_2)\cap  V(B_2\cup B_3)$ such that
$w\ne u_1$.
By symmetry, we may assume $w\in V(B_2-u_1)$ and choose $w$ so that $wB_2v_2$ is minimal. 

Then $U_2$ has a path $X$ between $x_2$ to $w$ and internally disjoint from $B_1\cup B_2\cup B_3$, 
and a path from $u_2$ to some $u_2'\in V(X)$ and internally disjoint
from $X\cup B_1\cup B_2\cup B_3$.
Since $G$ is 5-connected,  $U_2$ has four independent paths from
$u_2'$ to four distinct vertices in $V(U_2)\cap V(B_1\cup B_2\cup B_3)$ and
internally disjoint from $B_1\cup B_2\cup B_3$.  Thus, by Lemma~\ref{Perfect}, $U_2$ contains 
independent paths $L_1,L_2,L_3,L_4$ from $u_2'$ to $u_2,x_2,w,w'$, 
respectively, and internally disjoint from $B_1\cup B_2\cup B_3$, where $w'\in V(B_1\cup B_2\cup B_3)$. 

If $w'\in V(wB_2u_1-w)$ then  $T\cup (L_1\cup u_2x)\cup L_2\cup
(L_3\cup wB_2v_2\cup P_1)\cup (u_1B_2w'\cup L_4)\cup u_1x\cup (B_1\cup P_2)\cup (B_3\cup P_3)$ 
is a $TK_5$ in $G'$ with branch vertices $u_1,u_2',x,x_1,x_2$. (Note
we identify $x_1'$ with $x_1$ when we use $P_3$.)

If $w'\in V(B_1-u_1)$ then  
$T\cup Q_1'\cup Q_2'\cup (Q_4'\cup B_3\cup u_1x)\cup (L_1\cup u_2x)\cup L_2\cup (L_3\cup wB_2v_2\cup Q_5')\cup
(L_4\cup w'B_1v_1\cup Q_3')$
is a $TK_5$ in $G'$ with branch vertices $t,u_2',x,x_1,x_2$. 

If  $w'\in V(B_3-u_1)$ then  
$T\cup Q_1\cup Q_2\cup (Q_3\cup B_1\cup u_1x)\cup (L_1\cup u_2x)\cup L_2\cup (L_3\cup wB_2v_2\cup Q_4)\cup (L_4\cup w'B_3v_3\cup Q_5)$
is a $TK_5$ in $G'$ with branch vertices $t,u_2',x,x_1,x_2$. 
This completes the proof of Claim 5.

\medskip

Now let $z\in V(B_1\cap U_2)$ such that $zB_1v_1$ is minimal. Since $G$ is 5-connected, there exists a path $Y$ in $G_B-x$ 
from some $y\in V(zB_1u_1)-\{u_1,z\}$
to some $y'\in V(B_2\cup B_3)-\{u_1\}$ and internally disjoint from $U_2\cup B_1\cup B_2\cup B_3$. 

\medskip

{\it Claim} 6. We may assume that  $G[U_2 -B_1+z]$ has no independent paths from $u_2$ to $x_2,z$, respectively. 

For, suppose $G[U_2 -B_1+z]$ (and hence $G[U_2\cup zB_1u_1]$) has independent paths from $u_2$ to $x_2,z$,
respectively. Then by Lemma~\ref{Perfect}, $G[U_2\cup zB_1u_1]$ has independent paths $L_1,L_2,L_3,L_4$ from $u_2$ to distinct
vertices  $x_2,z,z_1,z_2$, respectively, and 
internally disjoint from $B_1$, where $u_1,z_2,z_1,z$ occur on $B_1$ in the order listed. Possibly, $u_1=z_2$. 

If $y'\in V(B_2-u_1)$ then $T\cup Q_1'\cup Q_2'\cup (Q_4'\cup B_3\cup
u_1x)\cup u_2x \cup L_1\cup (L_2\cup zB_1v_1\cup Q_3')\cup (L_3\cup z_1B_1y\cup Y\cup y'B_2v_2\cup Q_5')$ 
is a $TK_5$ in $G'$ with branch vertices $t, u_2,x,x_1,x_2$. 

If  $y'\in V(B_3-u_1)$ then 
$T\cup Q_1\cup Q_2\cup (Q_4\cup B_2\cup u_1x) \cup u_2x\cup L_1\cup (L_2\cup zB_1v_1\cup Q_3)\cup (L_3\cup z_1B_1y\cup Y\cup y'B_3v_3\cup Q_5)$ 
is a $TK_5$ in $G'$ with branch vertices $t, u_2,x,x_1,x_2$. 

\medskip

By Claim 6, $G[U_2-B_1 +z]$ has a 1-separation $(U_{21},U_{22})$ such that $u_2\in V(U_{21})-V(U_{22})$ and $\{x_2,z\}\subseteq V(U_{22})$. 
We choose this separation so that $U_{22}$ is minimal. 
Let $u_2'$ denote the unique vertex in $V(U_{21}\cap U_{22})$. 
By the minimality of $U_{22}$, we see that $U_{22}$ has  independent paths $L_1,L_2$ from $u_2'$ to $x_2,z$, respectively. 

\medskip

{\it Claim} 7.  We may assume that  $u_2'$ has exactly  two neighbors in $U_{22}$. 

For, otherwise, by the minimality of $U_{22}$, $G[U_{22}\cup zB_1u_1]-u_1$ has three independent paths from $u_2'$ to 
three distinct vertices in $V(zB_1u_1-u_1)\cup \{x_2\}$. So by Lemma~\ref{Perfect}, $G[U_{22}\cup zB_1u_1]-u_1$ has 
independent paths $L_1',L_2',L_3'$ from $u_2'$ to $x_2,z,z_1$, respectively, and 
internally disjoint from $B_1$, where 
$z,z_1,u_1$ occur on $B_1$ in order. Let $L$ be a path in $U_{21}$
from $u_2$ to $u_2'$.

If $y'\in V(B_2-u_1)$ then $T\cup Q_1'\cup Q_2'\cup (Q_4'\cup B_3\cup
u_1x)\cup (L\cup u_2x)\cup L_1'\cup (L_2'\cup zB_2v_1\cup Q_3')\cup (L_3'\cup z_1B_1y\cup Y\cup y'B_2v_2\cup Q_5')$ 
is a $TK_5$ in $G'$ with branch vertices $t, u_2',x,x_1,x_2$.

If  $y'\in V(B_3-u_1)$ then 
$T\cup Q_1\cup Q_2\cup (Q_4\cup B_2\cup u_1x) \cup (L\cup u_2x)\cup L_1'\cup (L_2'\cup zB_2v_1\cup Q_3)\cup (L_3'\cup z_1B_1y\cup Y\cup y'B_3v_3\cup Q_5)$ 
is a $TK_5$ in $G'$ with branch vertices $t, u_2',x,x_1,x_2$. This
completes the proof of Claim 7. 

\medskip

Since $G$ is $5$-connected, it follows from Claim 7 that $u_2'$ has at least two neighbors in
$U_{21}$. Since all paths from $u_2$ to $B_1\cup B_2\cup B_3\cup B_4$ must end on $B_4$, 
$G[U_{21}\cup zB_1u_1]-\{z,u_1\}$ has independent paths $L_3,L_4$ from $u_2'$ to $z_1, u_2$, respectively, 
and internally disjoint from $B_1$, where $z_1\in
V(zB_1u_1)-\{z,u_1\}$.  

If $y'\in V(B_2-u_1)$ then $T\cup Q_1'\cup Q_2'\cup (Q_4'\cup B_3\cup
u_1x)\cup (L_4\cup u_2x)\cup L_1\cup (L_2\cup zB_2v_1\cup Q_3')\cup (L_3\cup z_1B_1y\cup Y\cup y'B_2v_2\cup Q_5')$ 
is a $TK_5$ in $G'$ with branch vertices $t, u_2',x,x_1,x_2$.

If  $y'\in V(B_3-u_1)$ then 
$T\cup Q_1\cup Q_2\cup (Q_4\cup B_2\cup u_1x) \cup (L_4\cup u_2x)\cup L_1\cup (L_2\cup zB_2v_1\cup Q_3)\cup (L_3\cup z_1B_1y\cup Y\cup y'B_3v_3\cup Q_5)$ 
is a $TK_5$ in $G'$ with branch vertices $t, u_2',x,x_1,x_2$.\qed

\medskip

We conclude this section with another technical lemma, which deals with a
special case that occurs in the proof of Lemma~\ref{T_K_3_T'_K_3_second}. It is included in this section because its proof
also makes use of Lemmas~\ref{characterization}, \ref{simplerungs-1} and \ref{simplerungs}.

\begin{lem}\label{lemma_butterfly}
Let $G$ be \textsl{a} $5$-connected nonplanar graph and $x\in V(G)$.
Let $(T,S_T,A,B) \in  \mathcal{Q}_x$ such that  $|V(A)|$ is
minimum, and suppose there exists $(T',S_{T'},C,D) \in
\mathcal{Q}_x$ such that $T'\cong K_3$,  $T'\cap A\ne \emptyset$,
$V(A \cap C) = S_T \cap V(C) = V(B \cap D) = V(B) \cap S_{T'} =
\emptyset$, $|V(A) \cap S_{T'}| = |V(D) \cap S_{T}| =|V(D\cap T)|= 1$, and 
$|S_T \cap S_{T'}| =5$.  Suppose
for any $H\subseteq G$ with $x\in V(H)$ and $H\cong K_2$ or $H\cong
K_3$, we have $G/H$ is not 5-connected,  $|V(H\cap A)| \le 1$, and $H\cong
K_3$ when $H\cap A\ne \emptyset$. Then one of the following holds: 
\begin{itemize}
\item[$(i)$] $G$ has \textsl{a}  $TK_5$ in which $x$ is not \textsl{a} branch vertex.
\item[$(ii)$] $G$ contains $K_4^-$.
\item[$(iii)$] There exist $x_1,x_2,x_3\in N(x)$ such that, 
for any distinct $y_1,y_2\in N(x)-\{x_1,x_2,x_3\}$, $G':=G-\{xv:v\notin \{x_1,x_2,x_3,y_1,y_2\}\}$ contains $TK_5$.
\end{itemize}
\end{lem}

\pf  Note that $|S_T|=|S_T\cap S_{T'}|+|V(D\cap T)|=6$.  Let $V(T) = \{x,w,x_1\}$ and $T' = \{x,a,b\}$ such that 
$V(A) \cap S_{T'} = \{a\}$ and $V(D) \cap S_{T} = \{w\}$, and let $S_T \cap S_{T'} = \{x,x_1,b,z_1,z_2\}$. 
Then $|V(D)| = |V(A)|=|V(A\cap D)|+1$. Moreover, 
\begin{itemize}
\item [(1)] $|N(s)\cap V(A)| \geq 2$ for $s \in \{b,z_1,z_2\}$,
\end{itemize}
for,  otherwise, $(T,(S_T-\{s\})\cup (N(s)\cap V(A)),A-N(s),G[B+s]) \in \mathcal{Q}_x$,
contradicting the choice of $(T,S_T,A,B)$ that $|V(A)|$ is minimum.
We may assume that 
\begin{itemize}
\item [(2)] $G$ has no edge from $T-x$ to $T'-x$,
\end{itemize}
 as otherwise $G[T\cup T']$ contains $K_4^-$ and $(ii)$ holds. We may
 also assume 
\begin{itemize}
\item [(3)] $N(x_1)\cap V(D)\ne \{w\}$ and $N(w)\cap V(A)\ne \emptyset$,
\end{itemize}
for, otherwise, let $S:=S_T\setminus \{x_1\}$ and $B'=G[B+x_1]$ if
$N(x_1)\cap V(D)=\{w\}$, and let  $S:=S_T\setminus \{w\}$ and
$B'=G[B+w]$ if $N(w)\cap V(A)=\emptyset$;  then $(xw,S,A,B')\in {\cal
  Q}_x$, and $(ii)$ follows from Lemma~\ref{T_K_2}.  We may further assume that 
\begin{itemize}
\item [(4)]  for any $x' \in N(x) \cap V(A \cap D)$, $xx'z_1x$
 or $xx'z_2x$ is a triangle. 
\end{itemize}
For, let $x' \in N(x) \cap V(A \cap  D)$. By Lemma~\ref{existence_of_T'}, 
we may assume
that there exists $H\subseteq G$ with $x, x'\in V(H) $ and $H\cong K_2$ or $H\cong 
K_3$. By the assumption of this lemma,  $H\cong K_3$ and $V(H)\cap S_T\ne \{x\}$. If
$V(H)\cap \{b,x_1\}\ne\emptyset$ then  $H\cup T$ or
$H\cup T'$ contains $K_4^-$. So we may assume $V(H)\cap \{z_1,z_2\}\ne
\emptyset$ and, hence, $xx'z_1x$ or $xx'z_2x$ is a triangle.

\medskip 

We may  assume that 
\begin{itemize}
\item [(5)] $|N(x)\cap V(A\cap D)|\le 2$.
\end{itemize}
 For, otherwise, by (4), there
exist $i\in [2]$ and distinct $x',x''\in N(x)\cap V(A\cap D)\cap
N(z_i)$. So $G[x',x'',x,z_i]$ contains $K_4^-$, and $(ii)$ holds.  

\medskip

We now distinguish two cases.

\medskip

\medskip

{\it Case} 1. $z_i\notin N(x)$ for $i\in [2]$.

Then  by (4),  $N(x) \cap V(A \cap D) = \emptyset$. We prove that $(iii)$
holds with $x_2=w$ and $x_3=b$. Let $y_1,y_2\in N(x)-\{x_1,x_2,x_3\}$. 
Since $G$ is 5-connected and  $z_1,z_2\notin N(x)$, 
we may assume $y_1 \in V(B\cap C)$. Then   
$G_B:=G[B+\{b, x_1,z_1,z_2\}]$ has independent paths $Y_1,Y_2,Y_3,Y_4$ from $y_1$ to $z_1,z_2,x_1,b$, respectively. 

We may  assume that $wz_i\notin E(G)$ for $i\in [2]$. For, suppose $wz_1\in E(G)$. If  $G[A+\{b,w,x_1\}]$ has independent paths $Q_1,Q_2$ from 
$b$ to $x_1,w$, respectively, then $T\cup bx\cup Q_1\cup Q_2\cup y_1x\cup (Y_1\cup z_1w)\cup Y_3\cup Y_4$ is a $TK_5$ in $G'$ with 
branch vertices  $b,w,x,x_1,y_1$. So we may assume that such $Q_1,Q_2$ do  not exist.  
Then  $G[A+\{b,w,x_1\}]$ has a cut vertex $v$ separating $b$ from
$\{w,x_1\}$. Let $D$ denote the component of $G[A+\{b,w,x_1\}]-v$
containing $b$. Since  $|N(b)\cap V(A)|\ge 2$ (by (1)), $|V(D)|\ge 2$. Now
$\{b,v,x,z_1,z_2\}$ is a cut in $G$, and $G$ has a separation
$(G_1,G_2)$ such that  $V(G_1\cap G_2)=\{b,v,x,z_1,z_2\}$,
 $|V(G_1)|\ge 6$ and $\{a,b\}\subseteq V(G_1)$,
and $B+\{w,x_1\}\subseteq G_2$. By the choice of $(T,S_T,A,B)$ with
$|V(A)|$ minimum, $|V(G_1)|=6$. Let $u\in V(G_1)-V(G_2)$. If  $u=a$ then,
$V(G_1\cap G_2)\subseteq N(a)$ (since $G$ is 5-connected) and  $bv\in
E(G)$ (since $|N(b)\cap
V(A)|\ge 2$); so  $G[\{a,b,v,x\}]\cong K_4^-$, and $(ii)$ holds. So
assume $u\ne a$. Then $v=a$ and $G[\{b,u,v,x\}]$ contains
$K_4^-$; so $(ii)$ holds. 

We may assume that $G_A:=G[A+\{b, w,x_1,z_1,z_2\}]$ does not contain three
independent paths, with one from $x_1$ to $b$, one from $b$ to $w$, 
and one from $w$ to $z_i$ for some $i\in [2]$. 
For, otherwise, such three paths and $T\cup bx\cup y_1x \cup Y_i\cup Y_3 \cup Y_4$ form a $TK_5$ in $G'$ with 
branch vertices $b,w,x,x_1,y_1$.
 
We wish to apply Lemma~\ref{characterization}. 
Let $G_A'$ be the graph obtained from $G_A$ by
identifying $z_1$ and $z_2$ as $z'$, 
and duplicating $w,b$ with $w',b'$, respectively (adding edges from $w'$ 
to all vertices in $N(w)$, and from  $b'$ to all vertices in $N(b)$). Then any three disjoint paths in $G_A'$ from $\{w,x_1,w'\}$ to 
$\{b,z',b'\}$, if exist, must contain a path from $x_1$ to $z'$. 

Suppose $G_A'$ has a separation $(A_1,A_2)$ such that $|V(A_1\cap A_2)|\le 2$, $\{w,x_1,w'\}\subseteq V(A_1)$, and $\{b,z',b'\}\subseteq V(A_2)$. 
Since $w$ and $w'$ have the same set of neighbors in $G_A'$, we may assume  $\{w,w'\}\subseteq V(A_1\cap A_2)$ or   $\{w,w'\}\cap V(A_1\cap A_2)=\emptyset$. 
If $\{w,w'\}\subseteq V(A_1\cap A_2)$ then $V(A_1)=\{x_1\}\cup
V(A_1\cap A_2)$ as $\{x,x_1,w\}$ cannot be a cut in $G$; hence, $N(x_1)\cap V(D)=\{w\}$, contradicting (3).
So $\{w,w'\}\cap V(A_1\cap A_2)=\emptyset$. Suppose $\{b,b,',z'\}\cap
V(A_1\cap A_2)=\emptyset$. Then, since
$wz_i\notin E(G)$ for $i\in [2]$, $V(A_1\cap A_2)\cup
\{x_1,x\}$ is a cut in $G$ separating $w$ from  $B+\{b,z_1,z_2\}$,
contradicting the fact that $G$ is 5-connected. So  $\{b,b,',z'\}\cap
V(A_1\cap A_2)\ne \emptyset$. Note that $\{b,b'\}\not\subseteq  V(A_1\cap A_2)$;
as otherwise $\{b,x,x_1\}$ would be a cut in $G$. Thus, we may assume
that $b,b'\notin V(A_1\cap A_2)$ as $b$ and $b'$ have the same set of
neighbors in $G_A'$. Hence, $z'\in V(A_1\cap A_2)$. Now
$S:=\{x,x_1,z_1,z_2\}\cup (V(A_1\cap A_2)-\{z'\})$ is a cut in $G$
separating $w$ from $B+b$. Since $G$ is 5-connected, $x_1\notin
V(A_1\cap A_2)$. If
$|V(A_1-x_1-A_2)|\ge 2$ then $(xx_1,S, A_1-x_1-A_2, G-S-A_1)\in {\cal
  Q}_x$ which contradicts the choice of $(T,S_T,A,B)$ with $|V(A)|$
minimum.  So $V(A_1-x_1-A_2)=\{w\}$. Since $G$ is 5-connected, $wz_i\in
E(G)$ for $i\in [2]$, a contradiction.

Hence, by Lemma~\ref{characterization}, $G_A'$ has a separation $(J,L)$ such that $V(J\cap
L)=\{w_0,\ldots, w_n\}$, $(J,w_0,\ldots,w_n)$ is planar (since $G$ is 5-connected), 
$(L,(w,x_1,w'),(b,z',b'))$ is a ladder along a sequence $b_0\ldots
b_m$, where $b_0=x_1$, $b_m=z'$, and 
$w_0\ldots w_n$ is the reduced sequence of $b_0 \ldots b_m$.  
Moreover, we may assume that 
$L$ has disjoint induced paths $P_1,P_2,P_3$ from $w,x_1,w'$ to $b,z',b'$, respectively, and $J$ is a connected plane graph with $P_2$ as part of the
outer walk of $J$ and $w_0,\ldots, w_n$ occurring on $P_2$ in
order. (When $(ii)$ of Lemma~\ref{characterization} holds, we let
$J=P_2$.) Note that by Lemmas~\ref{simplerungs-1} and
\ref{simplerungs}, each rung of $(L,(w,x_1,w'),(b,z',b'))$
is  of type $(i)$--$(iv)$ as in Lemma~\ref{simplerungs}, with possible 
exceptions of those rungs containing $z'$.  Let
$(R_j,(a_{j-1},b_{j-1},c_{j-1}),(a_j,b_j,c_j))$, $j\in [m]$,  be the  rungs in 
$(L,(w,x_1,w'),(b,z',b'))$ such that $a_j\in V(P_1)$ and $c_j\in
V(P_3)$ for $j=0, 1, \ldots, m$.

We now show that there exists $t\in N(w)$ such that $t\in
V(P_2)-\{x_1,z'\}$.  
For, suppose such $t$ does not exist.
Choose the largest $j$ such that $\{w,w'\}\subseteq V(R_j)$ and
$(R_j,(a_{j-1},b_{j-1},c_{j-1})$, $(a_j,b_j,c_j))$ is not of type
$(ii)$ in Lemma~\ref{simplerungs}, which is
well defined as $w\ne b$. 
Since $G$ is 5-connected and $w$ and $w'$ have the same set
of neighbors in $G_A'$, $(R_j,(a_{j-1},b_{j-1},c_{j-1}), (a_j,b_j,c_j))$ cannot be of
type $(iii)$ as in Lemma~\ref{simplerungs}. 
Moreover,  $(R_j,(a_{j-1},b_{j-1},c_{j-1}),(a_j,b_j,c_j))$ is not of type $(iv)$ as 
in Lemma~\ref{simplerungs}, as otherwise $G$ contains $K_4^-$
(obtained from $R_j-\{b_{j-1},b_j\}$ after identifying $w$ with
$w'$). So $(R_j,(a_{j-1},b_{j-1},c_{j-1}),(a_j,b_j,c_j))$ is of type $(i)$ as 
in Lemma~\ref{simplerungs}.
Now $V(R_j) = \{a_{j},b_{j},c_{j},w,w'\} $, as otherwise
$\{a_j,b_j,c_j,w\}$ would be  a cut in $G$. Then $wb_j\in E(G)$; for
otherwise, $N(w) \subseteq \{a_j,c_j,x,x_1\}$, a contradiction. 
Hence $t:=b_j$ is as desired.

Without loss of generality, we may assume that the edge of $P_2$
incident with $z'$ corresponds to the edge of $G$ incident with $z_1$.
We view $P_3$ as a path in $G_A$ from $b$ to $w$.  
Then  $G_A-V(P_1\cup P_3)-z_2$ has independent paths from $t$ to $x_1, z_{1}$, respectively. 
Hence, by Lemma~\ref{Perfect}, $G_A$ has five independent paths
$Q_1,Q_2,Q_3,Q_4,Q_5$ from $t$ to $x_1,w,z_1,(V(P_1\cup P_3)-\{w\})\cup \{z_{2}\}$, 
respectively, with only $t$ in common, and internally disjoint from
$P_1\cup P_3$. Without loss of generality, we may assume that $Q_4$ ends at
$t'\in V(P_3)$. 

If $G_B-x$ contains disjoint paths $S_1,S_2$ from $z_1,b$ to
$y_1,x_1$, respectively, then $T\cup bx\cup P_1\cup S_2\cup Q_1\cup Q_2\cup (Q_3\cup S_1\cup y_1x)\cup 
(Q_4\cup t'P_3b)$ is $TK_5$ in $G'$ with branch vertices $b,t,w,x,x_1$. Hence, we may assume such $S_1,S_2$ 
do not exist. Then by Lemma~\ref{seymour-2dp}, there exists a collection ${\cal D}$ of subsets of
$(G_B-x)-\{z_1,b,y_1,x_1\}$  such that $(G_B-x,{\cal
  D},z_1,b,y_1,x_1)$ is 3-planar. 

If $(G_B-x,\{b,x_1,z_1,z_2\})$ is planar then the assertion of the
lemma follows from
Lemma~\ref{apex}, with the cut $\{b,x,x_1,z_1,z_2\}$ giving the required 
5-separation for Lemma~\ref{apex}. 

So we may assume that either ${\cal D}=\emptyset$ and
$z_2$ does not belong to the facial  walk of $G_B-x$ containing
$\{b,x_1,y_1,z_1\}$, or ${\cal D}=\{D\}$ for some $D\subseteq
V(G_B-x)-\{b,x_1,y_1,z_1\}$ and $z_2\in D$. Thus, since
$G$ is 5-connected and $(G_B-x, \{z_1,b,y_1,x_1\})$ is 3-planar, $G_B-x$
has disjoint paths $S_1',S_2'$ from $z_2,b$ to $y_1,x_1$,
respectively. 
Moreover, if $b$
has degree at least two in $G_B-x$ then $G_B-x$ has independent paths $Y,Y_2',Y_3',Y_4'$, with $Y$ from $b$ to $x_1$ and $Y_2',Y_3',Y_4'$ from 
$y_1$ to $z_2,x_1,b$, respectively. 

We may assume that  $G_A'-J$ contains a path $Z$ from $z_2$ to
some $z_2'\in V(P_1\cup P_3)-\{b,b'\}$ and internally disjoint from
$P_1\cup P_3$. For, suppose not. Then, since
$|N(z_2)\cap V(A)|\ge 2$ (by (1)), $z_2$ has at least two neighbors in $J-z'$. Then  $G_A'-V(P_1\cup P_3)-z_1$
has independent paths from $t$ to $x_1, z_{2}$, respectively;
for otherwise,  $G_A'-V(P_1\cup P_3)-z_1$ has a cut vertex $v\in V(tP_2x_2)$ separating $t$ from
$\{x_1,z_2\}$ and, hence,  $V(T)\cup \{v,z_1,z_2\}$ is a cut in $G$,
contradicting the choice of $S_T$ with $|V(A)|$ minimum. 
Hence, by Lemma~\ref{Perfect}, $G_A$ has five independent paths
$Q_1',Q_2',Q_3',Q_4',Q_5'$ from $t$ to $x_1,w,z_2,(V(P_1\cup P_3)-\{w,b'\})\cup \{z_{1}\}$, 
respectively, with only $t$ in common, and internally disjoint from
$P_1\cup (P_3-\{b',w'\})$. Without loss of generality, we may assume that $Q_4'$ ends at
$t''\in V(P_3)$. Then $T\cup bx\cup P_1\cup S_2'\cup Q_1'\cup Q_2'\cup (Q_3'\cup S_1'\cup y_1x)\cup 
(Q_4'\cup t''P_3b)$ is $TK_5$ in $G'$ with branch vertices $b,t,w,x,x_1$.

 Without loss of generality, we may assume that $z_2'\in
V(P_3)$. We may further assume that $b$ has only one neighbor in
$G_B-x$; for, otherwise, $T\cup bx\cup
P_1\cup Y\cup y_1x\cup (Y_2'\cup Z\cup z_2'P_3w)\cup Y_3'\cup Y_4'$ is a
$TK_5$ in $G'$ with branch vertices $b,w,x,x_1,y_1$. 

Thus, since $G$ is 5-connected and $bw\notin E(G)$ (by (2)), $b$ has a neighbor $u\in V(A)-V(P_1\cup
P_3)$. We choose $u$ and the rung
$(R_j,(a_{j-1},b_{j-1},c_{j-1}),(a_j,b_j,c_j))$ such that $b,b',u\in
V(R_j)$. Since $b$ and $b'$ have the same set of neighbors
in $G_A'$, $a_{j-1}=b$ if, and only if, $c_{j-1}=b'$. Moreover,  we must have
$b_j=z'$ because of the path $Z$.

First, suppose $b_{j-1}=z'$. Then $a_{j-1}\ne b$ and $c_{j-1}\ne
b'$. If $z_2$ has no neighbor in  $ V(G_A'-J-R_j)$ then $V(T)\cup \{a_{j-1},c_{j-1},z_1\}$ is a cut in $G$
separating $a_{j-1}P_1w\cup c_{j-1}P_3w\cup (J-z')$ from $B\cup
(R_j-\{b',z'\})$, contradicting the choice of $(T,S_T,A,B)$ with
$|V(A)|$ minimum. 
Thus, $z_2$ has a neighbor in $V(G_A'-J-R_j)$; so the above path $Z$ may be chosen to be disjoint from $R_j$.
Let $S$ be a path in $R_j-\{a_{j-1},c_{j-1}\}$ from $b$ to
$z_1$ (which must exist as otherwise $\{a_{j-1},c_{j-1},z_2\}\cup
V(T)$ is a cut in $G$ contradicting the choice of $(T,S_T,A,B)$ that
$|V(A)|$ is minimum). So $T\cup bx\cup P_1\cup (S\cup z_1P_2x_1) \cup y_1x\cup
(Y_2\cup Z\cup z_2'P_3w)\cup Y_3\cup Y_4$ is $TK_5$ in $G'$ with branch
vertices $b,w,x,x_1,y_1$.

Now assume $b_{j-1}\ne z'=b_j$.   Since $b_{j-1}\ne
b_j$ and since $b$ and $b'$ have the same set of neighbors in $G_A'$, we must have $a_{j-1}=b$ and $c_{j-1}=b'$. If $u\in
\{b_{j-1},b_j\}$ then, since $bz'\notin E(G_A')$, $u=b_{j-1}$; and
let $S=bb_{j-1}$. Now suppose $u\notin \{b_{j-1},b_j\}$. Then 
$\{b,b_{j-1},x,z_1,z_2\}$ is a cut in $G$ separating $u$ from
$(J-z')\cup B$. By the choice of $(T,S_T,A,B)$ that $|V(A)|$ is
minimum, $\{u\}=V(R_j)-\{b,b',b_{j-1},z'\}$. Since $G$ is 5-connected, 
$N(u)=\{b,b_{i-1},x,z_1,z_2\}$. Let $S=bub_{j-1}$. Since $|N(z_2)\cap V(A)|\ge 2$ (by (1)), the path $Z$ may be chosen to be 
disjoint from $R_j$. 
So $T\cup bx\cup P_1\cup (S\cup b_{j-1}P_2x_1) \cup y_1x\cup
(Y_2\cup Z\cup z_2'P_3w)\cup Y_3\cup Y_4$ is $TK_5$ in $G'$ with branch vertices $b,w,x,x_1,y_1$.

\medskip
{\it Case} 2.  $N(x) \cap \{z_1,z_2\}\ne \emptyset$. 

Without loss of generality, we may assume $xz_1 \in E(G)$. We may
further assume $z_1$ is not adjacent to any of $\{a,b,w,x_1\}$; for
otherwise, $G[T+z_1]$ or $G[T'+z_1]$ contains $K_4^-$, and  $(ii)$ holds. 
We wish to prove $(iii)$, with $x_2=b$ and $x_3=z_1$. Let  $y_1,y_2\in N(x)-\{b,x_1,z_1\}$ be distinct. 

\medskip

{\it Subcase} 2.1. There exists some $i\in [2]$ such that $y_i\in V(B)\cup \{z_2\}$. 

Without loss of generality, assume $y_1\in V(B)\cup \{z_2\}$ and, whenever
possible,  let $y_1\in V(B)$. Let $G_B:=G[B+\{b,x_1,z_1,z_2\}]$.  When
$y_1\in V(B)$ let $t=y_1$  and let $Y_1,Y_2,Y_3,Y_4,Y_5$ be
independent paths in $G[B]$ from $t$ to $z_1,y_1,b,x_1,z_2$, respectively.
When $y_1=z_2$ let $t=y_1$ and let  $Y_1,Y_2,Y_3,Y_4$ be independent
paths in $G[B]$ from $t$ to $z_1,y_1,b,x_1,z_2$, respectively. Let $G_A=G[A+\{b,w,x_1,z_1\}]$.

We may assume that there is no cycle in $G_A$ containing $\{b,x_1,z_1\}$. For, such a cycle and
$xb\cup xx_1\cup xz_1\cup Y_1\cup (Y_2\cup y_1x)\cup Y_3\cup Y_4$ is a $TK_5$ in $G'$ 
with branch vertices $b,t,x,x_1,z_1$. 

We may also assume that $G_A$ is 2-connected. To see this, we first  assume $N(x_1)\cap N(w)=\{x\}$; for otherwise, letting $u\in
(N(x_1)\cap N(w))-\{x\}$  we see that $G[T+u]$ contains
$K_4^-$ and $(ii)$ holds. Therefore, since $N(w)\cap V(A)\ne \emptyset\ne N(x_1)\cap V(A)$ (by (3)), it suffices to show that $G[A+\{b,z_1\}]$ is 2-connected. 
So assume for a contradiction that   there exists a separation $(A_1,A_2)$ in $G[A+\{b,z_1\}]$ such that $|V(A_1 \cap A_2)| \leq 1$. 
Without loss of generality, let $|\{b,z_1\} \cap V(A_1)| \leq 1$. Then $V(A_1)\not\subseteq V(A_2)\cup \{b,z_1\}$ as  
$|N(s)\cap V(A)| \geq 2$ for $s \in \{b,z_1\}$ (by (1)). 
Hence, $V(T)\cup (\{b,z_1\} \cap A_1) \cup V(A_1 \cap A_2) \cup 
\{z_2\}$ is a cut in $G$ of size at most $6$ which 
separates $A_1$ from the rest of $G$, contradicting the choice of
$(T,S_T,A,B)$ that $|V(A)|$ is minimum.

Then, since $G_A$ has no cycle containing $\{b,x_1,z_1\}$,  $(i)$, or $(ii)$, or $(iii)$ of Lemma~\ref{Watkins} holds for $G_A$ and $\{b,x_1,z_1\}$. 
So for each $u\in \{b,x_1,z_1\}$, $G_A$ has a 2-cut $S_u$ separating $u$ from $\{b,x_1,z_1\}-\{u\}$, and let $D_u$ denote a
union of components of $G_A-S_u$  such that $u\in V(D_u)$ for $u\in \{b,x_1,z_1\}$ and
$D_b, D_{x_1}, D_{z_1}$ are pairwise disjoint.   We choose $S_u$ and $D_u$, $u\in
\{b,x_1,z_1\}$, to maximize $D_b\cup D_{x_1} \cup D_{z_1}$.  Note that,
since $wx_1\in E(G)$, $w\notin V(D_b\cup D_{z_1})$. 
 
We claim that for $u \in \{ b,x_1,z_1 \}$, $V(D_u) = \{u\}$. For,
otherwise, $S:=S_u \cup \{u, x,z_2\}$ is a cut in $G$ separating $D_u-u$ from the rest of $G$. 
If $|V(D_u)|\ge 3$ then $(ux,S,D_u,G-S-D_u)\in {\cal Q}_x$ contradicts the choice of $(T,S_T,A,B)$ with $|V(A)|$ minimum. So let 
$V(D_u)=\{u,u'\}$ and let $S_u=\{s_u,t_u\}$. Since $G$ is 5-connected, $N(u')=\{s_u,t_u,u,x,z_2\}$. Since $|N(u)\cap V(A+w)| \geq 2$ 
(by (1) and (3)), we may assume that $us_u \in E(G)$. 
Then $G[\{s_u,u,u',x\}]$ contains $K_4^-$, and  $(ii)$ holds.

For $u \in \{ b,x_1,z_1 \}$, let $S_u = \{s_{u}, t_{u}\}$. 
Since $G_A$ is 2-connected, $\{us_{u}, ut_{u}\} \subseteq E(G)$.
Note $a\in \{s_b,t_b\}$; so we may assume $s_{b}t_{b} \not\in E(G)$
because otherwise $G[\{x,b,s_{b},t_{b}\}]$ contains $K_4^-$, and $(ii)$ holds.
Similarly, $w\in \{s_{x_1},t_{x_1}\}$ and we may assume
$s_{x_1}t_{x_1} \not\in E(G)$.
If $(i)$ of Lemma~\ref{Watkins} occurs then $ax_1 \in E(G)$,
contradicting (2). 
If $(iii)$ of Lemma~\ref{Watkins} occurs then let $R_1,R_2$ be the
components of $G_A- V(D_b\cup D_{x_1} \cup D_{z_1})$ and assume
without loss of generality that $s_{u} \in V(R_1)$ and $t_u\in V(R_2)$ for $u \in \{
b,x_1,z_1 \}$. By symmetry, assume $w\notin V(R_1)$. 
Hence,  $(xb, \{x,b,x_1,s_{z_1},z_2\}, R_1-s_{z_1}, G-R_1-\{x,b,x_1,z_2\}]) \in \mathcal{Q}_x$ 
with $2\le |V(R_1-s_{z_1})| < |V(A)|$, contradicting the choice of $(T,S_T,A,B)$. 

So we may assume that $(ii)$ of Lemma~\ref{Watkins} holds. Without
loss of generality let $R_1,R_2$ be the
components of $G- V(D_b\cup D_{x_1} \cup D_{z_1})$ containing $z =s_b = s_{x_1} = s_{z_1}$,
$\{t_b,t_{x_1},t_{z_1}\}$, respectively. By (2), $z\ne a$ and $z\ne w$. 
So $a =  t_{b}$ and $w =
t_{x_1}$.  Thus, we may assume $xz\notin E(G)$ as, otherwise,
$G[T+z]$ contains $K_4^-$ and $(ii)$ holds. Hence, $R_1=R_2$ (otherwise $z$ would have
degree at most 4 in $G$). 
By (1) and by the maximality of $D_b\cup D_{x_1} \cup D_{z_1}$, 
$G[R_2+z_2]$ is 2-connected (since $G$ is 5-connected).

We claim that there exist distinct  $t_1,t_2 \in \{a,w,t_{z_1}\}$ such
that $G[R_2+z_2]$ contains disjoint paths $P_1,P_2$ from $z,t_1$ to $z_2,t_2$, respectively. 
 For, suppose $\{a,w\}$ cannot serve as $\{t_1,t_2\}$. Then, by Lemma~\ref{seymour-2dp}, 
$(G[R_2+z_2],a,z_2,w,z)$ is 3-planar. Hence, $G[R_2+z_2]$ has disjoint paths 
from $z,a$ to $z_2,t_{z_1}$, respectively, or disjoint paths from $z,w$ to $z_2,t_{z_1}$, respectively. 

Suppose $z_2\ne y_1$.
Recall the definition of $t$ and the paths $Y_1,Y_2,Y_3,Y_4, Y_5$. If
$\{t_1,t_2\} = \{a,w\}$ then $bxx_1zb \cup xz_1z
\cup ( x_1w \cup P_2 \cup ab )\cup (Y_2\cup y_1x)\cup (Y_5 \cup P_1) \cup Y_3 \cup Y_4$ 
is a $TK_5$ in $G'$ with branch vertices $b,t,x,x_1,z$.
If $\{t_1,t_2\} = \{a,t_{z_1}\}$ then $bxz_1zb \cup xx_1z \cup
( z_1t_{z_1} \cup P_2 \cup ab ) 
\cup Y_1 \cup (Y_2 \cup y_1x) \cup Y_3\cup (Y_5 \cup P_1)$ is a $TK_5$ in $G'$ with branch vertices $b,t,x,z,z_1$.
If $\{t_1,t_2\} = \{w,t_{z_1}\}$ then $x_1xz_1zx_1\cup xbz\cup ( x_1w \cup
P_2 \cup t_{z_1}z_1 )  \cup Y_1 \cup (Y_2 \cup y_1x)\cup Y_4\cup (Y_5 \cup P_1)$ is a $TK_5$ in $G'$ with branch vertices $t,x,x_1,z,z_1$.

So assume $z_2=y_1$. Then $y_2\ne z_2$; and hence, by the choice of
$y_1$, we have $y_2\in V(A)\cup \{w\}$. 
If $R_2-z$ has  independent paths $S_1,S_2,S_3$ from $y_2$ to
$a,w,t_{z_1}$, respectively, then
$xbzx_1x\cup y_2x\cup (S_1\cup ab)\cup (S_2\cup wx_1)\cup Y_3\cup
Y_4\cup (Y_1\cup z_1t_{z_1}\cup 
S_3)\cup (Y_2\cup z_2x)$ is  a $TK_5$ in $G'$ with branch vertices
$b,t,x,x_1,y_2$. So assume such $S_1,S_2,S_3$ do not exist. Then $R_2$
has a separation $(A_1,A_2)$ such that $z\in V(A_1\cap A_2)$,
$|V(A_1\cap A_2)|\le 3$, $y_2\in V(A_1-A_2)$ and
$\{a,w,t_{z_1}\}\subseteq 
V(A_2)$. Thus $S:=\{x,z_2\}\cup V(A_1\cap A_2)$ is a 5-cut in $G$
separating $y_2$ from $B\cup A_2\cup \{b,x_1,z_1,z\}$. Hence, by the
choice of $(T,S_T,A,B)$ (with $|V(A)|$ minimum),
$V(A_1-A_2)=\{y_2\}$. Therefore, since $G$ is 5-connected,
$N(y_2)=S$. By the maximality of $D_b\cup D_{x_1}\cup D_{z_1}$, 
$R_2-\{y_2,z\}$ has a path $Q$ from $a$ to $w$. Then $bxx_1zb\cup
(ba\cup Q \cup wx_1)\cup zy_2x\cup (Y_1\cup z_1z)\cup (Y_2\cup z_2x)\cup Y_3\cup
Y_4$ is a $TK_5$ in $G'$ with branch vertices $b,t,x,x_1,z$.

\medskip
{\it Subcase} 2.2. $y_1,y_2\in V(A)\cup \{w\}$.

First, we show that we may assume $y_1=w$. 
For, suppose $y_1,y_2\in V(A)$. Then by Lemma~\ref{existence_of_T'},
for each $i\in [2]$ there exists $(T_i, S_{T_i},A_i,B_i)\in {\cal Q}_x$ such that
$x,y_i\in V(T_i)$ and $T_i\cong K_2$ or $T_i\cong K_3$.  By the
assumption of this lemma,
we have $T_i\cong K_3$ and $V(A)\cap S_{T_i}=\{y_i\}$. Hence,
$\{b,w,x_1,z_1,z_2\}\cap V(T_i)\ne \emptyset$ for $i\in [2]$. 
Without loss of generality, we may assume  that $y_1\ne a$. By the
symmetry between $z_1$ and $z_2$, we may
also assume  $z_1\in V(T_1)$; for, otherwise, $G[T+y_1]$ or $G[T'+y_1]$ contains $K_4^-$ and $(ii)$ holds. 
Therefore, we may choose $S_{T_1}=V(T_1)\cup \{b,x_1,z_2\}$.  Note the
symmetry between $T_1,S_{T_1}$ and $T,S_T$, and we may
choose $T_1,S_{T_1}$ as $T,S_T$, respectively. So we may  assume
$y_1=w$ (as $y_1$ now plays the role of $w$).

Let $t \in V(B)$, and let $L_1,L_2,L_3,L_4$ be independent paths in $G_B=G[B+\{b,x_1,z_1,z_2\}]$ from $t$ to $z_1,z_2, b,x_1$, respectively.
 Let $G_A:=G[A+\{b,w,x_1,z_2\}]$. Note that, by the same argument as in Subcase 2.1 (with $z_2$
in place of $z_1$), we may assume that $G_A$ is 2-connected.

We may assume that $G_A$ does not contain  independent paths  from
$z_2,w,b$ to $w,b,x_1$, respectively; for otherwise, 
these paths and  $T\cup bx\cup (L_1\cup z_1x)\cup L_2\cup L_3\cup L_4$ form a $TK_5$ in $G$ 
with branch vertices $b,t,w,x,x_1$.  

Hence, since $G_A$ is 2-connected,  $wz_2\notin E(G)$. We may assume that $wz_1\notin E(G)$; else $G[T+z_1]$ 
contains $K_4^-$ and $(ii)$ holds. Therefore, since $G$ is 5-connected, it follows from (2) that $$|N(w)\cap V(A\cap D)|\ge 3.$$ 
Let $G_A'$ be the graph obtained from $G_A$ by duplicating $w,b$ with $w',b'$, respectively, and adding all edges from $w'$ 
to $N(w)$, and from $b'$ to $N(b)$. 
Then any three disjoint paths in $G_A'$ 
from $\{b,b',z_2\}$ to $\{w,w',x_1\}$ must have a path from $z_2$ to $x_1$,  and we wish to apply Lemma~\ref{characterization}. 

First, we note that $G_A'$ has no cut of size at most 2 separating $\{x_1,w,w'\}$ from $\{b,b',z_2\}$. For, otherwise, $G_A'$ 
has a separation $(A_1,A_2)$ such that $|V(A_1\cap A_2)|\le 2$, 
$\{x_1,w,w'\}\subseteq V(A_1)$ and $\{b,b',z_2\}\subseteq
V(A_2)$. Note that $V(A_1\cap A_2)\ne \{w,w'\}$ as otherwise, $w$ would be a cut vertex in $G_A$. Further,
$\{w,w'\}\cap V(A_1\cap A_2)=\emptyset$; for, otherwise, since $w$ and
$w'$ have the same set of neighbors in $G_A'$, it follows from (3) that $V(A_1\cap
A_2)-\{w,w'\}$ would be a cut in $G_A$ of size at most one.
On the other hand, $V(A_1-A_2) \subseteq \{x_1,w\}$; otherwise
$(T,V(T)\cup \{z_1\}\cup  V(A_1 \cap A_2), (A_1-A_2)-w',G-(T\cup
A_1)) \in \mathcal{Q}_x$ with $1 \leq |(A_1-A_2)-w'| < |A|$,
contradicting the choice of $(T,S_T,A,B)$. 
However, this implies $|N(w) \cap V(A \cap D)| \le |V(A_1\cap A_2)|\le 2$, a contradiction.

Hence by Lemma~\ref{characterization},
$G_A'$ has a separation $(J,L)$ such that $V(J\cap L)=\{w_0,\ldots, w_n\}$, 
$(J,w_0,\ldots,w_n)$ is 3-planar,  $(L,(w,x_1,w'),(b,z_2,b'))$ is a ladder along some sequence $b_0\ldots
b_m$, where $b_0=z_2$, $b_m=x_1$, and $w_0\ldots w_n$ is the reduced sequence of $b_0 \ldots b_m$.
Let $P_1,P_2,P_3$ be three disjoint paths  in $L$ from $w,x_1,w'$ to
$b,z_2,b'$, respectively, and assume that they are induced in $G_A'$.  (Let $L=G_A'$ and
$J=P_2$ if $(ii)$ of Lemma~\ref{characterization} holds.) Let
$(R_i,(a_{i-1},b_{i-1},c_{i-1})$, $(a_i,b_i,c_i))$, $i\in [m]$, be the
rungs in $L$ with $a_i\in V(P_1)$ and $c_i\in V(P_3)$ for $i=0, 1,
\ldots, m$.

Since $|N(w)\cap V(A\cap D)|\ge 3$ and $P_1,P_3$ are induced paths in
$G_A'$, there exists  $w^*\in (N(w)\cap
V(A))-V(P_1\cup P_3)$. We show that there exists $u\in V(P_2)$ such
that $G[G_A+\{x,z_1\}]$  has  five independent paths $Q_1,Q_2,Q_3,Q_4,Q_5$ 
from $u$ to distinct vertices $x_1,w,z_2,u_1,u_2$, respectively, with $u_1,u_2\in V(P_1-w)\cup V(P_3-\{b',w'\})\cup \{x,z_1\}$, 
and internally disjoint from $P_1\cup (P_3-\{b',w'\})$.  
If $w^*\in V(P_2)$ then let $u=w^*$ and we see that  there exist
independent paths in $G_A-(V(P_1-w)\cup V(P_3-\{b',w'\}))$ from $u$ to
$x_1,w,z_2$, respectively; then the paths $Q_1,\ldots, Q_5$ exist by 
 Lemma~\ref{Perfect}. Now suppose $w^*\notin
V(P_2)$. Let $(R_i,(a_{i-1},b_{i-1},c_{i-1}),(w,b_i,w'))$ 
be the rung in $L$ containing $\{w,w',w^*\}$. 
Since $w$ and $w'$ have the same set of neighbors in $G_A'$, $w=a_{i-1}$ iff $w'=c_{i-1}$. If $w=a_{i-1}$ and 
$w'=c_{i-1}$ then $S_T^*:=V(T)\cup \{b_{i-1},b_i,z_1\}$ is a cut in $G$ of size at most 6, and $G-S_T^*$ has a component of size smaller than $|V(A)|$, 
contradicting the choice of $(T,S_T,A,B)$. So $w\ne a_{i-1}$ and $w'\ne c_{i-1}$. 
Suppose $R_i-x_1$ has a separation
$(R',R'')$ such that  $|V(R'\cap R'')|\le 2$, $w\in V(R'-R'')$, and $\{a_{i-1},c_{i-1},b_{i-1},b_i\}-\{x_1\} \subseteq V(R'')$. 
Then we may assume $w'\in V(R'- R'')$ as $w$ and $w'$ have the
same set of neighbors in $G_A'$. Therefore, since $|N(w)\cap V(A\cap D)|\ge
3$, $S_T^*:=V(T)\cup V(R'\cap R'')\cup \{z_1\}$  is a cut in $G$ of size at
most 6, and $G-S_T^*$ has a component of size smaller than $|V(A)|$,
contradicting the choice of $(T,S_T,A,B)$. 
Thus we may assume, by Lemma~\ref{Perfect},  $R_i-x_1$ contains three independent paths from 
$w$ to $a_{i-1},c_{i-1}, \{b_{i-1},b_i\} - \{x_1\}$, respectively, and
internally disjoint from $\{b_{i-1},b_i\}$.  
Again since $w$ and $w'$ have the same set of neighbors in $G_A'$,  the parts of $P_1,P_3$ inside $R$ can be modified 
so that the three paths in $R_i$ correspond to $wP_1a_{i-1}, w'P_3c_{i-1}$ and a 
path from $w$ to some $u\in \{b_{i-1},b_i\}-\{x_1\}$ and internally disjoint from $P_1\cup P_2\cup P_3$. 
Thus, there exist independent paths in $G_A-(V(P_1-w)\cup V(P_3-\{b',w'\}))$
from $u$ to $x_1,w,z_2$, respectively. 
Now the paths $Q_1,\ldots, Q_5$ exist by Lemma~\ref{Perfect},.

We may assume $u_1 = z_1$ and $u_2 = x$. For, otherwise, we may assume
by symmetry that $u_1\in V(P_1)$. 
If $G_B-x$ has disjoint paths $B_1,B_2$ from $z_1,b$ to $z_2,x_1$,
respectively, then  $T\cup bx \cup P_3 \cup B_2 \cup Q_1 \cup Q_2 \cup (Q_3 \cup B_1 \cup z_1x) \cup (Q_4 \cup u_1P_1b)$ 
is a $TK_5$ in $G$ with branch vertices $b,u,w,x,x_1$. (Here we view $P_3$ as a path in $G$ by identifying $b',w'$ with $b,w$, respectively.)
So we may
assume that such $B_1,B_2$ do not exist. Then by Lemma \ref{seymour-2dp},
$(G_B-x, z_1,b,z_2,x_1)$ is planar; so the assertion of the lemma
follows from Lemma~\ref{apex}.

We may also assume $|N(b)\cap V(B)|\le 1$. For, suppose $|N(b)\cap
V(B)|\ge 2$. Then, since $G$ is $5$-connected, $G[B+\{b,x_1,z_2\}]$ contains
independent paths $B_1,B_2$ from $b$ to $x_1,z_2$, respectively.
Hence, $T\cup bx\cup P_3 \cup B_1 \cup Q_1 \cup Q_2 \cup (Q_3 \cup B_2) \cup (Q_4 \cup z_1x)$ 
is a $TK_5$ in $G$ with branch vertices $b,u,w,x,x_1$, where we view $P_3$ as a path in $G'$ by identifying $b',w'$ with $b,w$, respectively. 

Then we may assume $|N(b)\cap V(A+z_2)|\ge 3$ as otherwise, $bz_1\in E(G)$ by (2); so $G[T'+z_1]$ contains $K_4^-$ and $(ii)$ holds.  
Let $b^*\in (N(b)\cap V(A+z_2))-V(P_1\cup P_3)$. 

If $b^*\in V(P_2)$ let $z=b^*$ and let $P=bz$ which is internally
disjoint from $P_1\cup P_2\cup  P_3$.
Now suppose $b^*\notin V(P_2)$. Let $(R_j,(b,b_{j-1},b'),(a_j,b_j,c_j))$ be the rung in $L$ containing 
$\{b,b,'b^*\}$. Since $b$ and $b'$ have the same set of neighbors in $G_A'$,
$b=a_j$ iff $b'=c_j$. If $b=a_j$ and $b'=c_j$ then, since $az_1\notin
E(G)$,  $S_T^*:=V(T')\cup \{b_{j-1},b_j,z_1\}$ is a cut in $G$ of size 6 and $G-S_T^*$ has a component of size smaller than $|V(A)|$, contradicting the choice of $(T,S_T,A,B)$. 
So $b\ne a_j$ and $b'\ne c_j$. 
We claim that $P_1\cap R_j$ and $P_3\cap R_j$ may be modified so that $G_A$ contains a path $P$ from $b$ to some $z\in V(P_2)$
and internally disjoint from $P_1\cup P_2\cup (P_3-\{b',w'\})$. 
If $R_j$ contains three independent paths from $b$ to $a_j,c_j,
\{b_{j-1},b_j\}$, respectively, and internally disjoint from
$\{a_j,c_j,b_{j-1},b_j\}$, then $P_1\cap R_j,P_3\cap R_j$ can be modified so that the three paths in $R_j$ correspond to
$bP_1a_j, b'P_3c_j$ and a path $P$ from $b$ to $z\in \{b_{j-1},b_j\}$
and internally disjoint from $P_1\cup P_2\cup (P_3-\{b',w'\})$.  So assume that such three paths
in $R_j$ do not exist.  Then by the existence of $bP_1a_j$ and $b'P_3c_j$ and by Lemma~\ref{Perfect}, $R_j$ has no
three independent paths from $b$ to $\{a_j,c_j,b_{j-1},b_j\}$ and
internally disjoint from $\{a_j,c_j,b_{j-1},b_j\}$. Thus  
$R_j$ has a separation $(A_1,A_2)$ with $|V(A_1\cap A_2)|\le 2$,   $V(A_1\cap A_2)\subseteq V(P_1\cup P_3)$, $b,b^*\in V(A_1-A_2)$ and 
$\{a_j,c_j,b_{j-1},b_j\}\subseteq V(A_2)$.  Since
$b'$ is a copy of $b$,  we may assume  $b'\in V(A_1-A_2)$. Now, since
$az_1\notin E(G)$, $V(A_1\cap A_2)\cup \{x,b,z_1\}$ 
is a cut in $G$; so  $V(A_1)=V(A_1\cap A_2)\cup \{b,b', b^*\}$
by the choice of  $(T,S_T,A,B)$ that $|V(A)|$ is minimum. 
Then $b^*x,b^*z_1\in E(G)$ (as $G$ is 5-connected); so
$G[\{x,b^*,b,z_1\}]$ contains $K_4^-$, and 
$(ii)$ holds.

Suppose $R_i \neq R_j$.
Since $G$ is $5$-connected, $G[B+\{b,x_1\}]$ has a path $B_1$ from $b$
to $x_1$. Since $Q_3$ is internally disjoint from $P_1\cup P_3$,  we
may assume that $z\in V(Q_3)$ and $P$ is also internally disjoint from $Q_3$. 
Hence, $T\cup bx \cup P_3 \cup B_1 \cup Q_1 \cup Q_2 \cup (uQ_3z \cup P) \cup 
(Q_4 \cup z_1x)$ is a $TK_5$ in $G'$ with branch vertices
$b,u,w,x,x_1$, where we view $P_3$ as a path in $G$ by identifying $b',w'$ with $b,w$, respectively.

So $R_i=R_j$. Then $a_{i-1}=b$ and $c_{i-1} = b'$.  Recall $bw\notin
E(G)$ (by (2)).
Since $w$ and $w' $ (respectively, $b$ and $b'$) have the same set of
neighbors in $G_A'$, it follows from Lemma~\ref{simplerungs} that  $b_{i-1} = b_i$. Then $\{b,b_i,w,x,z_1\}$
is a cut in $G$ separating $P_1\cup (P_3-\{b',w'\})$ from $B\cup
J$. Since $bw\notin E(G)$, $|V(P_1\cup (P_3-\{b',w'\}))|\ge 2$. This
contradicts the choice of $(T,S_T,A,B)$ that $|V(A)|$ is
minimum. \qed

\section{Interactions between  quadruples}

In this section, we explore the structure of $G$ by considering  a quadruple
$(T,S_T,A,B)$ with $|V(A)|$ minimum and a quadruple $(T',S_{T'},C,D)\in {\cal Q}_x$ with $T'\cap A\ne \emptyset$. 
The lemma  below allows us to assume  that if $T\cap C=
\emptyset$ then $A\cap C=\emptyset$. 

\begin{lem}
\label{T_K_3_T'_K_3_first}
Let $G$ be \textsl{a} 5-connected nonplanar graph and $x \in V(G)$.
Suppose for any $H\subseteq G$ with $x\in V(H)$ and $H\cong K_2$ or
$H\cong K_3$, $G/H$ is not 5-connected. 
Let $(T,S_T,A,B) \in {\cal Q}_x$ with $|V(A)|$ minimum and $(T',S_{T'},C,D) \in {\cal Q}_x$ with $T' \cap A \neq \emptyset$.
Suppose  $T\cap C=\emptyset$. Then $A\cap C=\emptyset$, or one of the following holds: 
\begin{itemize}
\item[$(i)$] $G$ contains \textsl{a} $TK_5$ in which $x$ is not \textsl{a} branch vertex.
\item[$(ii)$] $G$ contains $K_4^-$.
\item[$(iii)$] There exist $x_1,x_2,x_3\in N(x)$ such that for any $y_1,y_2\in N(x)-\{x_1,x_2,x_3\}$, 
$G-\{xv:v\notin \{x_1,x_2,x_3,y_1,y_2\}\}$ contains  $TK_5$. 
\end{itemize}
\end{lem}

\pf We may assume $T\cong K_3$ (by Lemma~\ref{T_K_2}) and $T'\cong K_3$
(by Lemma~\ref{T_K_3_T'_K_2}). 
Suppose $A \cap C \neq \emptyset $.  

Then $|(S_{T} \cup S_{T'}) - V(B\cup D) | \geq 7$; otherwise $(T', (S_{T'}
\cup S_T )- V(B \cup D), A \cap C, B \cup D) 
\in \mathcal{Q}_x$ and $1 \leq |V(A \cap C)| \leq |V(A - a)| < |V(A)|$, 
contradicting the choice of $(T,S_T,A,B)$ that $|V(A)|$ is minimum.  Hence $|(S_{T} \cup S_{T'}) - V(A\cup C) |=5$, as $|S_T|=|S_{T'}|=6$. 
Since $T\cap C=\emptyset$, $V(T)\subseteq (S_{T} \cup S_{T'} )- V(A\cup C)$.
 
Suppose $|V(B\cap D)|\ge 2$. Then $G$ has a separation
$(G_1,G_2)$ such that $VG_1\cap G_2)= (S_{T} \cup S_{T'})- V(A\cup C)$
and $|V(G_i)|\ge 7$. So the assertion of this lemma follows from 
Lemma~\ref{5cut_triangle}. 

Hence,  we may assume $|V(B\cap D)|\le 1$.
Therefore, by the minimality of $|V(A)|$, $|S_T\cap V(D)|\ge | S_{T'}\cap V(A)|$. But this implies that 
$|S_T|\ge |(S_{T} \cup S_{T'}) - V(B\cup D) |\ge 7$, a contradiction. \qed

\medskip 
 
We need a lemma for finding paths to deal with a special case when $A\cap C=\emptyset$ for quadruples $(T,S_T,A,B),(T',S_{T'},C,D)\in {\cal Q}_x$.

\begin{lem}
\label{2independentPaths-b-x1x2}
Let $G$ be \textsl{a} 5-connected nonplanar graph and $x \in V(G)$, and suppose for any $H\subseteq G$ with $x\in V(H)$ and $H\cong K_2$ or
$H\cong K_3$, $G/H$ is not 5-connected. 
Let $(T,S_T,A,B) \in {\cal Q}_x$ with $|V(A)|$ minimum and $(T',S_{T'},C,D) \in {\cal Q}_x$ with $T' \cap A \neq \emptyset$.
Let $V(T) =\{x,x_1,x_2\}$ and $V(T')=\{x,a,b\}$ with $a\in V(A)$. Suppose $A \cap C = \emptyset$, $|S_T| = 6=|S_{T'}|$, 
$V(T) \subseteq S_T - V(C)$, 
 $|(S_T \cup S_{T'}) -V( B\cup C)|=7$, and $(S_T \cup S_{T'}) -V( B\cup C
 \cup T \cup T') = \{x_3,x_4\}$. Then $G$ contains $K_4^-$, or the
 following statements hold:
\begin{itemize}
\item [$(i)$] $N(b)\cap V(A-a)\ne \emptyset$ and if $t \in
  N(b)\cap V(A -a)$ then $G[(A -a)+ \{b,x_1,x_2,x_3,x_4\}]$ has
  independent paths from $t$ to $b,x_1,x_2,x_3,x_4$, respectively, and 
\item [$(ii)$] if $b\in S_{T}$ then $G[A+\{b,x_1,x_2\}]$  has independent
  paths from $b$ to $x_1,x_2$, respectively.
\end{itemize}
\end{lem}

\pf 
First, we note that  $N(b) \cap V(A -a) \ne 
\emptyset$. For,  otherwise, $(T,(S_T\cup S_{T'})-V(B\cup C)-\{b\},A-a, G[B\cup C+b])\in {\cal Q}_x$. 
By the choice of $(T,S_T,A,B)$ that $|V(A)|$ is minimum, we must have
$V(A-a)=\emptyset$. So $G$ contains $K_4^-$ by Lemma~\ref{cut_4-}.

To complete the proof of  $(i)$, let $t \in N(b) \cap
V(A -a)$. If $G[(A -a) + \{x_1,x_2,x_3,x_4\}]$  has four independent
paths from $t$ to $x_1,x_2,x_3,x_4$, respectively, 
then these four paths and $tb$ give the desired five paths. 
So we may assume that such four paths do not exist. 
Then $G[(A -a) + \{x_1,x_2,x_3,x_4\}]$ has a separation $(G_1,G_2)$ such that $|V(G_1 \cap G_2)| \le 3$, $t \in V(G_1-G_2)$ and $\{x_1,x_2,x_3,x_4\} \subseteq V(G_2)$. 
Hence, $(T',V(T')\cup V(G_1\cap G_2), G_1 - G_2, G-T'-G_1)\in {\cal Q}_x$ and $1\le |V(G_1 - G_2)| \leq |V(A -a)|<|V(A)|$, contradicting the choice of $(T,S_T,A,B)$.

To prove $(ii)$, let $b\in S_T$ and  assume that the two paths in
$(ii)$ do not exist. Note that if $b\in V(T)$ then $T\cup T'$ contains $K_4^-$. So we may assume $b\notin V(T)$. 
 Then, $G[A+ \{b,x_1,x_2\}]$ has a separation $(G_1,G_2)$ 
such that $|V(G_1) \cap V(G_2)|\le 1$, $b \in V(G_1)-V(G_2)$ and
$\{x_1,x_2\} \subseteq V(G_2)$. Since $N(b)\cap V(A-a)\ne \emptyset$
and $|V(G_1) \cap V(G_2)|\le 1$, $|V(G_1-G_2)|\ge 2$. 
Let $S_{bx} = (S_T - \{x_1,x_2\}) \cup V(G_1\cap G_2)$, and let
$F=G_1-S_{bx}$. Then $|V(F)|\ge 1$ as $|V(G_1-G_2)|\ge 2$.    If $|V(F)|\ge 2$ then $(bx, S_{bx}, F, G-S_{bx}-F) 
\in \mathcal{Q}_x$ with $2\le |V(F)|< |V(A)|$, contradicting the choice of $(T,S_T,A,B)$ 
that $|V(A)|$ is minimum. So assume $|V(F)|=1$ and let $v\in V(F)$. Since
$G$ is 5-connected, $v$ is adjacent to all vertices in $S_{bx}$. If
$v\ne a$ then $V(G_1\cap G_2)=\{a\}$; so $G[\{a,b,v,x\}]$ contains
$K_4^-$. Now assume $v=a$. Let $w\in V(G_1\cap G_2)$. Since
$N(b)\cap V(A-a)\ne \emptyset$, $bw\in E(G)$. So  $G[\{a,b,w,x\}]$ contains $K_4^-$.
\qed

\medskip

In the next two lemmas, we consider the case when quadruples $(T,S_T,A,B)$ and $(T',S_{T'},C,D)$ may be chosen so that $|V(T'\cap A)|=2$.

\begin{lem}\label{reduction}
Let $G$ be \textsl{a} 5-connected nonplanar graph and $x \in V(G)$.
Suppose for any $H\subseteq G$ with $x\in V(H)$ and $H\cong K_2$ or
$H\cong K_3$, $G/H$ is not 5-connected. 
Let $(T,S_T,A,B) \in {\cal Q}_x$ with $|V(A)|$ minimum.
Suppose there exists  $(T',S_{T'},C,D)\in {\cal Q}_x$ such that $T'\cong K_3$ and $|V(T'\cap A)|=
2$. Then  one of the following holds:
\begin{itemize}
\item [$(i)$] $G$ contains \textsl{a} $TK_5$ in which $x$ is not \textsl{a} branch vertex. 
\item [$(ii)$] $G$ contains $K_4^-$. 
\item [$(iii)$] There exist $x_1,x_2,x_3\in N(x)$ such that for any $y_1,y_2\in N(x)-\{x_1,x_2,x_3\}$, 
$G-\{xv:v\notin \{x_1,x_2,x_3,y_1,y_2\}\}$ contains  $TK_5$. 
\item [$(iv)$]  $|S_T\cap S_{T'}|=1$, $|S_{T'}\cap V(B)|= 2$, and either $|S_T\cap
V(C)|=2$ and $T\cap C=\emptyset$ or $|S_{T}\cap V(D)|=2$ and $T\cap
D=\emptyset$.
\end{itemize}
\end{lem}

\pf We may assume $T\cong K_3$ (by Lemma~\ref{T_K_2}).  We may also assume that
$|S_T|=|S_{T'}|=6$; for, otherwise, $(i)$ or $(ii)$ or $(iii)$ follows
from Lemma~\ref{5cut_triangle}. We may further assume $|V(A)|\ge 5$; as otherwise, by  Lemma~\ref{cut_4-}, $G$ contains $K_4^-$ 
and $(ii)$ holds.
 
Let $T'=\{a,b,x\}$ with $a,b\in V(A)$.  By symmetry, assume $T\cap
C=\emptyset$. Then,  by Lemma~\ref{T_K_3_T'_K_3_first}, we may assume $A\cap
C=\emptyset$. Now $B\cap C\ne \emptyset$; for, otherwise, $|V(C)|=|S_T\cap
V(C)|\le 3$, contradicting the choice of $(T,S_T,A,B)$ that $|V(A)|$ is minimum.  
Hence, $S_T\cap V(C)\ne \emptyset$ as $S_{T'}-\{a,b\}$ is not a cut in
$G$. Moreover, $A \cap D \neq \emptyset$; for otherwise, $|V(A)\cap
S_{T'}|=5$ and, hence, $|S_{T'}\cap S_T|=1$ and $|S_{T'}\cap V(B)|=0$;
so $(S_T\cup S_{T'})-V(A\cup D)$ is a cut 
in $G$ of size at most 4 and separating $B\cap C$ from $A\cup D$, a contradiction.

We claim that $|(S_T \cup S_{T'}) - V(B \cup C)| = 7$ and $|(S_T \cup
S_{T'}) -V(A \cup D)| = 5$. 
First, note that  $|(S_T \cup S_{T'}) - V(B \cup C)| \ge 7$; 
otherwise, 
 $(T' , (S_T \cup S_{T'}) - V(B \cup C), A \cap D, B \cup C) \in
 \mathcal{Q}_x$ and $1\le |V(A \cap D)| \leq |V(A - a)| < |V(A)|$, contradicting 
the choice of $(T,S_T,A,B)$ that $|V(A)$ is minimum. Also note that $|(S_T \cup S_{T'}) -V(A
\cup D)| \geq 5$ since $B\cap C\ne \emptyset$ and $G$ is $5$-connected. 
Thus the claim follows from the fact that $|(S_T \cup S_{T'}) - V(B \cup
C)| + |(S_T \cup S_{T'}) -V(A \cup D)|=|S_T| + |S_{T'}| =12$.

 We may assume that $|S_T\cap V(C)|\ne 1$ or $|S_{T'}\cap V(A)|\ne 2$. For, suppose  
$S_T\cap V(C)=\{c\}$ and $S_{T'}\cap V(A)=\{a,b\}$.  If $a,b\in N(c)$ then $G[T'+c]$ contains $K_4^-$ and $(ii)$
holds. So by the symmetry between $a$ and $b$, we may assume that $ca
\notin E(G)$. Then $(T, (S_T-c)\cup \{b\}, A-b, G[B+c])\in {\cal Q}_x$, contradicting the choice of $(T,S_T, A, B)$ that $|V(A)|$
is minimum. 

We may also assume $T\cap D\ne\emptyset$; for, otherwise, since $A\cap D\ne
\emptyset$, $(i)$ or $(ii)$ or $(iii)$ follows from
Lemma~\ref{T_K_3_T'_K_3_first}. 
Therefore,  $S_T\cap V(D)\ne \emptyset$. Note that $1\le |S_T\cap S_{T'}|\le 4$, and we
distinguish four cases according to $|S_T\cap S_{T'}|$.

Suppose $|S_T\cap S_{T'}|=4$. Then $S_{T'}\cap V(B)=\emptyset$ and
$|S_T\cap V(C)|=|S_T\cap V(D)|=1$. 
Therefore, by the minimality of $|V(A)|$, 
$B\cap D\ne \emptyset$. Hence, $S_T-V(C)$ is a 5-cut in $G$ and
$V(T)\subseteq S_T-V(C)$. By the choice of $(T,S_T,A,B)$ that $|V(A)|$ is  minimum,  $|V(B\cap D)|\ge 5$. Now $(i)$ or $(ii)$
or $(iii)$ follows from Lemma~\ref{5cut_triangle}.

Consider $|S_T\cap S_{T'}|=3$. Suppose for the moment $S_{T'}\cap
V(B)=\emptyset$.  Then $|S_T\cap V(C)|= 2$  as $|(S_T \cup
S_{T'})-V(A \cup D)| = 5$.  So $B\cap D=\emptyset$ as otherwise $S_T-V(C)$ would be 
a 4-cut in $G$. However, this implies $|V(D)|<|V(A)|$, contradicting
the choice of $(T,S_T,A,B)$ that $|V(A)|$ is minimum.  So $S_{T'}\cap
V(B)\ne \emptyset$. 
Therefore, since  $|S_{T'}|=6$, we have $|S_{T'}\cap V(B)|=1$ and  $S_{T'}\cap V(A)=\{a,b\}$. Since  $|(S_T \cup
S_{T'})-V(A \cup D)| = 5$, 
$|S_T\cap V(C)|=1$. This is a contradiction, as we have $|S_T\cap
V(C)|\ne 1$ or $|S_{T'}\cap V(A)|\ne 2$.

Now let $|S_T\cap S_{T'}|=2$. First, assume $|S_T\cap V(C)|=1$. Then $|S_{T'}\cap
V(B)|=2$ (as $|(S_T\cup S_{T'})-V(A\cup D)|=5$) and, hence, $|S_{T'}\cap
V(A)|=2$ (as $|S_{T'}|=6$), a contradiction. 
So we may assume that $|S_T\cap V(C)|\ge 2$, which implies $|S_{T'}\cap V(B)|\le
1$ as  $|(S_T\cup S_{T'})-V(A\cup D)|=5$. Hence, since
$|S_T|=|S_{T'}|=6$, $|S_{T'}\cap V(A)|\ge 3$ and $|S_T\cap V(D)|\le 2$. 
Therefore, by the minimality of $|V(A)|$, $B\cap D\ne \emptyset$. Thus $(S_T\cap S_{T'})-V(A\cup C)$ is a 
5-cut in $G$ and contains $V(T)$. So $|V(B\cap D)|\ge 5$ by the minimality of $|V(A)$. 
Now $(i)$ or $(ii)$ or $(iii)$  follows from Lemma~\ref{5cut_triangle}. 

Finally, assume $|S_T\cap S_{T'}|=1$. If $|S_{T'}\cap V(B)|= 2$ then $|S_T\cap V(C)|=2$ (as $|(S_T\cup
S_{T'})-V(A\cup D)|=5$);  so $(iv)$ holds.  
If $|S_{T'}\cap V(B)|=3$ then
$|S_T\cap V(C)|=1$ (since $|(S_T\cup
S_{T'})-V(A\cup D)|=5$) and $S_{T'}\cap V(A)=\{a,b\}$ (as $|S_{T'}|=6$), a contradiction. 
Hence, we may assume $|S_{T'}\cap V(B)| \leq 1$.
Then $|S_{T}\cap V(C)| \geq 3$ (since $|(S_T\cup
S_{T'})-V(A\cup D)|=5$), $|S_{T'}\cap V(A)|\ge 4$, and $|(S_T\cup
S_{T'})-V(A\cup C)|\le 4$. Hence, since $G$ is 5-connected, $B\cap
D=\emptyset$; so 
$|V(D)|<|V(A)|$. However, this shows that $(T',S_{T'},D,C)$
contradicts the choice of $(T,S_T,A,B)$. \qed

\medskip

Next, we take care of the case when $(iv)$ of Lemma~\ref{reduction} holds.

\begin{lem}
\label{T_K_3_T'_K_3_third}
Let $G$ be \textsl{a} 5-connected nonplanar graph and $x \in V(G)$, and suppose for any $H\subseteq G$ with $x\in V(H)$ and $H\cong K_2$ or
$H\cong K_3$, $G/H$ is not 5-connected. 
Let $(T,S_T,A,B) \in {\cal Q}_x$ with $|V(A)|$ minimum and $(T',S_{T'},C,D) \in {\cal Q}_x$ with $T' \cap A \neq \emptyset$.
Suppose $T\cap C=\emptyset$, $S_T\cap S_{T'}=\{x\}$ and $|S_T\cap V(C)|=|S_{T'}\cap V(B)|=2$. 
Then one of the following holds: 
\begin{itemize}
\item[(i)] $G$ contains \textsl{a} $TK_5$ in which $x$ is not \textsl{a} branch vertex.
\item[(ii)] $G$ contains  $K_4^-$.
\item[(iii)] There exist $x_1,x_2,x_3\in N(x)$ such that, for any $y_1,y_2\in N(x)-\{x_1,x_2,x_3\}$,
$G':=G-\{xv:v\notin \{x_1,x_2,x_3,y_1,y_2\}\}$ contains  $TK_5$. 
\end{itemize}
\end{lem}

\pf  We may assume $T\cong K_3$ (by Lemma~\ref{T_K_2}) and $T'\cong
K_3$ (by Lemma~\ref{T_K_3_T'_K_2}). By Lemma~\ref{cut_4-}, we may assume $|V(A)|\ge 5$. We may further assume that
$|S_T|=|S_{T'}|=6$; for, otherwise, the assertion follows
from Lemma~\ref{5cut_triangle}.

Let $V(T)=\{x,x_1,x_2\}$, $V(T')=\{x,a,b\}$,  $S_T \cap V(C) =
\{p_1,p_2\}$, $S_{T'}\cap V(B) = \{c_1,c_2\}$,  $S_{T'} \cap V(A)= \{a,b,q\}$, and $S_T\cap V(D) = \{x_1,x_2,w\}$. 
Since $T\cap C=\emptyset$, we may assume by Lemma~\ref{T_K_3_T'_K_3_first} that 
$A\cap C=\emptyset$. Then $B \cap C \neq \emptyset$ by the minimality of $|V(A)|$.  

We may assume $N(p_1)\cap V(A)=\{a,q\}$ and $N(p_2)\cap
V(A)=\{b,q\}$.  To see this, for $i\in [2]$,
let $S_i:=(S_T-\{p_i\})\cup (N(p_i)\cap \{a,b,q\})$ which is a cut in
$G$ and containing $V(T)$. If $N(p_i)\cap
\{a,b,q\}=\emptyset$ then 
$|S_i|=5$ and the assertion of this lemma
follows from Lemma~\ref{5cut_triangle}.  If  $|N(p_i)\cap
\{a,b,q\}|=1$ then $(T,S_i, A-(N(p_i)\cap \{a,b,q\}), S_i,
G[B+p_i])\in {\cal Q}_x$, contradicting the choice of $(T,S_T,A,B)$ that
$|V(A)|$ is minimum. Hence, we may assume that $|N(p_i)\cap \{a,b,q\}|\ge 2$ for $i\in [2]$.
We may assume
$\{a,b\}\not\subseteq N(p_i)$ for $i\in [2]$; as otherwise, $G[T'+p_i]$ contains $K_4^-$ and $(ii)$ holds. 
Moreover, $N(p_1)\cap \{a,b,q\}\ne N(p_2)\cap \{a,b,q\}$, as otherwise,  
$S:=(S_T-\{p_1,p_2\})\cup (N(p_1)\cap \{a,b,q\})$ is a cut in 
$G$ containing $V(T)$; so $(T,S,A-(N(p_1)\cap \{a,b,q\}),
G[B+\{p_1,p_2\}])\in {\cal Q}_x$, contradicting the choice of
$(T,S_T,A,B)$ with $|V(A)|$ minimum.
Hence, we may assume $N(p_1)\cap V(A)=\{a,q\}$ and $N(p_2)\cap
V(A)=\{b,q\}$.

Note that $N(x_i)\cap V(B)\ne \emptyset$ for $i\in [2]$; for,
otherwise, $S:=V(T')\cup \{q,x_{3-i},w\}$ is a cut in $G$, and
$(T',S,G[(A\cap D)+x_i],G[B+\{p_1,p_2\}])\in {\cal Q}_x$, contradicting the
choice of $(T,S_T,A,B)$ that $|V(A)|$ is minimum. 
Moreover, we may assume $N(w)\cap V(B)\ne\emptyset$; as otherwise,
$S_T-\{w\}$ is a 5-cut in $G$ and $V(T)\subseteq S_T-\{w\}$, and the
assertion of this lemma follows from  Lemma~\ref{5cut_triangle}. 

We wish to  prove $(iii)$ with $x_3=b$. Let  $y_1,y_2\in N(x)-\{x_1,x_2,x_3\}$ be distinct.  Choose $v \in  \{y_1,y_2\} - \{a\}$. 
We may assume $v \not\in \{p_1,p_2\} $, as  otherwise $G[T'+v]$ contains $K_4^-$ and $(ii)$ holds.
 By Lemma \ref{2independentPaths-b-x1x2}, we may choose $t \in
 N(b)\cap V(A-a)$ such that 
$G[(A -a) +\{b,q,x_1,x_2,w\} ]$ has independent paths $P_1,P_2,P_3,P_4,P_5$ from $t$ to $b,x_1,x_2,w,q$ respectively. 
We distinguish four cases according to the location of $v$.

\medskip

{\it Case} 1. $v \in V(B)$. 

Let $W$ be the component of $B$ containing  $v$. First, suppose 
$N(x_i)\cap W \neq \emptyset $ for $i\in [2]$. Then there exists $v^*\in V(W)$ such that 
$G[W+\{x_1,x_2\}]$ has three independent paths from $v^*$ to $v,x_1,x_2$, respectively. Hence by Lemma~\ref{Perfect}, $G[W +(S_T-\{x\})]$ 
has  independent paths $Q_1,Q_2,Q_3,Q_4$ from $v^*$ to $v,x_1,x_2, u$,
respectively, and internally disjoint from $S_T$, where $u\in
S_T-\{x,x_1,x_2\}$. 
If $u=w$ then  $T \cup (P_1\cup bx) \cup P_2 \cup P_3 \cup 
(Q_1\cup vx) \cup Q_2 \cup Q_3 \cup (Q_4\cup  P_4)$
is a $TK_5$ in $G'$ with branch vertices $t,v^*,x,x_1,x_2$. If 
$u=p_i$ for some $i\in [2]$ then $T  \cup (P_1\cup bx) \cup P_2 \cup P_3 \cup (Q_1\cup vx)\cup Q_2 \cup Q_3\cup (Q_4\cup p_iq\cup  P_5)$
is a $TK_5$ in $G'$ with branch vertices $t,v^*,x,x_1,x_2$.

Thus, we may assume that $N(x_1)\cap W = \emptyset$.  
Since $G$ is 5-connected, $N(x_2)\cap W \ne \emptyset$. So $G[W +
(S_T-\{x_1\})]$ has independent paths $Q_1,Q_2,Q_3,Q_4,Q_5$ from $v$
to $x,x_2,w,p_1,p_2$, respectively. Clearly, we may assume that $Q_1=vx$.  Since $N(x_1)\cap V(B)\ne \emptyset$, let $W'$ be a component of $B$ with $N(x_1)\cap V(W')\ne \emptyset$. 
Since $G$ is 5-connected, there exists $i\in [2]$ such that
$N(p_i)\cap V(W')\ne \emptyset$.
Hence, $G[W'+\{x_1,p_i\}]$ has a path $R$ from $x_1$ to $p_i$, and, by symmetry, assume 
$R$ is from $x_1$ to $p_1$. 
Now $T \cup (P_1\cup bx) \cup P_2 \cup P_3 \cup Q_1 \cup Q_2  \cup
(Q_3\cup P_4) \cup (Q_4\cup R)$ is a $TK_5$ in $G'$ with branch vertices $t,v,x,x_1,x_2$.

\medskip

{\it Case} 2. $v \in V(A \cap D) $. 

First, we show that $G[(A \cap D) +\{q,w,x,x_1,x_2\}]$ has 
independent paths $P_1',P_2',P_3',P_4',$ $P_5'$ from $v$ to
$q,x,x_1,x_2,w$, respectively (and we may assume
that $P_2'=vx$).
This is clear if $G[(A \cap D) + \{q,w,x_1,x_2\}]$ has independent paths from $v$ to $q,x_1,x_2,w$, respectively. 
So we may assume that  $G[(A \cap D) + \{q,w,x_1,x_2\}]$ has a separation $(G_1,G_2)$ such that $|V(G_1 \cap G_2)| \leq 3$, 
$v \in V(G_1-G_2)$ and $\{q,w,x_1,x_2\} \subseteq V(G_2)$. 
Then $S:=V(T')\cup V(G_1\cap G_2)$ is a cut in $G$, and
$(T',S,G_1-G_2,G-S-G_1)\in {\cal Q}_x$, contradicting the choice of
$(T,S_T,A,B)$ that $|V(A)|$ is minimum.

Suppose $B$ has a component $W$ such that  $N(x_i)\cap W \neq
\emptyset$ for $i\in [2]$. 
Then there exists $z \in V(W)$ such that $G[W + \{x_1,x_2\}]$ has
independent paths from $z$ to $x_1, x_2$, respectively. 
Hence  by Lemma~\ref{Perfect},  $G[W + (S_T-\{x\})]$ has four
independent paths $Q_1,Q_2,Q_3,Q_4$ from $z$ to $x_1,x_2,u_1,u_2$, 
respectively, and internally disjoint from $S_T$, where $u_1,u_2\in
\{w,p_1,p_2\}$ are distinct. 
If $\{u_1,u_2\} = \{w,p_1\}$ then  we may assume $u_1=w$ and
$u_2=p_1$; now $T\cup P_2'\cup P_3'\cup P_4'\cup Q_1\cup Q_2\cup (Q_3\cup P_5')
\cup (Q_4\cup p_1abx)$ is a $TK_5$ in 
$G'$ with branch vertices $v, x,x_1,x_2,z$.
If $\{u_1,u_2\} = \{w,p_2\}$ then  we may assume $u_1=w$ and
$u_2=p_2$; now $T\cup P_2'\cup P_3'\cup P_4'\cup Q_1\cup Q_2\cup (Q_3\cup P_5')\cup
(Q_4\cup p_2bx)$ is a $TK_5$ in 
$G'$ with branch vertices $v, x,x_1,x_2,z$. So assume $\{u_1,u_2\} =
\{p_1,p_2\}$. We may further assume $u_i=p_i$ for $i\in [2]$. Then 
$T\cup P_2'\cup P_3'\cup P_4'\cup Q_1\cup Q_2 \cup (Q_3\cup p_1q\cup
P_1')\cup (Q_4\cup p_2bx)$ is a $TK_5$ in $G'$ with branch vertices $v,x,x_1,x_2,z$.

Hence, we may assume that no component of $B$ contains neighbors of both $x_1$ and $x_2$. Since $G$ is 5-connected, we may assume  by symmetry that $Z$ is a 
component of $B$ such that $N(x_1)\cap V(Z)=\emptyset$ and $N(x_2)\cap V(Z)\ne \emptyset$. Again, since $G$ is 5-connected, $G[Z+(S_T-\{x_1\})]$ 
has five independent paths $Q_1,Q_2,Q_3,Q_4,Q_5$ from some $z\in V(Z)$ to
$x_2,w,p_1,p_2,x$, respectively. Since $N(x_1)\cap V(B)\ne \emptyset$,
let $Z'$ be a component of $B$ with $N(x_1)\cap Z'\ne\emptyset$. Then $N(x_2)\cap V(Z')=\emptyset$. So $G[Z'+\{x_1,p_1\}]$ contains a path $R$ from $x_1$ to $p_1$. 
Now $T\cup P_2'\cup P_3'\cup P_4'\cup (Q_4\cup p_2bx)\cup Q_1 \cup
(Q_3\cup R)\cup (Q_2\cup P_5')$ is a $TK_5$ in $G'$ with branch vertices $v, x,x_1,x_2,z$.

\medskip

{\it Case} 3. $v = q $.

Suppose $B$ has a component $Z$ such that $\{w,x_1,x_2\}\subseteq N(Z)$.
Then there exists $z \in V(Z)$ such that $G[Z+ \{w,x_1,x_2\} ]$ has independent paths from $z$ to $w,x_1, x_2$, respectively. 
By Lemma~\ref{Perfect}, $G[Z+ (S_T-\{x\})]$ has independent paths
$Q_1,Q_2,Q_3,Q_4$ from $z$ to $x_1, x_2,w,u$, respectively,  and 
internally disjoint from $S_T$, where $u\in \{p_1,p_2\}$.
Let $S = Q_4\cup p_1abx$ if $u = p_1$ and $S = Q_4\cup p_2bx$ if $u = p_2$. 
Then $T\cup Q_1 \cup Q_2 \cup S \cup (P_4\cup Q_3) \cup P_2 \cup P_3
\cup (P_5\cup qx) $ is a $TK_5$ in $G'$ with branch vertices $t, x,x_1,x_2,z$.

So we may assume that no component of $B$ is adjacent to all of
$x_1,x_2$ and $w$. Since $N(w)\cap V(B)\ne \emptyset$, there exists a
component $Z$ of $B$ such that 
$N(w)\cap V(Z)\ne\emptyset$. 
Since $G$ is 5-connected, we may assume by symmetry that 
$N(x_2)\cap V(Z)\ne\emptyset$. Then $N(x_1)\cap V(Z)=\emptyset$.
 Since $G$ is $5$-connected, $G[Z+(S_T-\{x_1\})]$ has independent
 paths $Q_1,Q_2,Q_3,Q_4,Q_5$ from some $z\in V(Z)$ to $x_2,w,p_1,p_2,x$,
 respectively. 
Since $N(x_1)\cap V(B)\ne \emptyset$, there exists some component $Z'$
of $B$ with $N(x_1)\cap V(Z')\ne \emptyset$. Hence, $N(x_2)\cap
V(Z')=\emptyset$ or $N(w)\cap V(Z')=\emptyset$; so 
$G[Z'+\{x_1,p_1\}]$ contains a path $R$ from $x_1$ to $p_1$. Now 
$T \cup Q_1\cup (Q_3\cup R) \cup (Q_4\cup p_2bx) \cup (P_4\cup Q_2) \cup P_2 \cup P_3 \cup (P_5\cup qx) $ is  a $TK_5$ in $G'$ with branch vertices $t,x,x_1,x_2,z$.

\medskip

{\it Case} 4. $v = w $.

Suppose $B$ has a component $Z$ such that  $\{w,x_1,x_2\}\subseteq N(Z)$. 
Then there exists $z \in V(Z)$ such that $G[Z+\{w,x_1,x_2\} ]$ has three independent paths from $z$ to $w,x_1, x_2$, respectively. 
Hence, by Lemma~\ref{Perfect}, $G[Z+ (S_T-\{x\})]$ has independent
paths $Q_1,Q_2,Q_3,Q_4$ from $z$ to $x_1, x_2,w,u$, respectively, and
internally disjoint from $S_T$,  where $u=p_i$ for some $i\in [2]$.
Then $T \cup Q_1\cup Q_2 \cup (Q_3\cup wx)  \cup (P_1\cup bx)\cup P_2 \cup P_3\cup (P_5\cup qp_i\cup Q_4) $ is a $TK_5$ in $G'$ with branch vertices $t,x,x_1,x_2,z$.

Hence, we may assume that no component of $B$ is adjacent to all of $w,x_1,x_2$. 
 Since $N(w)\cap V(B)\ne \emptyset$, $B$ has a component $Z$ such that 
$N(w)\cap V(Z)\ne\emptyset$. Since $G$ is 5-connected, we may assume by symmetry that 
$N(x_2)\cap V(Z)\ne\emptyset$. Then $N(x_1)\cap V(Z)=\emptyset$.
 Since $G$ is $5$-connected, $G[Z+(S_T-\{x_1\})]$ has five independent paths $Q_1,Q_2,Q_3,Q_4,Q_5$ from $z$ to $x_2,w,p_1,p_2,x$, respectively. 
Since $N(x_1)\cap V(B)\ne \emptyset$, $B$ has  a component  $Z'$ such
that  $N(x_1)\cap V(Z')\ne \emptyset$. Then $N(x_2)\cap
V(Z')=\emptyset$ or $N(w)\cap V(Z')=\emptyset$; so $G[Z'+\{x_1,p_1\}]$ contains a path $R$ from $x_1$ to $p_1$. Now 
$T \cup Q_1\cup (Q_2\cup wx) \cup (Q_3\cup R) \cup (P_1\cup bx) \cup
P_2 \cup P_3 \cup (P_5\cup qp_2\cup Q_4) $ is  a $TK_5$ in $G'$ with branch vertices $t,x,x_1,x_2,z$. \qed

\medskip
 We end this section with the following lemma which deals with another
 special  case when $(T,S_T,A,B) \in {\cal Q}_x$ with $|V(A)|$ minimum,
 $(T',S_{T'},C,D) \in {\cal Q}_x$ with $T' \cap A\ne\emptyset$, and
 $A \cap C = \emptyset$.

\begin{lem}
\label{T_K_3_T'_K_3_second}
Let $G$ be \textsl{a} 5-connected nonplanar graph and $x \in V(G)$ such that
for any $H\subseteq G$ with $x\in V(H)$ and $H\cong K_2$ or $H\cong
K_3$, $G/H$ is not 5-connected. 
Let $(T,S_T,A,B) \in {\cal Q}_x$ with $|V(A)|$ minimum, and $(T',S_{T'},C,D) \in {\cal Q}_x$ with $T' \cap A \neq \emptyset$.
Suppose $A \cap C = \emptyset$, $|S_T| = 6$, $|S_{T'}| = 6$, $V(T') \cap
S_{T} = \{x,b\}$, $V(T'\cap A)=S_{T'}\cap V(A) = \{a\}$ and $V(C) \cap S_{T} =
\emptyset$. Then, one of the following holds: 
\begin{itemize}
\item[$(i)$] $G$ contains \textsl{a} $TK_5$ in which $x$ is not \textsl{a} branch vertex.
\item[$(ii)$] $G$ contains $K_4^-$.
\item[$(iii)$] There exist distinct $x_1,x_2\in N(x)$ such that for any
  distinct $y_1,y_2\in N(x)-\{b,x_1,x_2\}$, 
$G':=G-\{xv:v\notin \{x_1,x_2,b,y_1,y_2\}\}$ contains  $TK_5$.  
\end{itemize}
\end{lem}

\pf  By assumption, $V(T')=\{a,b,x\}$ with $a\in V(A)$ and $b,x\in
S_T\cap S_{T'}$. Let $V(T) =\{x,x_1,x_2\}$ and  $S_T = \{b, x,x_1,x_2,
x_3,x_4\}$. We wish to prove $(iii)$ with $x_3=b$; so let $y_1,y_2\in
N(x)-\{b,x_1,x_2\}$ be distinct.  Let $v \in  \{y_1,y_2\} -\{a\}$.

Note that $B \cap C \neq \emptyset $ as $S_{T'}$ is a cut. 
So $|(S_T \cup S_{T'}) - V(A \cup D) | \geq 5 $. 
Moreover, we may assume $A \cap D \neq \emptyset $ by Lemma \ref{cut_4-}.
So $|(S_T \cup S_{T'}) - V(B \cup C) | \geq 7 $ by the minimality of $|V(A)|$. 
Since $|S_T| = |S_{T'}| = 6$, 
$$|(S_T \cup S_{T'} )- V(A \cup D) | = 5 \mbox{ and } |(S_T \cup S_{T'}) - V(B \cup C)| = 7.$$

We may assume that $N(x_i) \cap V(B) \neq \emptyset$ for $i  \in
[2]$. For, suppose this is not true and by symmetry assume $N(x_1)
\cap V(B) = \emptyset $. Let $S = (S_T - \{x_1\}) \cup \{a\} $,
$C'=B$, and $D'=G[(A-a)+x_1]$. Then $(T',S,C',D') \in {\cal Q}_x$.
 We now apply Lemma~\ref{lemma_butterfly} to $(T,S_T,A,B)$ and $(T',S,C',D')$. Note
that $|S\cap S_T| = 5$, $V(A \cap C') = S_T \cap V(C') = S \cap V(B) =V(B \cap D') =
\emptyset$, and $|S\cap V(A)| = |S_T \cap V(D')|=|V(T\cap D')| = 1$.  
To verify the other condition in Lemma~\ref{lemma_butterfly}, let
$(H,S_H,C_H,D_H)\in {\cal Q}_x$ such that $H\cong K_2$ or
$H\cong K_3$. Then we may assume that $H\cong K_3$  when $H\cap A\ne
\emptyset$ (by Lemma~\ref{T_K_3_T'_K_2}) and that $|V(H\cap A)|\le 1$ (by
Lemmas~\ref{reduction} and \ref{T_K_3_T'_K_3_third}).
Therefore, the assertion of this lemma follows from Lemma~\ref{lemma_butterfly}. 
Hence, we may assume $N(x_i) \cap B \neq \emptyset$ for $i  \in
[2]$.

We may assume that for any component $W$ of $B$, $N(b)\cap W \neq \emptyset$; for,
otherwise, $S_T - \{b\}$ is a $5$-cut  in $G$, and the assertion of
this lemmas follows from  
Lemma~\ref{5cut_triangle}. 
We consider three cases according to the location of $v$.

\medskip

{\it Case} 1. $v \in V(B)$. 

Let $B_v$ be the component of $B$ containing  $v$.  First, suppose
$N(x_i) \cap V(B_v) \neq \emptyset $ for $i\in [2]$. 
Then $G[B_v+\{x_1,x_2\}]$ has independent paths from some  $v^*\in
V(B_v)$ to $v,x_1,x_2$, respectively.
Thus, by Lemma~\ref{Perfect},  $G[B_v +S_T-x]$ has independent paths
$P_1,P_2,P_3,P_4$ from $v^*$ to $v,x_1,x_2,u$, respectively, and
internally disjoint from $S_T$,  where $u\in \{b,x_3,x_4\}$. 
Suppose $u=b$. By Lemma~\ref{2independentPaths-b-x1x2}, we may assume
that  $G[A+\{b,x_1,x_2\}]$ contains independent paths $R_1,R_2$ from $b$ to
$x_1,x_2$, respectively. Then $T\cup R_1\cup R_2\cup bx\cup (P_1\cup
vx)\cup P_2\cup P_3\cup P_4$ is a $TK_5$ in $G'$ with branch vertices
$b,v^*,x,x_1,x_2$. So we may assume by symmetry that $u=x_3$. 
By Lemma~\ref{2independentPaths-b-x1x2} again,  we may choose $t \in
N(b)\cap  V(A-a)$ and let  $Q_1,Q_2,Q_3,Q_4,Q_5$ be independent paths in $G[(A-a) + \{b,x_1,x_2,x_3,x_4\}]$ from 
$t$ to $b,x_1,x_2,x_3,x_4$, respectively. 
Then, $T \cup (Q_1\cup bx)\cup  Q_2 \cup Q_3 \cup
(P_1\cup vx)\cup P_2 \cup P_3 \cup (P_4\cup Q_4)$ is a $TK_5$ in $G'$ with branch vertices $t, v^*,x,x_1,x_2$.

Therefore, we may assume by symmetry that  $N(x_1) \cap V(B_v) = \emptyset $. 
Since $G$ is $5$-connected, $G[B_v +S_T -x_1]$ has  independent paths
$P_1,P_2,P_3,P_4,P_5$ from $v$ to $x,b,x_2,x_3,x_4$, respectively, and 
we may assume that $P_1=vx$. Since $N(x_1) \cap V(B) \neq \emptyset$, $B$ has a component  
$B_{x_1}$ such that $N(x_1) \cap V(B_{x_1}) \neq \emptyset$. 
Again, since $G$ is 5-connected, $N(x_j)\cap V(B_{x_1})\ne \emptyset$
for some $j\in \{3,4\}$, and we may assume $j=3$. Then $G[B_{x_1}+\{x_1,x_3\}]$ contains a path $Q$ from $x_1$ to $x_3$.
Let $t \in N(b)\cap V(A -a)$. By Lemma~\ref{2independentPaths-b-x1x2},
we may assume that  $G[(A -a) +\{b,x_1,x_2,x_3,x_4\}]$ has independent paths 
 $Q_1,Q_2,Q_3,Q_4,Q_5$  
from $t$ to $b,x_1,x_2,x_3,x_4$, respectively. 
Then $T \cup  (Q_1\cup bx)\cup Q_2 \cup Q_3 \cup (P_5\cup Q_5) \cup
(P_4\cup Q)\cup P_1 \cup P_3 $ is a $TK_5$ in $G'$ with branch vertices
$t, v,x,x_1,x_2$.

\medskip

{\it Case 2}. $v \in V(A \cap D)$. 

We claim that $G[(A -a) +\{x,x_1,x_2,x_3,x_4\}]$ has  independent
paths $P_1,P_2,P_3,P_4,P_5$ from $v$ to $x,x_1,x_2,x_3,x_4$,
respectively (and we may assume $P_1 = vx$). This is clear if
$G[(A -a) + \{x_1,x_2,x_3,x_4\}]$ has  independent paths from $v$ to $x_1,x_2,x_3,x_4$, respectively;
so we may assume such paths do not exist.
Then there exists a separation $(G_1,G_2)$ in $G[(A-a) +
\{x_1,x_2,x_3,x_4\}]$ such that $|V(G_1 \cap G_2)| \leq 3$, $v \in
V(G_1-G_2)$, and $\{x_1,x_2,x_3,x_4\} \subseteq V(G_2)$. 
Let $S:=V(G_1\cap G_2)\cup V(T')$, which is a cut  in $G$ of size at
most 6. Since $G$ is $5$-connected, $|V(G_1 \cap G_2)| \ge 2$. 
Then, $(T', S, G_1-G_2, (G-S)-G_1)\in {\cal Q}_x$ and $1\le |V(G_1-G_2)| \leq |V(A - a)|<|V(A)|$, contradicting the choice of $(T,S_T,A,B)$ that $|V(A)|$ is minimum.

Suppose that $B$ has a component $W$ such that $N(x_i)\cap V(W) \neq
\emptyset $ for $i\in [2]$. Then there exists $w \in V(W)$ such that $G[W+b]$ has independent paths from $w$ to $x_1,x_2,b$, respectively. 
By Lemma~\ref{Perfect}, $G[B + S_T]$ has independent paths
$Q_1,Q_2,Q_3,Q_4,$ $Q_5$ from $w$ to $x_1,x_2,b,u_1,u_2$,
respectively, and internally disjoint from $S_T$,  where 
$u_1,u_2\in \{x,x_3,x_4\}$ are distinct. By symmetry, we may assume $u_1 = x_3$. 
Then $T\cup P_1\cup P_2\cup P_3\cup Q_1\cup Q_2\cup (Q_3\cup bx)\cup (Q_4\cup P_4)$ is 
a $TK_5$ in $G'$ with branch vertices $v,w,x,x_1,x_2$. 

Hence, we may assume that no component of $B$ is adjacent to
both $x_1$ and $x_2$. Let $W$ be a component of $B$ such that $N(x_2) \cap V(W) \ne \emptyset$. 
Then $N(x_1) \cap V(W) = \emptyset$. 
Since $G$ is 5-connected, $G[W+ S_T - x_1 ]$ has  independent paths
$Q_1,Q_2,Q_3,Q_4,Q_5$ from some $w\in V(W)$ to $b,x_2,x_3,x_4,x$, respectively. 
Since  $N(x_1) \cap V(B) \neq \emptyset$,
$B$ has a component $B_x$ such that $N(x_1) \cap V(B_x) \neq \emptyset$. Then $N(x_2)\cap V(B_x)= \emptyset$. 
Again, since $G$ is 5-connected, $G[B_x+\{x_1,x_3\}]$ contains a path $R$ from $x_1$ to $x_3$. 
Now $T\cup P_1 \cup P_2 \cup P_3\cup (Q_1\cup bx) \cup Q_2 \cup
(Q_3\cup R) \cup (Q_4\cup P_5) $ is a $TK_5$ in $G'$ with branch vertices $v,w,x,x_1,x_2$.

\medskip

{\it Case} 3. $v \in S_T$. 

We may assume that $v = x_3$. 
By Lemma \ref{2independentPaths-b-x1x2}, we may assume $t \in N(b)\cap
V(A -a)$ and $G[(A -a) + \{b,x_1,x_2,x_3,x_4\}]$ has independent paths $P_1,P_2,P_3,P_4,P_5$
from $t$ to $b,x_1,x_2,$ $x_3,x_4$,  respectively, with $P_1=tb$. Also
by Lemma~\ref{2independentPaths-b-x1x2}, we may assume that  $G[A +
\{b,x_1,x_2\}]$ has independent paths $Q_1,Q_2$ from $b$ to $x_1,x_2$,  respectively.

Suppose $B$ has a component $W$ such that $\{x_1,x_2\}\subseteq N(W)$.  
Then there exists $w \in V(W)$ such that $G[W + \{b,x_1,x_2\}]$ has
independent paths from $w$ to $b,x_1,x_2$, respectively. 
So by Lemma~\ref{Perfect}, $G[B + S_T]$ has independent paths
$R_1,R_2,R_3,R_4,R_5$ from $w$ to $x_1,x_2,b,u_1,u_2$, respectively,
and internally disjoint from $S_T$, where $u_1,u_2 \in \{x,x_3,x_4\}$ are distinct. 
Assume by symmetry that $u_1\in \{x_3,x_4\} $.
If $u_1 = x_3$, then $T\cup bx\cup Q_1 \cup Q_2 \cup R_1 \cup R_2 \cup R_3 \cup (R_4\cup x_3x) $ is a $TK_5$ in $G'$ with branch vertices $b, w,x,x_1,x_2$. 
If $u_1 = x_4$, then $T\cup (P_4\cup x_3x) \cup P_2 \cup P_3 \cup R_1 
\cup R_2 \cup (R_3\cup bx) \cup (R_4\cup P_5)$ is  a $TK_5$ in $G'$ with branch vertices $t,w,x,x_1,x_2$. 

Thus, we may assume that no component of $B$ is adjacent to both $x_1$ and $x_2$. 
Since $G$ is $5$-connected, we may assume by symmetry that $W$ is a component of $B$ such that $N(x_2) \cap V(W) \neq \emptyset$ and $N(x_1) \cap V(W) =\emptyset$. 
Let $w \in V(W)$. Since $G$ is 5-connected, $G[W +S_T-x_1]$  has independent paths $R_1,R_2,R_3,R_4,R_5$ from $w$ to $x,x_2,x_3,x_4,b$, respectively.
Since $N(x_1) \cap B \neq \emptyset$, $B$ has a component $B_x$ such
that  $N(x_1)\cap V(B_x)\ne \emptyset$.  Then $N(x_2) \cap V(B_x) = \emptyset$. 
Since $G$ is 5-connected, $G[B_x+\{x_1,x_4\}]$ contains a path $R$ from $x_1$ to $x_4$. 
Now $T\cup bx \cup Q_1 \cup Q_2 \cup R_2 \cup (R_3\cup x_3x) \cup R_5 \cup (R_4\cup R)$ is a $TK_5$ in $G'$ with branch vertices $b,w,x,x_1,x_2$. 
\qed


\section{Proof of Theorem \ref{main}}

In this section, we complete the proof of Theorem~\ref{main}, using
the lemmas we have proved so far. 
Let $G$ be a 5-connected nonplanar graph. We proceed to find a $TK_5$
in $G$. 
By Lemma~\ref{K4-}, we may assume that 
\begin{itemize}
\item [(1)] $G$ contains no $K_4^-$. 
\end{itemize}

Let $M$ denote a maximal connected subgraph of $G$ such that 
$$\mbox{$H:=G/M$ is 5-connected and nonplanar, and contains no
$K_4^-$.}$$ 
Note that $|V(M)|=1$ (i.e., $H=G$) is possible. Let $x$ denote the vertex
of $H$ resulting from the contraction of $M$. Then, for any $T\subseteq
H$ with $x\in V(T)$ and  $T\cong K_2$ or $T\cong K_3$,  one of the
following holds: $$ H/T \mbox{  contains } K_4^-, \mbox{ or } H/T \mbox{  is
planar, or } H/T  \mbox{ is not 5-connected}.$$ 
For convenience, we will use $x_T$ to denote
the vertex of $H/T$ resulting from the contraction of $T$.
We may assume that 

\begin{itemize}
\item [(2)]  for any $T\subseteq
H$ with $x\in V(T)$ and  $T\cong K_2$ or $T\cong K_3$, if $F$ is a
$TK_5$ in $H/T$ then $x_T$ is a branch vertex of $F$. 
\end{itemize}
For, suppose that $F$ is a $TK_5$ in $H/T$ in which
$x_T$ is not a branch vertex. If $x_T\notin V(F)$ then $F$ is also
$TK_5$ in $G$. So assume $x_T\in V(T)$. Let $u,v\in V(F)$ such that
$x_Tu,x_Tv\in E(F)$. Since $M$ is connected, $G[M+\{u,v\}]$ has a path
$P$ from $u$ to $v$. Thus, $(F-x)\cup P$ is a $TK_5$ in $G$. So we may
assume (2).
\medskip

Suppose there exists  $T\subseteq V(H)$ with $x\in V(T)$ and $T\cong K_2$
or $T\cong K_3$, such that $H/T$ is 5-connected and  planar. Then by
Lemma~\ref{contraction}, $H-T$ contains $K_4^-$, contradicting  (1). So 
\begin{itemize}
\item [(3)] for any  $T\subseteq H$ with $x\in V(T)$ and $T\cong K_2$ or
$T\cong K_3$, if $H/T$ is 5-connected then $H/T$ is nonplanar.
\end{itemize}

\medskip
We now show that 

\begin{itemize}
\item [(4)] if $T\subseteq H$ with $x\in V(T)$ and $T\cong K_2$ or
$T\cong K_3$ and if $x_1,x_2,x_3\in N_{H/T}(x_T)$ such that $H/T-\{x_Tv:
  v\notin \{u_1,u_2,x_1,x_2,x_3\}\}$ contains  $TK_5$ for every choice
  of distinct $u_1,u_2\in N_{H/T}(x_T)-\{x_1,x_2,x_3\}$, then $G$ contains
  $TK_5$. 
\end{itemize}
 To prove (4), let $A=N_G(M\cup T)=N_{H/T}(x_T)$. Consider the
 subgraph $G[M\cup T+A]$. Since $M\cup T$ is connected, there is a vertex
$v\in V(M\cup T)$ such that
$G[M\cup T+\{x_1,x_2,x_3\}]$ has  independent paths from $v$ to
$x_1,x_2,x_3$, respectively. Since $G$ is 5-connected, $G[M\cup T+
A]$ has five independent paths from $v$ to $A$ with only $v$ in common
and internally disjoint from $A$. Hence, by Lemma~\ref{Perfect}, there
exist distinct  $u_1,u_2\in
A-\{x_1,x_2,x_3\}$ such that $G[M\cup T+A]$
has five independent paths $P_1,P_2,P_3,P_4,P_5$ from $v$ to $x_1,x_2,x_3,u_1,u_2$,
respectively, and internally disjoint from $A$. Now suppose $F$ is a $TK_5$ in $H/T - \{x_Tv: v \not\in
\{x_1,x_2,x_3,u_1,u_2\} \}$. Then $F-x_T$ and the four paths among 
$P_1,P_2,P_3,P_4,P_5$ corresponding to the four edges at $x_T$ in $F$ form a $TK_5$ in
$G$. Hence, we may assume  (4). 

\medskip

By (3), we have two cases:   for some $T\subseteq H$ with $x\in V(T)$ and $T\cong K_2$ or
$T\cong K_3$, $H/T$ is 5-connected and nonplanar  but contains $K_4^-$; or  for every $T\subseteq H$ with $x\in V(T)$ and $T\cong K_2$ or
$T\cong K_3$, $H/T$ is not 5-connected.

\medskip

{\it Case} 1.  There exists $T\subseteq H$ with $x\in V(T)$ and $T\cong K_2$ or
$T\cong K_3$ such that $H/T$ is 5-connected and nonplanar, and $H/T$ contains $K_4^-$.

Let  $K\subseteq H/T$ such that $K\cong K_4^-$, and let $V(K)=\{x_1,x_2,y_1,y_2\}$ with $y_1y_2\notin E(H)$.
By (1), $x_T\in V(K)$.

\medskip
{\it Subcase} 1.1. $x_T$ has degree 2 in $K$. 

Then we may assume that the
notation is chosen so that $x_T=y_2$. By Lemma~\ref{y_2}, one of the following holds:
\begin{itemize}
\item [$(i)$] $H/T$ contains a $TK_5$ in which $x_T$ is not a branch vertex.
\item[$(ii)$] $H/T-x_T$ contains $K_4^-$.
\item[$(iii)$] $H/T$ has a $5$-separation $(G_1,G_2)$ such that $V(G_1\cap G_2)=\{x_T, a_1,a_2,a_3,a_4\}$, and  $G_2$ is 
the graph obtained from the edge-disjoint union of the $8$-cycle $a_1b_1a_2b_2a_3b_3a_4b_4a_1$ 
and the $4$-cycle $b_1b_2b_3b_4b_1$ by adding $x_T$ and the edges $x_Tb_i$ for $i\in [4]$.
\item[$(iv)$] For $w_1,w_2,w_3\in N_{H/T}(x_T)-\{x_1,x_2\}$, $H/T-\{x_Tv: v\notin \{w_1,w_2,w_3,x_1,x_2\}\}$ contains  $TK_5$.
\end{itemize}

Note that $(i)$ does not occur because of  (2), and $(ii)$ does not
occur because of (1).

Now suppose $(iii)$ occurs.  First, assume  $|V(G_1)|\ge 7$. Then by Lemma~\ref{8cycle}, 
 for any $u_1, u_2 \in N(x_T) - \{b_1,b_2,b_3\}$, $H/T- \{x_Tv: v \not\in
 \{b_1,b_2,b_3,u_1,u_2\} \}$ contains  $TK_5$. Hence, by (4)  (with
 $x_i$ as $b_i$ for $i\in [3]$), $G$ contains $TK_5$. So we may assume
 that $|V(G_1)|=6$, and let $v\in V(G_1-G_2)$.  By (1), $a_ia_{i+1}\notin E(G)$ for $i\in [4]$,
 where $a_5=a_1$. Hence, since  $G$ is 5-connected, $a_1a_3,a_2a_4\in
 E(G)$. Now $(H-x_T)-\{a_1v,a_1b_4,a_4v,a_4b_4\}$ is a $TK_5$ with
 branch vertices $a_2,a_3,b_1,b_2,b_3$, contradicting (2).

Finally, suppose $(iv)$ holds. Then, by (4) (with $w_1,w_2,w_3$ as
$x_3,u_1,u_2$, respectively), we see that $G$ contains $TK_5$.

\medskip

{\it Subcase 1.2}. $x_T$ has degree 3 in $K$.

Then we may assume that the notation is chosen so that $x_T=x_1$. By
Lemma~\ref{x1a}, one of the following holds:
\begin{itemize}
\item [$(i)$] $H/T$ contains a $TK_5$ in which $x_T$ is not a branch vertex.
\item [$(ii)$] $H/T-x_T$ contains $K_4^-$, or $H/T$ contains a $K_4^-$ in which $x_T$ is of degree 2. 
\item [$(iii)$] $x_2,y_1,y_2$ may be chosen so that for any distinct $z_0, z_1\in N_{H/T}(x_T)-\{x_2,y_1,y_2\}$, 
$H/T-\{x_Tv:v\notin \{z_0, z_1, x_2,y_1,y_2\}\}$ contains $TK_5$. 
\end{itemize}

By (2), $(i)$ does not occur. If $(ii)$ holds then, by (1), $H/T$ contains $K_4^-$ in which
$x_T$ is of degree 2; and we are back in Subcase 1.1. If $(iii)$ holds
then $G$ contains $TK_5$ by (4). 

\medskip

{\it Case} 2. $H/T$ is not 5-connected for each $T\subseteq H$ with $x\in V(T)$ and $T\cong K_2$ or
$T\cong K_3$. 

Let $\mathcal{Q}_x$ denote the set of all quadruples $(T,S_T,A,B)$,  such that 
\begin{itemize}
\item $T\subseteq V(H)$, $x\in V(T)$, and $T\cong K_2$ or $T\cong K_3$, 
\item $S_T$ is a cut in $H$ with $V(T) \subseteq S_T$,  $A$ is a nonempty union
  of components of $H - S_T$, and $B = H-S_T-A\ne \emptyset$,
\item if $T\cong K_3$ then $5\le |S_T| \le 6$, and 

\item if $T\cong K_2$ then $|S_T| = 5$, $|V(A)|\ge 2$, and $|V(B)|\ge 2$.
 \end{itemize}
Among all the quadruples in $\mathcal{Q}_x$, we select $(T,S_T,A,B)$
such that $|V(A)|$ is minimum. 

Since $K_4^-\not\subseteq H$,  $T\cong K_3$ (by Lemma~\ref{T_K_2})
and there exists $a \in V(A)$ such that $ax
\in E(H)$ (by Lemma \ref{existence_of_a} and by (2) and (4)). 
By Lemma~\ref{existence_of_T'}, there exists $T'\subseteq
H$ such that $x\in V(T')$ and $T'\cong K_2$ or $T'\cong K_3$, and there
exists $(T',S_{T'},C,D)\in {\cal Q}_x$. Again since
$K_4^-\not\subseteq H$, $T'\cong K_3$ by Lemma~\ref{T_K_3_T'_K_2} and
by (2) and (4).

We may assume, without loss of generality, that $T\cap C=\emptyset$. 
Hence, by Lemma~\ref{T_K_3_T'_K_3_first} and by (2) and (4),  $A \cap C = \emptyset $ (since
$K_4^-\not\subseteq H$). 
We may assume $B\cap C\ne \emptyset$; for otherwise,
$|V(A)|\le |V(C)|=|V(C)\cap S_T|\le 3$ and,  by Lemma~\ref{cut_4-},  $H$
contains $K_4^-$, a contradiction. 
  
We may assume that $|V(T')\cap S_T|=2$ for any
choice of $(T',S_{T'},C,D)\in {\cal Q}_x$ with $T'\cap A\ne
\emptyset$; otherwise, by Lemmas~\ref{reduction} and
\ref{T_K_3_T'_K_3_third}, we derive a contradiction to (2), or (4), or the fact
$K_4^-\not\subseteq H$.  Hence, since
$K_4^-\not\subseteq H$,  we have $A\cap D\ne \emptyset$  by Lemma~\ref{cut_4-}. 

Note that $|S_T|=|S_{T'}|=6$; for otherwise, by
Lemma~\ref{5cut_triangle}, we derive a contradiction to (2), or (4), or the fact
$K_4^-\not\subseteq H$. 
We claim that $|(S_T \cup S_{T'}) - V(B \cup C)| = 7$ and $|(S_T \cup S_{T'}) -V(A \cup D)| = 5$.  First, note that  $|(S_T \cup S_{T'}) - V(B \cup C)| \ge 7$; 
otherwise, $(T' , (S_T \cup S_{T'}) - V(B \cup C), A \cap D, G[B \cup C])
\in \mathcal{Q}_x$ and $1\le |V(A \cap D)| < |V(A)|$,  
contradicting the choice of $(T,S_T,A,B)$ with $|V(A)|$ minimum. Since $H$ is 5-connected and $B\cap C\ne \emptyset$,  $|(S_T \cup S_{T'}) -V(A \cup D)| \geq 5$. So the 
claim follows from the fact that $|(S_T \cup S_{T'}) - V(B \cup C)|+|(S_T \cup S_{T'}) -V(A \cup D)=|S_T| + |S_{T'}| = 12$.

If $S_T\cap V(C)=\emptyset$ for some choice $(T',S_{T'},C,D)$ then
$|S_{T'}\cap V(A)|=1$ as $|S_{T'}|=6$ and $|(S_T \cup S_{T'}) -V(A \cup D)|=5$; so by
Lemma~\ref{T_K_3_T'_K_3_second}, we derive a contradiction to (2), or
(4), or the fact
$K_4^-\not\subseteq H$. 

Hence,  we may assume that $$S_T\cap V(C)\ne
\emptyset$$ for any choice 
of $(T',S_{T'},C,D)\in {\cal Q}_x$ with $T'\cap A\ne \emptyset$. Then $2\le |S_T\cap S_{T'}|\le 4$ as $|(S_T\cup
S_{T'})-V(A\cup D)|=5$.

Suppose  $|S_T\cap S_{T'}|=4$. Then $|S_{T'}\cap V(B)|=0$ and $|S_T\cap
V(C)|=1$, as $|(S_T\cup S_{T'})-V(A\cup D)|=5$. Since $|S_T|=|S_{T'}|=6$,
$|S_T\cap V(D)|=1$ and $|S_{T'}\cap V(A)|=2$. Hence, $B\cap D\ne
\emptyset$ (since $|V(D)|\ge V(A)|$). So  $S_T-V(C)$ is a 5-cut in $H$
and $V(T)\subseteq S_T-V(C)$. Note $|V(B\cap D)|\ge 2$; for otherwise,
since  $H$ is 5-connected, $H[T\cup (B\cap D)]$ contains $K_4^-$,
a contradiction. Hence, by Lemma~\ref{5cut_triangle}, we derive a
contradiction to (2), or (4), or the fact $K_4^-\not\subseteq H$. 

Now assume  $|S_T\cap S_{T'}|= 3$. Then, $|S_{T'}\cap V(B)|\le 1$ as
$|(S_T\cup S_{T'})-V(A\cup D)|=5$ and $|S_T\cap V(C)|>0$. Suppose  $|S_{T'}\cap V(B)|=0$. Then
$|S_{T'}\cap V(A)|=3$ as $|S_{T'}|=6$.  So $|S_T\cap V(D)|=1$ since
$|(S_T\cup S_{T'})-V(B\cup C)|=7$. Thus, since $H$ is 5-connected,
$B\cap D= \emptyset$.  However, this implies that $|V(D)|<|V(A)|$, a
contradiction. So  $|S_{T'}\cap V(B)|=1$. Then $|S_{T'}\cap V(A)|=2$
as $|S_{T'}|=6$, and $|S_T\cap V(C)|=1$ as $|(S_T\cup S_{T'})-V(A\cup D)|=5$.
Let $q\in S_{T'}\cap V(A-T')$, $S':=(S_{T'}-\{q\})\cup (S_T\cap V(C))$, 
$C':=B\cap C$, and $D'=G[D+q]$. Then $(T',S',C',D')\in {\cal Q}_x$ with $T'\cap A\ne \emptyset$ and 
$T\cap C'=\emptyset$, However,  $S_T\cap V(C')=\emptyset$, a contradiction. 

Finally, assume $|S_T\cap S_{T'}|=2$. Suppose $|S_T\cap V(C)|\ge
2$. Then  $|S_{T'}\cap V(B)|\le 1$ (as $|(S_T\cup S_{T'})-V(A\cup D)|=5$),
and  $|S_{T'}\cap V(A)|\ge 3$ (as $|S_{T'}|=6$). So  $B\cap D\ne
\emptyset$ as $|V(D)|\ge |V(A)|$. 
Hence, $(S_T\cup S_{T'})-V(A\cup C)$ is a 5-cut in $H$ and contains
$V(T)$. If $|V(B\cap D)|=1$ then, since $H$ is 5-connected,  $H[T\cup
(B\cap D)]$ contains $K_4^-$, a contradiction. So $|V(B\cap D)|\ge
2$. Then, by 
Lemma~\ref{5cut_triangle}, we derive a contradiction to (2), or (4), or the fact
$K_4^-\not\subseteq H$. Therefore,  we may assume $|S_T\cap V(C)|=1$.
Hence, $|S_T\cap V(D)|=3$ (as $|S_T|=6$), $|S_{T'}\cap V(B)|=2$ (as
$|(S_T\cup S_{T'})-V(A\cup D)|=5$), and $|S_{T'}\cap V(A)|=2$ (as
$|S_{T'}|=6$).  Let $q\in S_{T'}\cap V(A-T')$, $S':=(S_{T'}-\{q\})\cup (S_T\cap V(C))$, 
$C':=B\cap C$, and $D'=G[D+q]$. Then $(T',S',C',D')\in {\cal Q}_x$ with $T'\cap A\ne \emptyset$ and 
$T\cap C'=\emptyset$, However,  $S_T\cap V(C')=\emptyset$, a contradiction. 
\qed

\section{Concluding remarks}

We have shown that every 5-connected nonplanar graph contains
$TK_5$. Thus, if a graph contains no $TK_5$ then it is planar, or admits a cut of 
size at most 4. This is a step towards a more useful structural description of the class of graphs containing 
no $TK_5$. There is a nice result for graphs containing no $TK_{3,3}$ due to Wagner \cite{Wa37a}: 
Every such graph is planar, or is a $K_5$, or  admits a  cut of size at most 2.

Mader~\cite{Ma98}  conjectured that every simple graph with  minimum degree at least 5 and no $K_4^-$ contains $TK_5$, and he
also asked the following. 
\begin{ques}\label{maderK23}
Does every simple
graph on $n\ge 4$ vertices with more than $12(n-2)/5$ edges contain $K_4^-$,
$K_{2,3}$, or $TK_5$?
\end{ques}
In a recent paper \cite{KMY15},  it is shown that an affirmative answer to 
Question~\ref{maderK23} implies the Kelmans-Seymour conjecture. As an independent approach to resolve the Kelmans-Seymour conjecture, 
Kawarabayashi, Ma and Yu planned  to find a contractible cycle in a 5-connected nonplanar graph containing no $K_4^-$ or 
$K_{2,3}$, and then use such a cycle to find a $TK_5$ by applying
augmenting path arguments. This plan (if successful), combined with the results in \cite{MY13, KMY15}, would give an alternative (and cleaner) solution to the 
Kelmans-Seymour conjecture. 

\medskip

One of the motivations for us to work on the Kelmans-Seymour conjecture 
was the following conjecture of Haj\'{o}s  (see e.g., \cite{Th05}) which, if true, would generalize the Four Color Theorem. 
\begin{conj}\label{hajos-5} 
Graphs containing no $TK_5$ are 4-colorable.
\end{conj}
It is known that Conjecture~\ref{hajos-5} holds for graphs with large girth (see K\"{u}hn and Osthus \cite{KO02}). Let $G$ be a possible
counterexample to Conjecture~\ref{hajos-5} with $|V(G)|$ minimum. Then our result on the Kelmans-Seymour conjecture 
implies that $G$  has connectivity at most 4.  
By a standard coloring argument, it is easy to show that $G$ must be 3-connected. 
It is shown in \cite{YZ04} that $G$ must be 4-connected. It is further
shown in \cite{SY15} that for every 4-cut $T$ of $G$, $G-T$ has exactly two components. 
The work in \cite{YZ04,SY15} suggests that  $G$ should be ``close'' to being 5-connected.

Haj\'{o}s actually made a more general conjecture in the 1950s: For
any positive integer $k$, every graph containing no $TK_{k+1}$ is $k$-colorable. This is easy to
verify for $k\le 3$ (see \cite{Di52}), 
and disproved in \cite{Ca79} for $k\ge 6$.  
However, it remains open for $k=4$ (Conjecture~\ref{hajos-5}) and $k=5$.  Thomassen \cite{Th05} pointed out connections 
between Haj\'{o}s' conjecture and Ramsey numbers, maximum cuts, and perfect graphs.  We refer the reader to \cite{Th05}
for other work and references related to Haj\'{o}s' conjecture and topological minors. 

In fact,  Erd\H{o}s and Fajtlowicz \cite{EF81} showed that the above general Haj\'{o}s' conjecture for $k\ge 6$ fails for almost all graphs.
Let $H(n):=\max\{\chi(G)/\sigma(G): G \mbox{ is a graph with } |V(G)|=n\}$, where $\chi(G)$ denotes the chromatic number of $G$ and 
$\sigma(G)$ denotes the largest $t$ such that $G$ contains $TK_t$.
Erd\H{o}s and Fajtlowicz \cite{EF81} showed that $H(n)=\Omega(\sqrt{n}/\log n)$, and conjectured that $H(n)=\Theta(\sqrt{n}/\log n)$. 
This conjecture was verified by Fox, Lee and Sudakov 
\cite{FLS12}, by studying $\sigma(G)$ in terms of independence number $\alpha(G)$. The following conjecture of Fox, Lee and Sudakov \cite{FLS12} 
is interesting. 
\begin{conj}\label{fls}
There is a constant $c>0$ such that every graph $G$ with $\chi(G)=k$ satisfies $\sigma(G)\ge c\sqrt{k\log k}$. 
\end{conj}

A key idea in \cite{MY10, MY13, HWY15I, HWY15II, HWY15III}  for finding $TK_5$ in graphs containing $K_4^-$ 
is to find a nonseparating path in a graph that avoids two given vertices. 
Let $G$ be a 5-connected nonplanar graph and $x_1,x_2,y_1,y_2\in V(G)$ such that $\{x_1,x_2,y_1,y_2\}$ induces a $K_4^-$ 
in which $x_1,x_2$ are of degree 3. We used an induced path $X$ in $G$ between $x_1$ and $x_2$ such that $G-X$ is 2-connected and 
$\{y_1,y_2\}\not\subseteq V(X)$, and in certain cases we need $X$ to
contain a special edge at $x_1$ (for example, in Section 6, $x_1=x$ is the special vertex 
representing the contraction of $M$). If we could find such $X$ that $G-X$ is 3-connected then our proofs would have been much simpler.  
This is related to the following conjecture of Lov\'{a}sz \cite{Lo75}. 

\begin{conj}\label{lovasz}
There exists an integer valued function $f(k)$ such that for any $f(k)$-connected graph $G$ and 
for any $A\subseteq V(G)$ with $|A|=2$, there exist vertex disjoint subgraphs
$G_1,G_2$ of $G$ such that $V(G_1)\cup V(G_2)=V(G)$, 
$G_1$ is a path between the vertices in $A$, and $G_2$ is $k$-connected. 
\end{conj}
A classical result of Tutte \cite{Tu63} implies $f(1)=3$. That
$f(2)=5$ was proved by Kriesell \cite{Kr01} and, independently, by
Chen, Gould and Yu \cite{CGY03}. Despite much effort from the research
community, Conjecture~\ref{lovasz}
remains open for $k\ge 3$. Variations of Conjecture~\ref{lovasz} for $k=2$ are
used in \cite{MY10, MY13, HWY15I, HWY15II,
  HWY15III} to resolve the Kelmans-Seymour conjecture. 
An edge version of Conjecture~\ref{lovasz} was conjectured
by Kriesell and proved by Kawarabayashi {\it et al.}
\cite{KLRW08}. Thomassen \cite{Th83b} conjectured a statement that is
more general than Conjecture~\ref{lovasz}  by allowing $|A|\ge 2$ and requiring $A\subseteq V(G_1)$ and 
$G_1$ be $k$-connected.



\end{document}